\documentclass[11pt]{article}

\def\sqr#1#2{{\vcenter{\vbox{\hrule height.#2pt
              \hbox{\vrule width.#2pt height#1pt \kern#1pt \vrule width.#2pt}
              \hrule height.#2pt}}}}
\def\signed #1{{\unskip\nobreak\hfil\penalty50
              \hskip2em\hbox{}\nobreak\hfil#1
              \parfillskip=0pt \finalhyphendemerits=0 \par}}
\def\endpf{\signed {$\sqr69$}}

\def\dbR{{\mathop{\rm l\negthinspace R}}}

\def\3n{\negthinspace \negthinspace \negthinspace }
\def\2n{\negthinspace \negthinspace }
\def\1n{\negthinspace }

\def\dbE{{\mathop{\rm l\negthinspace E}}}
\def\dbF{{\mathop{\rm l\negthinspace F}}}

\def\dbH{{\mathop{\rm l\negthinspace H}}}

\def\dbM{{\mathop{\rm l\negthinspace M}}}

\def\dbP{{\mathop{\rm l\negthinspace P}}}
\def\dbR{{\mathop{\rm l\negthinspace R}}}
\def\={\buildrel \triangle \over =}

\def\ds{\displaystyle}

\def\ns{\noalign{\ss}}

\def\b{\beta}
\def\g{\gamma}
\def\d{\delta}
\def\e{\varepsilon}
\def\z{\zeta}

\def\l{\lambda}
\def\m{\mu}

\def\si{\sigma}
\def\t{\tau}
\def\f{\varphi}
\def\th{\theta}
\def\o{\omega}

\def\G{\Gamma}
\def\D{\Delta}
\def\Th{\Theta}

\def\O{\Omega}

\def\cA{{\cal A}}
\def\cB{{\cal B}}

\def\cF{{\cal F}}
\def\cG{{\cal G}}
\def\cH{{\cal H}}

\def\cM{{\cal M}}

\def\cU{{\cal U}}

\def\cl{{\cal l}}
\def\ss{\smallskip}
\def\ms{\medskip}

\def\q{\quad}
\def\qq{\qquad}
\def\hb{\hbox}

\def\lan{\mathop{\langle}}
\def\ran{\mathop{\rangle}}

\def\esssup{\mathop{\rm esssup}}

\def\h{\widehat}
\def\wt{\widetilde}
\def\cd{\cdot}
\def\cds{\cdots}

\def\ae{\hbox{\rm a.e.{ }}}
\def\as{\hbox{\rm a.s.{ }}}

\def\cl{\overline}

\def\({\Big (}
\def\){\Big )}
\def\[{\Big[}
\def\]{\Big]}
\def\bde{\begin{definition}}
\def\ede{\end{definition}}
\def\be{\begin{equation}}
\def\bel{\begin{equation}\label}
\def\ee{\end{equation}}
\def\bt{\begin{theorem}}
\def\et{\end{theorem}}
\def\bc{\begin{corollary}}
\def\ec{\end{corollary}}
\def\bl{\begin{lemma}}
\def\el{\end{lemma}}
\def\bp{\begin{proposition}}
\def\ep{\end{proposition}}
\def\bas{\begin{assumption}}
\def\eas{\end{assumption}}
\def\br{\begin{remark}}
\def\er{\end{remark}}
\def\ba{\begin{array}}
\def\ea{\end{array}}
\def\ed{\end{document}}

\def\square#1{\vbox{\hrule\hbox{\vrule height#1%
     \kern#1\vrule}\hrule}}
\def\rectangle#1#2{\vbox{\hrule\hbox{\vrule height#1%
     \kern#2\vrule}\hrule}}

\font\tenbb=msbm10 \font\sevenbb=msbm7 \font\fivebb=msbm5

\newfam\bbfam
\scriptscriptfont\bbfam=\fivebb \textfont\bbfam=\tenbb
\scriptfont\bbfam=\sevenbb

\newtheorem{lemma}{Lemma}[section]
\newtheorem{remark}{Remark}[section]

\newtheorem{theorem}{Theorem}[section]
\newtheorem{corollary}{Corollary}[section]

\newtheorem{definition}{Definition}[section]
\newtheorem{proposition}{Proposition}[section]
\newtheorem{assumption}{Assumption}[section]

\makeatletter
   
   \@addtoreset{equation}{section}
\makeatother

\begin{document}
\title{\bf Mean-Field Backward Stochastic\\
 Volterra Integral Equations\footnote{This work is supported
in part by National Natural Science Foundation of China (Grants
10771122 and 11071145), Natural Science Foundation of Shandong Province of China (Grant Y2006A08),
Foundation for Innovative Research Groups of National Natural Science Foundation of China (Grant 10921101),
National Basic Research Program of China (973 Program, No.
2007CB814900), Independent Innovation Foundation of Shandong
University (Grant 2010JQ010), Graduate Independent Innovation
Foundation of Shandong University (GIIFSDU), and the NSF grant
DMS-1007514.}}

\author{Yufeng Shi\footnote{School of Mathematics, Shandong University, Jinan 250100,
China},~~Tianxiao Wang\footnote{School of Mathematics, Shandong
University, Jinan 250100, China},~~and~~Jiongmin
Yong\footnote{Department of Mathematics, University of Central
Florida, Orlando, FL 32816, USA.}}

\maketitle

\begin{abstract}

Mean-field backward stochastic Volterra integral equations
(MF-BSVIEs, for short) are introduced and studied. Well-posedness of
MF-BSVIEs in the sense of introduced adapted M-solutions is
established. Two duality principles between linear mean-field
(forward) stochastic Volterra integral equations (MF-FSVIEs, for
short) and MF-BSVIEs are obtained. Several comparison theorems for
MF-FSVIEs and MF-BSVIEs are proved. A Pontryagin's type maximum
principle is established for an optimal control of MF-FSVIEs.

\end{abstract}

\ms

\bf Keywords. \rm Mean-field stochastic Volterra integral equation,
mean-field backward stochastic Volterra integral equation, duality
principle, comparison theorem, maximum principle.

\ms

\bf AMS Mathematics subject classification. \rm 60H20, 93E20, 35Q83.

\ms

\section{Introduction.}\label{1}

Throughout this paper, we let $(\O,\cF,\dbF,\dbP)$ be a complete
filtered probability space on which a one-dimensional standard
Brownian motion $W(\cd)$ is defined with $\dbF=\{\cF_t\}_{t\ge0}$
being its natural filtration augmented by all the $\dbP$-null sets.
Let us begin with the following stochastic differential equation
(SDE, for short) in $\dbR$:
\bel{1.1}\left\{\ba{ll}
\ns\ds dX(t)=b(X(t),\m(t))dt+dW(t),\qq t\in[0,T],\\
\ns\ds X(0)=x,\ea\right.\ee
where
\bel{b}\ba{ll}
\ns\ds b(X(t),\m(t))=\int_\O
b(X(t,\o),X(t;\o'))\dbP(d\o')\\
\ns\ds\qq\qq\q~\equiv\int_\dbR
b(\xi,y)\m(t;dy)\Big|_{\xi=X(t)}\equiv\dbE[b(\xi,X(t))]
\Big|_{\xi=X(t)},\ea\ee
where $b:\dbR\times\dbR\to\dbR$ is a (locally) bounded Borel
measurable function and $\m(t;\cd)$ is the probability distribution
of the unknown process $X(t)$:
\bel{m}\m(t;A)=\dbP(X(t)\in A),\qq\forall A\in\cB(\dbR).\ee
Here $\cB(\dbR^n)$ is the Borel $\si$-field of $\dbR^n$ ($n\ge1$).
Equation (\ref{1.1}) is called a {\it McKean--Vlasov} SDE. Such an
equation was suggested by Kac \cite{Kac 1956} as a stochastic toy
model for the Vlasov kinetic equation of plasma and the study of
which was initiated by McKean \cite{McKean 1966}. Since then, many
authors made contributions on McKean--Vlasov type SDEs and
applications, see, for examples, Dawson \cite{Dawson 1983},
Dawson--G\"artner \cite{Dawson-Gartner 1987}, G\'artner
\cite{Gartner 1988}, Scheutzow \cite{Scheutzow 1987}, Sznitman
\cite{Sznitman 1989}, Graham \cite{Graham 1992}, Chan \cite{Chan
1994}, Chiang \cite{Chiang 1994}, Ahmed--Ding \cite{Ahmed-Ding
1995}. In recent years, related topics and problems have attracted
more and more attentions, see, for examples, Veretennikov
\cite{Veretennikov 2003}, Huang--Malham\'e--Caines
\cite{Huang-Malhame-Caines 2006}, Ahmed \cite{Ahmed 2007},
Mahmudov--McKibben \cite{Mahmudov-Mckibben 2007}, Lasry--Lions
\cite{Lasry-Lions 2007}, Borkar--Kumar \cite{Borkar-Kumar 2010},
Crisan--Xiong \cite{Crisan-Xiong 2010}, Kotelenez--Kurtz
\cite{Kotelenez-Kurtz 2010}, Park--Balasubramaniam--Kang
\cite{Park-Balasubramaniam-Kang 2008}, Andersson--Djehiche
\cite{Andersson-Djehice 2011}, Meyer-Brandis--Oksendal--Zhou
\cite{Meyer-Brandis-Oksandal-Zhou 2011}, and so on.

\ms

Inspired by (\ref{1.1}), one can consider the following more general
SDE:
\bel{MF-SDE}\left\{\ba{ll}
\ns\ds dX(t)=b(t,X(t),\dbE[\th^b(t,\xi,X(t))]_{\xi=X(t)})dt\\
\ns\ds\qq\qq\qq\qq+\si
(t,X(t),\dbE[\th^\si(t,\xi,X(t))]_{\xi=X(t)})dW(t),\qq t\in[0,T],\\
\ns\ds X(0)=x.\ea\right.\ee
where $\th^b$ and $\th^\si$ are some suitable maps. We call the
above a {\it mean-field} (forward) stochastic differential equation
(MF-FSDE, for short). From (\ref{b}) and (\ref{MF-SDE}), we see that
(\ref{1.1}) is a special case of (\ref{MF-SDE}). Note also that
(\ref{MF-SDE}) is an extension of classical It\^o type SDEs. Due to
the dependence of $b$ and $\si$ on
$\dbE[\th^b(t,\xi,X(t))]_{\xi=X(t)}$ and
$\dbE[\th^\si(t,\xi,X(t))]_{\xi=X(t)}$, respectively, MF-FSDE
(\ref{MF-SDE}) is {\it nonlocal} with respect to the event
$\o\in\O$.

\ms

It is easy to see that the equivalent integral form of
(\ref{MF-SDE}) is as follows:
\bel{1.3}\ba{ll}
\ns\ds X(t)=x+\int_0^tb(s,X(s),\dbE[\th^b(s,\xi,X(s))]_{\xi=X(s)})ds\\
\ns\ds\qq\qq\qq+\int_0^t\si(s,X(s),\dbE[\th^\si(s,\xi,X(s))]_{\xi
=X(s)})dW(s),\qq t\in[0,T].\ea\ee
This suggests a natural extension of the above to the following:
\bel{MF-FSVIE}\ba{ll}
\ns\ds X(t)=\f(t)+\int_0^tb(t,s,X(s),\dbE[\th^b(t,s,\xi,X(s))]_{\xi=X(s)})ds\\
\ns\ds\qq\qq\qq+\int_0^t\si(t,s,X(s),\dbE[\th^\si(t,s,\xi,X(s))]_{\xi
=X(s)})dW(s),\qq t\ge0.\ea\ee
We call the above a mean-field (forward) stochastic Volterra
integral equation (MF-FSVIE, for short). It is worthy of pointing
out that when the drift $b$ and diffusion $\si$ in (\ref{MF-FSVIE})
are independent of the nonlocal terms
$\dbE[\th^b(t,s,\xi,X(s))]_{\xi=X(s)}$ and
$\dbE[\th^\si(t,s,\xi,X(s))]_{\xi=X(s)}$, respectively,
(\ref{MF-FSVIE}) is reduced to a so-called (forward) stochastic
Volterra integral equation (FSVIEs, for short):
\bel{FSVIE}\ba{ll}
\ns\ds
X(t)=\f(t)+\int_0^tb(t,s,X(s))ds+\int_0^t\si(t,s,X(s))dW(s),\qq
t\ge0.\ea\ee
Such kind of equations have been studied by a number of researchers,
see, for examples, Berger--Mizel \cite{Berger-Mizel 1980}, Protter
\cite{Protter 1985}, Pardoux--Protter \cite{Pardoux-Protter 1990},
Tudor \cite{Tudor 1989}, Zhang \cite{Zhang 2010}, and so on.
Needless to say, the theory for (\ref{MF-FSVIE}) is very rich and
have a great application potential in various areas.

\ms

On the other hand, a general (nonlinear) backward stochastic
differential equation (BSDE, for short) introduced in Pardoux--Peng
\cite{Pardoux-Peng 1990} is equivalent to the following:
\bel{}Y(t)=\xi+\int_t^Tg(s,Y(s),Z(s))ds-\int_t^TZ(s)dW(s),\qq
t\in[0,T].\ee
Extending the above, the following general stochastic integral
equation was introduced and studied in Yong \cite{Yong 2006, Yong
2007, Yong 2008}:
\bel{BSVIE}Y(t)=\psi(t)+\int_t^Tg(t,s,Y(s),Z(t,s),Z(s,t))ds-\int_t^TZ(t,s)dW(s),\qq
t\in[0,T].\ee
Such an equation is called a backward stochastic Volterra integral
equation (BSVIE, for short). A special case of (\ref{BSVIE}) with
$g(\cd)$ independent of $Z(s,t)$ and $\psi(t)\equiv\xi$ was studied
by Lin \cite{Lin 2002} and Aman--N'zi \cite{Aman-N'Zi 2005} a little
earlier. Some relevant studies of (\ref{BSVIE}) can be found in
Wang--Zhang \cite{Wang-Zhang 2007}, Wang--Shi \cite{Wang-Shi 2010},
Ren \cite{Ren 2010}, and Anh--Grecksch--Yong \cite{Anh-Grecksch-Yong
2011}. Inspired by BSVIEs, it is very natural for us to introduce
the following stochastic integral equation:
\bel{MF-BSVIE}\ba{ll}
\ns\ds Y(t)=\psi(t)+\int_t^Tg(t,s,Y(s),Z(t,s),Z(s,t),\Gamma
(t,s,Y(s),Z(t,s),Z(s,t)))ds\\
\ns\ds\qq\qq\qq-\int_t^TZ(t,s)dW(s),\qq t\in[0,T],\ea\ee
where $(Y(\cd),Z(\cd\,,\cd))$ is the pair of unknown processes,
$\psi(\cd)$ is a given {\it free term} which is $\cF_T$-measurable
(not necessarily $\dbF$-adapted), $g(\cd)$ is a given mapping, called
the {\it generator}, and
\bel{G}\G(t,s,Y,Z,\h Z)=\dbE\[\th(t,s,y,z,\hat z,Y,Z,\h
Z)\]_{(y,z,\hat z)=(Y,Z,\h Z)}\ee
with $(Y,Z,\h Z)$ being some random variables, for some mapping
$\th(\cd)$ (see the next section for precise meaning of the above).
We call (\ref{MF-BSVIE}) a {\it mean-field backward stochastic
Volterra integral equation} (MF-BSVIE, for short). Relevant to the
current paper, let us mention that in Buckdahn--Djehiche--Li--Peng
\cite{Buckdahn-Djehiche-Li-Peng 2009}, mean-field backward
stochastic differential equations (MF-BSDEs, for short) were
introduced and in Buckdahn--Li--Peng \cite{Buckdahn-Li-Peng 2009} a
class of nonlocal PDEs are studied with the help of an MF-BSDE and a
McKean-Vlasov forward equation.

\ms

We see that MF-BSVIE (\ref{MF-BSVIE}) not only includes MF-BSDEs
(which, of course, also includes standard BSDEs) introduced in
\cite{Buckdahn-Djehiche-Li-Peng 2009,Buckdahn-Li-Peng 2009}, but
also generalizes BSVIEs studied in \cite{Yong 2006, Yong 2008,
Wang-Shi 2010}, etc. in a natural way. Besides, investigating
MF-BSVIEs allows us to meet the need in the study of optimal control
for MF-FSVIEs. As a matter of fact, in the statement of Pontryagin
type maximum principle for optimal control of a forward
(deterministic or stochastic) control system, the adjoint equation
of variational state equation is a corresponding (deterministic or
stochastic) backward system, see \cite{Yong-Zhou 1999} for the case
of classical optimal control problems, \cite{Andersson-Djehice 2011,
Buckdahn-Djehiche-Li 2011, Meyer-Brandis-Oksandal-Zhou 2011} for the
case of MF-FSDEs, and \cite{Yong 2006, Yong 2008} for the case of
FSVIEs. When the state equation is an MF-FSVIE, the adjoint equation
will naturally be an MF-BSVIE. Hence the study of well-posedness for
MF-BSVIEs is not avoidable when we want to study optimal control
problems for MF-BSVIEs.

\ms

The novelty of this paper mainly contains the following: First,
well-posedness of general MF-BSVIEs will be established. In doing
that, we discover that the growth of the generator and the nonlocal
term with respect to $Z(s,t)$ plays a crucial role; a better
understanding of which enables us to have found a neat way of
treating term $Z(s,t)$. Even for BSVIEs, our new method will
significantly simplify the proof of well-posedness of the equation
(comparing with \cite{Yong 2008}). Second, we establish two slightly
different duality principles, one starts from linear MF-FSVIEs, and
the other starts from linear MF-BSVIEs. We found that ``{\sl Twice
adjoint of a linear MF-FSVIE is itself}'', whereas, ``{\sl Twice
adjoint of a linear MF-BSVIE is not necessarily itself}''. Third,
some comparison theorems will be established for MF-FSVIEs and
MF-BSVIEs. It turns out that the situation is surprisingly different
from the differential equation cases. Some mistakes found in
\cite{Yong 2006, Yong 2007} will be corrected. Finally, as an
application of the duality principle for MF-FSVIEs, we establish a
Pontryagin type maximum principle for an optimal control problem of
MF-FSVIEs.

\ms

The rest of the paper is organized as follows. Section 2 is devoted
to present some preliminary results. In Section 3, we prove the
existence and uniqueness of adapted M-solutions to MF-BSVIE
(\ref{MF-BSVIE}). In Section 4 we obtain duality principles.
Comparison theorems will be presented in Section 5. In Section 6, we
deduce a maximum principle of optimal controls for MF-FSVIEs.

\section{Preliminary Results.}

In this section, we will make some preliminaries.

\subsection{Formulation of MF-BSVIEs.}

Let us first introduce some spaces. For $H=\dbR^n$, etc., and $p>1$,
$t\in[0,T]$, let
$$\ba{ll}
\ns\ds L^p(0,T;H)=\Big\{x:[0,T]\to
H\Bigm|\int_0^T|x(s)|^pds<\infty\Big\},\\
\ns\ds L_{\cF_t}^p(\O;H)=\Big\{\xi:\O\to H\Bigm|\xi\hb{ is
$\cF_t$-measurable, }\dbE|\xi|^p<\infty\Big\},\\
\ns\ds L_{\cF_t}^p(0,T;H)=\Big\{X:[0,T]\times\O\to H\Bigm|X(\cd)\hb{
is
$\cF_t$-measurable, }\dbE\int_0^T|X(s)|^pds<\infty\Big\},\\
\ns\ds L_\dbF^p(0,T;H)=\Big\{X:[0,T]\times\O\to H\Bigm|X(\cd)\hb{ is
$\dbF$-adapted, }\dbE\int_0^T|X(s)|^pds<\infty\Big\},\\
\ns\ds L^p_{\dbF}(\Omega;L^2(0,T;H))=\Big\{X:[0,T]\times\O\to
H\Bigm| X(\cdot)\hb{ is $\dbF$-adapted, }
\dbE\(\int_0^T|X(s)|^2ds\)^{p\over2}<\infty\Big\}.\ea$$
Also, let (with $q\ge1$)
$$\ba{ll}
\ns\ds L^p(0,T;L^q_\dbF(0,T;H))=\Big\{Z:[0,T]^2\times\O\to
H\Bigm|Z(t,\cd)\hb{ is $\dbF$-adapted for almost all $t\in[0,T]$,
}\\
\ns\ds\qq\qq\qq\qq\qq\qq\dbE\int_0^T\(\int_0^T|Z(t,s)|^qds\)^{p\over q}dt<\infty\Big\},\\
\ns\ds C_\dbF^p([0,T];H)=\Big\{X:[0,T]\times\O\to H\Bigm|X(\cd)\hb{
is $\dbF$-adapted, $t\mapsto X(t)$ is continuous}\\
\ns\ds\qq\qq\qq\qq\qq\qq\hb{ from $[0,T]$ to $L^p_{\cF_T}(\O;H)$,
}\sup_{t\in[0,T]}\dbE[|X(t)|^p]<\infty\Big\}.\ea$$
We denote
$$\ba{ll}
\ns\ds\cH^p[0,T]=L^p_\dbF(0,T;H)\times
L^p(0,T;L_\dbF^2(0,T;H)),\\
\ns\ds\dbH^p[0,T]=C^p_\dbF([0,T];H)\times
L^p_{\dbF}(\Omega;L^2(0,T;H)).\ea$$
Next, let $(\O^2,\cF^2,\dbP^2)=(\O\times\O,\cF\otimes
\cF,\dbP\otimes\dbP)$ be the completion of the product probability
space of the original $(\O,\cF,\dbP)$ with itself, where we define
the filtration as $\dbF^2=\{\cF_t\otimes\cF_t,~t\in[0,T]\}$ with
$\cF_t\otimes\cF_t$ being the completion of $\cF_t\times\cF_t$. It
is worthy of noting that any random variable $\xi=\xi(\o)$ defined
on $\O$ can be extended naturally to $\O^2$ as
$\xi'(\o,\o')=\xi(\o)$, with $(\o,\o')\in\O^2$. Similar to the
above, we define
$$L^1(\O^2,\cF^2,\dbP^2;H)=\Big\{\xi:\O^2\1n\to\1n H\Bigm|\xi\hb{ is
$\cF^2$-measurable, }
\dbE^2|\xi|\1n\equiv\2n\int_{\O^2}|\xi(\o',\o)|\dbP(d\o')\dbP(d\o)<\infty\Big\}.$$
For any $\eta\in L^1(\O^2,\cF^2,\dbP^2;H)$, we denote
$$\dbE'\eta(\o,\cd)=\int_\O\eta(\o,\o')\dbP(d\o')\in
L^1(\O,\cF,\dbP).$$
Note that if $\eta(\o,\o')=\eta(\o')$, then
$$\dbE'\eta=\int_\O\eta(\o')\dbP(d\o')=\int_\O\eta(\o)\dbP(d\o)=\dbE\eta.$$
In what follows, $\dbE'$ will be used when we need to distinguish
$\o'$ from $\o$, which is the case when both $\o$ and $\o'$ appear
at the same time. Finally, we denote
$$\D=\Big\{(t,s)\in[0,T]^2\Bigm|t\le s\Big\},\q\D^*=\Big\{(t,s)\in[0,T]^2\Bigm|t\ge s\Big\}\equiv
\cl{\D^c}.$$
Let
\bel{g,f}
g:\D\times\O\times\dbR^{3n}\times\dbR^m\to\dbR^n,\qq\th:\D\times\O^2\times\dbR^{6n}\to\dbR^m,\ee
be some suitable maps (see below for precise conditions) and define
\bel{G2}\ba{ll}
\ns\ds\G(t,s,Y,Z,\h Z)=\dbE'\[\th(t,s,y,z,\h z,Y,Z,\h
Z)\]_{(y,z,\hat
z)=(Y,Z,\h Z)}\\
\ns\ds=\int_\O\th(t,s,\o,\o',Y(\o),Z(\o),\h Z(\o),Y(\o'),Z(\o'),\h
Z(\o'))\dbP(d\o'),\ea\ee
for all reasonable random variables $(Y,Z,\h Z)$. This gives the
precise meaning of (\ref{G}). Hereafter, when we talk about MF-BSVIE
(\ref{MF-BSVIE}), the mapping $\G$ is defined by (\ref{G2}). With
such a mapping, we have
$$\ba{ll}
\ns\ds\G(t,s,Y(s),Z(t,s),Z(s,t))\equiv\G(t,s,\o,Y(s,\o),Z(t,s,\o),Z(s,t,\o))\\
\ns\ds=\int_\O\th(t,s,\o,\o',Y(s,\o),Z(t,s,\o),Z(s,t,\o),Y(s,\o'),
Z(t,s,\o'),Z(s,t,\o'))\dbP(d\o').\ea$$
Clearly, the operator $\G$ is {\it nonlocal} in the sense that the
value $\G(t,s,\o,Y(s,\o),Z(t,s,\o),Z(s,t,\o))$ of
$\G(t,s,Y(s),Z(t,s),Z(s,t))$ at $\o$ depends on the whole set
$$\{(Y(s,\o'),Z(t,s,\o'),Z(s,t,\o'))\bigm|\o'\in\O\},$$
not just on $(Y(s,\o),Z(t,s,\o),Z(s,t,\o))$. To get some feeling
about such an operator, let us look at a simple but nontrivial
special case.

\ms

\bf Example 2.1. \rm Let
$$\ba{ll}
\ns\ds\th(t,s,\o,\o',y,z,\hat z,y',z',\hat
z')=\th_0(t,s,\o)+A_0(t,s,\o)y+B_0(t,s,\o)z+C_0(t,s,\o)\hat
z\\
\ns\ds\qq\qq\qq\qq\qq\qq\q+A_1(t,s,\o,\o')y'+B_1(t,s,\o,\o')z'+C_1(t,s,\o,\o')\hat
z'.\ea$$
We should carefully distinguish $\o'$ and $\o$ in the above. Then
(suppressing $\o$)
$$\ba{ll}
\ns\ds\G(t,s,Y(s),Z(t,s),Z(s,t))=\th_0(t,s)+A_0(t,s)Y(s)+B_0(t,s)Z(t,s)
+C_0(t,s)Z(s,t)\\
\ns\ds\qq\qq\qq\qq\qq\qq\q+\dbE'[A_1(t,s)Y(s)]+\dbE'[B_1(t,s)Z(t,s)]+\dbE'[C_1(t,s)Z(s,t)],\ea$$
where, for example,
$$\dbE'[B_1(t,s)Z(t,s)]=\int_\O B_1(t,s,\o,\o')Z(t,s,\o')\dbP(d\o').$$
For such a case,
$(Y(\cd),Z(\cd\,,\cd))\mapsto\G(\cd\,,\cd\,,Y(\cd),Z(\cd\,,\cd),Z(\cd\,,\cd))$
is affine.

\ms

Having some feeling about the operator $\G$ from the above, let us
look at some useful properties of the operator $\G$ in general. To
this end, we make the following assumption.

\ms

{\bf(H0)$_q$} The map $\th:\D\times\O^2\times\dbR^{6n}\to\dbR^m$ is
measurable and for all $(t,y,z,\hat z,y',z',\hat
z')\in[0,T]\times\dbR^{6n}$, the map
$(s,\o,\o')\mapsto\th(t,s,\o,\o',y,z,\hat z,y',z',\hat z')$ is
$\dbF^2$-progressively measurable on $[t,T]$. Moreover, there exist
constants $L>0$ and $q\ge2$ such that
\bel{th-Lip}\ba{ll}
\ns\ds|\th(t,s,\o,\o',y_1,z_1,\hat z_1,y_1',z_1',\hat
z_1')-\th(t,s,\o,\o',y_2,z_2,\hat z_2,y_2',z_2',\hat z_2')|\\
\ns\ds\le L\(|y_1-y_2|+|z_1-z_2|+|\hat z_1-\hat
z_2|+|y_1'-y_2'|+|z_1'-z_2'|+|\hat z_1'-\hat z_2'|\),\\
\ns\ds\qq\qq\qq\forall(t,s,\o,\o')\in\D\times\O^2,~(y_i,z_i,\hat
z_i,y_i',z_i',\hat z_i')\in\dbR^{6n},i=1,2,\ea\ee
and
\bel{th-growth}\ba{ll}
\ns\ds|\th(t,s,\o,\o',y,z,\hat z,y',z',\hat z')|\le
L\(1+|y|+|z|+|\hat z|^{2\over q}+|y'|
+|z'|+|\hat z'|^{2\over q}\),\\
\ns\ds\qq\qq\qq\forall(t,s,\o',\o)\in\D\times\O^2,~(y,z,\hat
z,y',z',\hat z')\in\dbR^{6n}.\ea\ee

In the above, we may replace constant $L$ by some function $L(t,s)$
with certain integrability (similar to \cite{Yong 2008}). However,
for the simplicity of presentation, we prefer to take a constant
$L$. Also, we note that $(\hat z,\hat
z')\mapsto\th(t,s,\o,\o',y,z,\hat z,y',z',\hat z')$ is assumed to
grow no more than $|\hat z|^{2\over q}+|\hat z'|^{2\over q}$. If
$q=2$, then the growth is linear and if $q>2$, the growth is
sublinearly. This condition is very subtle in showing that the
solution $(Y(\cd),Z(\cd\,,\cd))$ of MF-BSVIE belongs in
$\cH^q[0,T]$. We would like to mention that (H0)$_\infty$ is
understood as that (\ref{th-growth}) is replaced by the following
\bel{th-growth2}\ba{ll}
\ns\ds|\th(t,s,\o,\o',y,z,\hat z,y',z',\hat z')|\le
L\(1+|y|+|z|+|y'|+|z'|\),\\
\ns\ds\qq\qq\qq\forall(t,s,\o',\o)\in\D\times\O^2,~(y,z,\hat
z,y',z',\hat z')\in\dbR^{6n}.\ea\ee

Under (H0)$_q$, for any $(Y(\cd),Z(\cd\,,\cd))\in\cH^q[0,T]$, we see
that for each $t\in[0,T]$, the map
$$\ba{ll}
\ns\ds(s,\o)\mapsto\G(t,s,\o,Y(s),Z(t,s),Z(s,t))\\
\ns\ds\equiv\int_\O\th(t,s,\o,\o',Y(s,\o),Z(t,s,\o),Z(s,t,\o),Y(s,\o'),Z(t,s,\o'),Z(s,t,\o'))\dbP(d\o')\ea$$
is $\dbF$-progressively measurable on $[t,T]$. Also,
\bel{2.6}\ba{ll}
\ns\ds|\th(t,s,\o,\o',Y(s,\o),Z(t,s,\o),Z(s,t,\o),y,z,\hat z)|\\
\ns\ds\le L\(1+|Y(s,\o)|+|Z(t,s,\o)|+|Z(s,t,\o)|^{2\over
q}+|y|+|z|+|\hat z|^{2\over q}\).\ea\ee
Consequently,
\bel{2.7}\ba{ll}
\ns\ds|\G(t,s,Y(s),Z(t,s),Z(s,t))|\\
\ns\ds\le L\(1+|Y(s)|+|Z(t,s)|+|Z(s,t)|^{2\over
q}+\dbE|Y(s)|+\dbE|Z(t,s)|+\dbE|Z(s,t)|^{2\over q}\).\ea\ee
Likewise, for any
$(Y_1(\cd),Z_1(\cd\,,\cd)),(Y_2(\cd),Z_2(\cd\,,\cd))\in\cH^q[0,T]$,
we have
\bel{2.8}\ba{ll}
\ns\ds|\G(t,s,Y_1(s),Z_1(t,s),Z_1(s,t))-\G(t,s,Y_2(s),Z_2(t,s),Z_2(s,t))|\\
\ns\ds\le
L\(|Y_1(s)-Y_2(s)|+|Z_1(t,s)-Z_2(t,s)|+|Z_1(s,t)-Z_2(s,t)|\\
\ns\ds\qq+\dbE|Y_1(s)-Y_2(s)|+\dbE|Z_1(t,s)-Z_2(t,s)|+\dbE|Z_1(s,t)-Z_2(s,t)|\).\ea\ee
The above two estimates will play an interesting role later. We now
introduce the following definition.

\ms

\bf Definition 2.2. \rm A pair of
$(Y(\cd),Z(\cd\,,\cd))\in\cH^p[0,T]$ is called an {\it adapted
M-solution} of MF-BSVIE (\ref{MF-BSVIE}) if (\ref{MF-BSVIE}) is
satisfied in the It\^o sense and the following holds:
\bel{M}Y(t)=\dbE Y(t)+\int_0^tZ(t,s)dW(s),\qq\qq t\in[0,T].\ee

\ms

It is clear that (\ref{M}) implies
\bel{M2}Y(t)=\dbE[Y(t)\bigm|\cF_S]+\int_S^tZ(t,s)dW(s),\qq0\le S\le
t\le T.\ee
This suggests us define $\cM^p[0,T]$ as the set of all elements
$(y(\cd),z(\cd\,,\cd))\in\cH^p[0,T]$ satisfying:
\bel{2.10}y(t)=\dbE\[y(t)\bigm|\cF_S\]+\int_S^tz(t,s)dW(s),\qq
t\in[S,T],\q S\in[0,T).\ee
Obviously $\cM^p[0,T]$ is a closed subspace of $\cH^p[0,T]$. Note
that for any $(y(\cd),z(\cd\,,\cd))\in\cM^2[0,T]$,
\bel{2.11}\dbE|y(t)|^2=\(\dbE\[y(t)\bigm|\cF_S\]\)^2+\dbE\int_S^t|z(t,s)|^2ds\ge\dbE\int_S^t|z(t,s)|^2ds.\ee
Relation (\ref{2.11}) can be generalized a little bit more. To see
this, let us present the following lemma.

\ms

\bf Lemma 2.3. \sl Let $0\le S<t\le T$, $\eta\in
L^p_{\cF_S}(\O;\dbR^n)$ and $\z(\cd)\in
L^p_\dbF(\O;L^2(S,t;\dbR^n))$. Then
\bel{2.12}\dbE\[|\eta|^p+\(\int_S^t|\z(s)|^2ds\)^{p\over2}\]\le
K\dbE\Big|\eta+\int_S^t\z(s)dW(s)\Big|^p.\ee
Hereafter, $K>0$ stands for a generic constant which can be
different from line to line.

\ms

\it Proof. \rm For fixed $(S,t)\in\D$ (which means $0\le S\le t\le
T$) with $S<t$, let
$$\xi=\eta+\int_S^t\z(s)dW(s),$$
which is $\cF_t$-measurable. Let $(Y(\cd),Z(\cd))$ be the adapted
solution to the following BSDE:
$$Y(r)=\xi-\int_r^tZ(s)dW(s),\qq r\in[S,t].$$
Then it is standard that
\bel{2.13}\dbE\[\sup_{r\in[S,t]}|Y(r)|^p+\(\int_S^t|Z(s)|^2ds\)^{p\over2}\]\le
K\dbE|\xi|^p.\ee
Now,
$$Y(S)+\int_S^tZ(s)dW(s)=\xi=\eta+\int_S^t\z(s)dW(s).$$
By taking conditional expectation $\dbE[\cd\,|\,\cF_S]$, we see that
$$Y(S)=\eta.$$
Consequently,
$$\int_S^t\(Z(s)-\z(s)\)dW(s)=0,$$
which leads to
$$Z(s)=\z(s),\qq s\in[S,t],~\as$$
Then (\ref{2.12}) follows from (\ref{2.13}). \endpf

\ms

We have the following interesting corollary for elements in
$\cM^p[0,T]$ (comparing with (\ref{2.11})).

\ms

\bf Corollary 2.4. \sl For any $(y(\cd),z(\cd\,,\cd))\in\cM^p[0,T]$,
the following holds:
\bel{2.14}\dbE\(\int_S^t|z(t,s)|^2ds\)^{p\over2}\le
K\dbE|y(t)|^p,\qq\forall S\in[0,t].\ee

\ms

\it Proof. \rm Applying (\ref{2.12}) to (\ref{2.10}), we have
$$\dbE\(\int_S^t|z(t,s)|^2ds\)^{p\over2}\le\dbE\[\big|\dbE[y(t)\bigm|\cF_S]\big|^p+
\(\int_S^t|z(t,s)|^2ds\)^{p\over2}\]\le K\dbE|y(t)|^p.$$
This proves the corollary. \endpf

\ms

From the above, we see that for any
$(y(\cd),z(\cd , \cd))\in\cM^p[0,T]$, and any $\b>0$,
\bel{}\ba{ll}
\ns\ds K\dbE\int_0^Te^{\b t}|y(t)|^pdt\ge\dbE\int_0^Te^{\b
t}\[|\dbE y(t)|^p+\(\int_0^t|z(t,s)|^2ds\)^{p\over2}\]dt\\
\ns\ds\qq\qq\qq\qq\ge\dbE\int_0^Te^{\b
t}\(\int_0^t|z(t,s)|^2ds\)^{p\over2}dt.\ea\ee
Hence,
$$\ba{ll}
\ns\ds\|(y(\cd),z(\cd\,,\cd))\|_{\cH^p[0,T]}^p\equiv\dbE\[\int_0^T|y(t)|^pdt+\int_0^T\(\int_0^T
|z(t,s)|^2ds\)^{p\over2}dt\]\\
\ns\ds\le
K\dbE\[\int_0^T|y(t)|^pdt+\int_0^T\(\int_0^t|z(t,s)|^2ds\)^{p\over2}dt
+\int_0^T\(\int_t^T|z(t,s)|^2ds\)^{p\over2}dt\]\\
\ns\ds\le K\dbE\[\int_0^Te^{\b t}|y(t)|^pdt+\int_0^Te^{\b
t}\(\int_0^t|z(t,s)|^2ds\)^{p\over2}dt+
\int_0^Te^{\b t}\(\int_t^T|z(t,s)|^2ds\)^{p\over2}dt\]\\
\ns\ds\le K\dbE\[\int_0^Te^{\b t}|y(t)|^pdt+\int_0^Te^{\b
t}\(\int_t^T|z(t,s)|^2ds\)^{p\over2}dt\]\le
K\|(y(\cd),z(\cd\,,\cd))\|_{\cH^p[0,T]}^p.\ea$$
This means that we can use the following as an equivalent norm in
$\cM^p[0,T]$:
$$\|(y(\cd),z(\cd\,,\cd))\|_{\cM^p[0,T]}\equiv\left\{\dbE\int_0^Te^{\b
t}|y(t)|^pdt+\dbE\int_0^Te^{\b
t}\(\int_t^T|z(t,s)|^2ds\)^{p\over2}dt\right\}^{1\over p}.$$
Sometimes we use $\cM^p_\b[0,T]$ for $\cM^p[0,T]$ to emphasize the
involved parameter $\b$.

\ms

To conclude this subsection, we state the following corollary of
Lemma 2.3 relevant to BSVIEs, whose proof is straightforward.

\ms

\bf Corollary 2.5. \sl Suppose $(\eta(\cd),\z(\cd ,\cd))$ is an
adapted M-solution to the following BSVIE:
\bel{2.15}\eta(t)=\xi(t)+\int_t^Tg(t,s)ds-\int_t^T\z(t,s)dW(s),\q
t\in[0,T],\ee
for $\xi(\cd)\in L^p_{\cF_T}(0,T;\dbR^n)$ and $g(\cd\,,\cd)\in
L^p(0,T;L^1_\dbF(0,T;\dbR^n))$. Then
\bel{2.16}\dbE\[|\eta(t)|^p+\(\int_t^T|\z(t,s)|^2ds\)^{p\over2}\]\le
K\dbE\[|\xi(t)|^p+\(\int_t^T|g(t,s)|ds\)^p\], \q \forall t\in[0,T].\ee

\rm

\ms

\subsection{Mean-field forward stochastic Volterra integral
equations.}

In this subsection, we study the following MF-FSVIE:
\bel{MF-FSVIE1}\ba{ll}
\ns\ds X(t)=\f(t)+\int_0^tb(t,s,X(s),\G^b(t,s,X(s)))ds\\
\ns\ds\qq\qq\qq+\int_0^t\si(t,s,X(s),\G^\si(t,s,X(s)))dW(s),\qq
t\in[0,T],\ea\ee
where
\bel{2.6}\left\{\ba{ll}
\ns\ds\G^b(t,s,X)=\dbE'\[\th^b(t,s,\xi,X,)\]_{\xi=X}\equiv\int_\O\th^b(t,s,\o,\o',X(\o),X(\o'))
\dbP(d\o'),\\
\ns\ds\G^\si(t,s,X)=\dbE'\[\th^\si(t,s,\xi,X)\]_{\xi=X}\equiv\int_\O\th^\si(t,s,\o,\o',X(\o),X(\o'))
\dbP(d\o'). \ea\right.\ee
We see that MF-FSVIE (\ref{MF-FSVIE1}) is slightly more general than
MF-FSVIE (\ref{MF-FSVIE}) because of the above definition
(\ref{2.6}) of the operators $\G^b$ and $\G^\si$.

\ms

An $\dbF$-adapted process $X(\cd)$ is called a solution to
(\ref{MF-FSVIE1}) if (\ref{MF-FSVIE1}) is satisfied in the usual
It\^o sense. To guarantee the well-posedness of (\ref{MF-FSVIE1}),
let us make the following hypotheses.

\ms
{\bf(H1)} The maps $b:
\D^*\times\O\times\dbR^n\times\dbR^{m_1}\to\dbR^n$ and
$\si:\D^*\times\O\times\dbR^n\times\dbR^{m_2}\to\dbR^n$ are
measurable, and for all
$(t,x,\g,\g')\in[0,T]\times\dbR^n\times\dbR^{m_1}\times\dbR^{m_2}$,
the map
$$(s,\o)\mapsto(b(t,s,\o,x,\g),\si(t,s,\o,x,\g'))$$
is $\dbF$-progressively measurable on $[0,t]$. Moreover, there
exists some constant $L>0$ such that
\bel{b-si-Lip}\ba{ll}
\ns\ds
|b(t,s,\o,x_1,\g_1)-b(t,s,\o,x_2,\g_2)|+|\si(t,s,\o,x_1,\g'_1)-\si(t,s,\o,x_2,\g'_2)|\\
\ns\ds\qq\le
L(|x_1-x_2|+|\g_1-\g_2|+|\g'_1-\g'_2|),\\
\ns\ds\qq\qq\qq(t,s,\o)\in\D^*\times\O,~(x_i,\g_i,\g_i')\in\dbR^n\times\dbR^{m_1}\times\dbR^{m_2},~
i=1,2.\ea\ee
Moreover,
\bel{b-si-growth}\ba{ll}
\ns\ds|b(t,s,\o,x,\g)|+|\si(t,s,\o,x,\g')|\le
L(1+|x|+|\g|+|\g'|),\\
\ns\ds\qq\qq\qq(t,s,\o,x,\g,\g')\in\D^*\times\O\times\dbR^n\times\dbR^{m_1}\times\dbR^{m_2}.\ea\ee

{\bf(H2)} The maps
$\th^b:\D^*\times\O^2\times\dbR^{2n}\to\dbR^{m_1}$ and
$\th^\si:\D^*\times\O^2\times\dbR^{2n}\to\dbR^{m_2}$ are measurable,
and for all $(t,x,x')\in[0,T]\times\dbR^{2n}$, the map
$$(s,\o,\o')\mapsto(\th^b(t,s,\o,\o',x,x'),\th^\si(t,s,\o,\o',x,x'))$$
is $\dbF^2$-progressively measurable on $[0,t]$. Moreover, there
exists some constant $L>0$ such that
\bel{th-b-Lip}\ba{ll}
\ns\ds
|\th^b(t,s,\o,\o',x_1,x'_1)-\th^b(t,s,\o,\o',x_2,x'_2)|+|\th^\si(t,s,\o,\o',x_1,x'_1)
-\th^\si(t,s,\o,\o',x_2,x'_2)|\\
\ns\ds\qq\le
L(|x_1-x_2|+|x'_1-x'_2|),\qq(t,s,\o,\o')\in\D^*\times\O^2,
(x_i,x_i')\in\dbR^{2n},i=1,2,\ea\ee
and
\bel{th-b-growth}\ba{ll}
\ns\ds|\th^b(t,s,\o,\o',x,x')|+|\th^\si(t,s,\o,\o',x,x')|\le
L(1+|x|+|x'|),\\
\ns\ds\qq\qq\qq\qq\qq\qq(t,s,\o,\o')\in\D^*\times\O^2,
x,x'\in\dbR^n.\ea\ee

We will also need the following assumptions.

\ms

{\bf(H1)$'$} In addition to (H1), the map
$t\mapsto(b(t,s,\o,x,\g),\si(t,s,\o,x,\g'))$ is continuous on
$[s,T]$.

\ms

{\bf(H2)$'$} In addition to (H2), the map
$t\mapsto(\th^b(t,s,\o,\o',x,x'),\th^\si(t,s,\o,\o',x,x'))$ is
continuous on $[s,T]$.

\ms

Now, let us state and prove the following result concerning MF-FSVIE
(\ref{MF-FSVIE1}).

\ms

\bf Theorem 2.6. \sl Let {\rm(H1)--(H2)} hold. Then for any $p\ge2$,
and $\f(\cd)\in L^p_\dbF(0,T;\dbR^n)$, MF-FSVIE $(\ref{MF-FSVIE1})$
admits a unique solution $X(\cd)\in L^p_\dbF(0,T;\dbR^n)$, and the
following estimate holds:
\bel{|X|Lp-estimate}\dbE\int_0^T|X(t)|^pdt\le
K\(1+\dbE\int_0^T|\f(t)|^pdt\).\ee
Further, for $i=1,2$, let $X_i(\cd)\in L^p_\dbF(0,T;\dbR^n)$ be the
solutions of $(\ref{MF-FSVIE1})$ corresponding to $\f_i(\cd)\in
L^p_\dbF(0,T;\dbR^n)$ and
$b_i(\cd),\si_i(\cd),\th_i^b(\cd),\th_i^\si(\cd)$ satisfying
{\rm(H1)--(H2)}. Let
$$\left\{\ba{ll}
\ns\ds\G_i^b(t,s,X)=\dbE'\[\th_i^b(t,s,\xi,X)\]_{\xi=X}\equiv\int_\O\th_i^b(t,s,\o,\o',X(\o),X(\o'))
\dbP(d\o'),\\
\ns\ds\G^\si_i(t,s,X)=\dbE'\[\th_i^\si(t,s,\xi,X)\]_{\xi=X}\equiv\int_\O\th_i^\si(t,s,\o,\o',X(\o),
X(\o'))\dbP(d\o'),\ea\right.\q i=1,2.$$
Then the following stability estimate
holds:
\bel{|X-X|Lp-estimate}\ba{ll}
\ns\ds\dbE\int_0^T|X_1(t)-X_2(t)|^p\le K\Big\{\dbE\int_0^T|\f_1(t)-\f_2(t)|^pdt\\
\ns\ds\q+\dbE\int_0^T\(\int_0^t|b_1(t,s,X_1(s),\G^b_1(t,s,X_1(s)))-b_2(t,s,X_1(s),\G^b_2(t,s,X_1(s)))|ds\)^pdt\\
\ns\ds\q+\dbE\int_0^T\(\int_0^t|\si_1(t,s,X_1(s),\G^b_1(t,s,X_1(s)))-\si_2(t,s,X_1(s),\G^b_2(t,s,X_1(s)))|^2
ds\)^{p\over2}dt\Big\}.\ea\ee
Moreover, let {\rm(H1)$'$--(H2)$'$} hold. Then for any $p\ge2$, and
any $\f(\cd)\in C^p_\dbF([0,T];\dbR^n)$, the unique solution
$X(\cd)\in C^p_\dbF([0,T];\dbR^n)$, and estimate
$(\ref{|X|Lp-estimate})$ is replaced by the following:
\bel{|X|-estimate}\sup_{t\in[0,T]}\dbE|X(t)|^p\le
K\Big\{1+\sup_{t\in[0,T]}\dbE|\f(t)|^p\Big\}.\ee
Also, for $i=1,2$, let $X_i(\cd)\in L^p_\dbF(0,T;\dbR^n)$ be the
solutions of $(\ref{MF-FSVIE1})$ corresponding to $\f_i(\cd)\in
L^p_\dbF(0,T;\dbR^n)$ and
$b_i(\cd),\si_i(\cd),\th_i^b(\cd),\th_i^\si(\cd)$ satisfying
{\rm(H1)$'$--(H2)$'$}. Then $(\ref{|X-X|Lp-estimate})$ is replaced
by the following:
\bel{|X-X|-estimate}\ba{ll}
\ns\ds\sup_{t\in[0,T]}\dbE|X_1(t)-X_2(t)|^p\le K\Big\{\sup_{t\in[0,T]}\dbE|\f_1(t)-\f_2(t)|^p\\
\ns\ds+\sup_{t\in[0,T]}\dbE\(\int_0^t|b_1(t,s,X_1(s),\G^b_1(t,s,X_1(s)))-b_2(t,s,X_1(s),\G^b_2(t,s,X_1(s)))|ds\)^p\\
\ns\ds+\sup_{t\in[0,T]}\dbE\(\int_0^t|\si_1(t,s,X_1(s),\G^b_1(t,s,X_1(s)))-\si_2(t,s,X_1(s),\G^b_2(t,s,X_1(s)))|^2
ds\)^{p\over2}\Big\}.\ea\ee

\ms

\it Proof. \rm By (H2), similar to (\ref{2.7})--(\ref{2.8}), making
use of (\ref{th-b-growth}), for any $X(\cd)\in
L^p_\dbF(0,T;\dbR^n)$, we have
\bel{}|\G^b(t,s,X(s))|+|\G^\si(t,s,X(s))|\le
L\(1+\dbE|X(s)|+|X(s)|\).\ee
Thus, if $X(\cd)\in L^p_\dbF(0,T;\dbR^n)$ is a solution to
(\ref{MF-FSVIE1}) with $\f(\cd)\in L^p_\dbF(0,T;\dbR^n)$, then by
(\ref{b-si-growth}),
\bel{2.29}\ba{ll}
\ns\ds\dbE|X(t)|^p\le
3^{p-1}\dbE\Big\{|\f(t)|^p+\Big|\int_0^tb(t,s,X(s),\G^b(t,s,X(s)))ds\Big|^p\\
\ns\ds\qq\qq\qq\qq+\Big|\int_0^t\si(t,s,X(s),\G^\si(t,s,X(s)))dW(s)\Big|^p\Big\}\\
\ns\ds\le
3^{p-1}\Big\{\dbE|\f(t)|^p+\dbE\(\int_0^tL\[1+|X(s)|+|\G^b(t,s,X(s))|\]ds\)^p\\
\ns\ds\qq\qq+\dbE\(\int_0^tL^2\[1+|X(s)|+|\G^\si(t,s,X(s))|\]^2ds\)^{p\over2}\Big\}\\
\ns\ds\le K\Big\{1+\dbE|\f(t)|^p+\int_0^t|X(s)|^pds\Big\}.\ea\ee
Consequently,
$$\int_0^t\dbE|X(r)|^pdr\le
K\Big\{1+\int_0^t\dbE|\f(r)|^pdr+\int_0^t\[\int_0^r|X(s)|^pds\]dr\Big\},\qq0\le
t\le T.$$
Using Gronwall's inequality, we obtain (\ref{|X|Lp-estimate}).

\ms

Now, let $\d>0$ be undetermined. For any $x(\cd)\in
L^p_\dbF(0,\d;\dbR^n)$, define
$$\ba{ll}
\ns\ds\cG(x(\cd))(t)=\f(t)+\int_0^tb(t,s,x(s),\G^b(t,s,x(s)))ds\\
\ns\ds\qq\qq\qq+\int_0^t\si(t,x,x(s),\G^\si(t,s,x(s)))dW(s),\qq
t\in[0,\d].\ea$$
Then we have
$$\ba{ll}
\ns\ds\dbE\int_0^\d|\cG(x(\cd))(t)|^pdt\le
K\dbE\Big\{\int_0^\d|\f(t)|^pdt+\int_0^\d\(\int_0^t\(1+|x(s)|+|\G^b(t,s,x(s))|\)ds\)^p\\
\ns\ds\qq\qq\qq\qq\qq\qq+\int_0^\d\Big|\int_0^t\si(t,s,x(s),\G^\si(t,s,x(s)))dW(s)\Big|^pdt
\Big\}\\
\ns\ds\le
K\Big\{\dbE\int_0^\d|\f(t)|^pdt+\dbE\int_0^\d|x(t)|^pdt\Big\}.\ea$$
Thus, $\cG:L^p_\dbF(0,\d;\dbR^n)\to L^p_\dbF(0,\d;\dbR^n)$. Next,
for any $x_1(\cd),x_2(\cd)\in L^p_\dbF(0,\d;\dbR^n)$, we have
(making use of (\ref{b-si-Lip}) and (\ref{th-b-Lip}))
$$\ba{ll}
\ns\ds\dbE\int_0^\d|\cG(x_1(\cd))(t)-\cG(x_2(\cd))(t)|^pdt\\
\ns\ds\le2^{p-1}\Big\{\dbE\int_0^\d\[\int_0^tL\(|x_1(s)-x_2(s)|+|\G^b(t,s,x_1(s))-\G^b(t,s,x_2(s))|\)ds\]^pdt\\
\ns\ds\q+\int_0^\d\dbE\[\int_0^tL^2\(|x_1(s)-x_2(s)|^2+|\G^\si(t,s,x_1(s))-\G^\si(t,s,x_2(s))|^2\)
ds\]^{p\over2}dt\Big\}\\
\ns\ds\le K_0\d\dbE\int_0^\d|x_1(t)-x_2(t)|^pdt,\ea$$
with $K_0>0$ being an absolute constant (only depending on $L$ and
$p$). Then letting $\d={1\over2K_0}$, we see that
$\cG:L^p_\dbF(0,\d;\dbR^n)\to L^p_\dbF(0,\d;\dbR^n)$ is a
contraction. Hence, MF-FSVIE (\ref{MF-FSVIE1}) admits a unique
solution $X(\cd)\in L^p_\dbF(0,\d;\dbR^n)$.

\ms

Next, for $t\in[\d,2\d]$, we write (\ref{MF-FSVIE}) as
\bel{MF-FSVIE2}\ba{ll}
\ns\ds X(t)=\h\f(t)+\int_\d^tb(t,s,X(s),\G^b(t,s,X(s)))ds\\
\ns\ds\qq\qq\qq+\int_\d^t\si(t,s,X(s),\G^\si(t,s,X(s))dW(s),\ea\ee
with
$$\ba{ll}
\ns\ds\h\f(t)=\f(t)+\int_0^\d
b(t,s,X(s),\G^b(t,s,X(s)))ds+\int_0^\d\si(t,s,X(s),\G^\si(t,s,X(s))dW(s).\ea$$
Then a similar argument as above applies to obtain a unique solution
of (\ref{MF-FSVIE2}) on $[\d,2\d]$. It is important to note that the
step-length $\d>0$ is uniform. Hence, by induction, we obtain the
unique solvability of (\ref{MF-FSVIE1}) on $[0,T]$.

\ms

Now, for $i=1,2$, let $X_i(\cd)\in L^p_\dbF(0,T;\dbR^n)$ be the
solutions of (\ref{MF-FSVIE1}) corresponding to $\f_i(\cd)\in
L^p_\dbF(0,T;\dbR^n)$ and
$b_i(\cd),\si_i(\cd),\th_i^b(\cd),\th_i^\si(\cd)$ (satisfying
(H1)--(H2)). Then
$$\ba{ll}
\ns\ds\dbE|X_1(t)-X_2(t)|^p\le3^{p-1}\dbE\Big\{|\f_1(t)-\f_2(t)|^p\\
\ns\ds\qq+\(\int_0^t|b_1(t,s,X_1(s),\G^b_1(t,s,X_1(s)))-b_2(t,s,X_2(s),\G^b_2(t,s,X_2(s)))|ds\)^p\\
\ns\ds\qq+\Big|\int_0^t\(\si_1(t,s,X_1(s),\G^b_1(t,s,X_1(s)))-\si_2(t,s,X_2(s),\G^b_2(t,s,X_2(s)))\)
dW(s)\Big|^p\Big\}\\
\ns\ds\le
K\Big\{\dbE|\f_1(t)-\f_2(t)|^p+\dbE\int_0^t|X_1(s)-X_2(s)|^pds\\
\ns\ds\qq+\dbE\(\int_0^t|b_1(t,s,X_1(s),\G^b_1(t,s,X_1(s)))-b_2(t,s,X_1(s),\G^b_2(t,s,X_1(s)))|ds\)^p\\
\ns\ds\qq+\dbE\Big|\int_0^t\(\si_1(t,s,X_1(s),\G^b_1(t,s,X_1(s)))-\si_2(t,s,X_1(s),\G^b_2(t,s,X_1(s)))\)
dW(s)\Big|^p\Big\}.\ea$$
Then we can obtain estimate (\ref{|X-X|Lp-estimate}).

\ms

The conclusions under (H1)$'$--(H2)$'$ are easy to obtain. \endpf

\ms

\subsection{Linear MF-FSVIEs and MF-BSVIEs.}

Let us now look at linear MF-FSVIEs, by which we mean the
following:
\bel{L-MF-FSVIE}\ba{ll}
\ns\ds X(t)=\f(t)+\int_0^t\(A_0(t,s)X(s)+\dbE'\[C_0(t,s)X(s)\]\)ds\\
\ns\ds\qq\qq\q+\int_0^t\(A_1(t,s)X(s)+\dbE'\[C_1(t,s)X(s)\]\)dW(s),\qq
t\in[0,T].\ea\ee
For such an equation, we introduce the following hypotheses.

\ms

{\bf(L1)} The maps
$$A_0,A_1:\D^*\times\O\to\dbR^{n\times n},\q
C_0,C_1:\D^*\times\O^2\to\dbR^{n\times n},$$
are measurable and uniformly bounded. For any $t\in[0,T]$,
$s\mapsto(A_0(t,s),A_1(t,s))$ is $\dbF$-progressively measurable on
$[0,t]$, and $s\mapsto(C_0(t,s),C_1(t,s))$ is $\dbF^2$-progressively
measurable on $[0,t]$.

\ms

{\bf(L1)$'$} In addition to (L1), the maps
$$t\mapsto(A_0(t,s,\o),A_1(t,s,\o),C_0(t,s,\o,\o'),C_1(t,s,\o,\o'))$$
is continuous on $[s,T]$.

\ms

Clearly, by defining
$$\left\{\ba{ll}
\ns\ds
b(t,s,\o,x,\g)=A_0(t,s,\o)x+\g,\q\th^b(t,s,\o,\o',x,x')=C_0(t,s,\o,\o')x',\\
\ns\ds\si(t,s,\o,x,\g')=A_1(t,s,\o)x+\g',\q\th^\si(t,s,\o,\o',x,x')=C_1(t,s,\o,\o')x',\ea\right.$$
we see that (\ref{L-MF-FSVIE}) is a special case of
(\ref{MF-FSVIE1}). Moreover, (L1) implies (H1)--(H2), and (L1)$'$
implies (H1)$'$--(H2)$'$. Hence, we have the following corollary of
Theorem 2.6.

\ms

\bf Corollary 2.7. \sl Let {\rm(L1)} hold, and $p\ge2$. Then for any
$\f(\cd)\in L^p_\dbF(0,T;\dbR^n)$, $(\ref{L-MF-FSVIE})$ admits a
unique solution $X(\cd)\in L^p_\dbF(0,T;\dbR^n)$, and estimate
$(\ref{|X|Lp-estimate})$ holds. Further, let $p>2$. If for $i=1,2$,
$A_0^i(\cd),A_1^i(\cd),C_0^i(\cd)$, $C_1^i(\cd)$ satisfy {\rm(L1)},
$\f_i(\cd)\in L^p_\dbF(0,T;\dbR^n)$, and $X_i(\cd)\in
L^p_\dbF(0,T;\dbR^n)$ are the corresponding solutions to
$(\ref{L-MF-FSVIE})$, then for any $r\in(2,p)$,
\bel{stability3}\2n\ba{ll}
\ns\ds\dbE\int_0^T|X_1(t)-X_2(t)|^rdt\le
K\dbE\int_0^T|\f_1(t)-\f_2(t)|^rdt+K\(1+\dbE\int_0^T|\f_1(t)|^pdt\)^{r\over
p}\\
\ns\ds\cd\2n\int_0^T\3n\Big\{\[\dbE\(\1n\int_0^t\3n|A^1_0(t,s)\1n-\1n
A^2_0(t,s)|^{r\over r-1}ds\1n\)\1n^{(r-1)p\over p-r}\]\1n^{p-r\over
p}\3n+\2n\[\dbE^2\1n\(\1n\int_0^t\2n|C_0^1(t,s)\1n-\1n
C_0^2(t,s)|^{r\over
r-1}ds\1n\)\1n^{(r-1)p\over p-r}\]^{p-r\over p}\\
\ns\ds+\1n\[\dbE\(\1n\int_0^t\2n|A^1_1(t,s)\1n-\1n
A_1^2(t,s)|^{2r\over
r-2}ds\1n\)\1n^{(r-2)p\over2(p-r)}\]\1n^{p-r\over
p}\3n+\2n\[\dbE^2\(\1n\int_0^t\2n|C_1^1(t,s)\1n-\1n
C^2_1(t,s)|^{2r\over r-2}ds\)^{(r-2)p\over2(p-r)}\]^{p-r\over
p}\1n\Big\}dt.\ea\ee
Moreover, let {\rm(L1)$'$} hold. Then for any $\f(\cd)\in
C_\dbF^p([0,T];\dbR^n)$, $(\ref{L-MF-FSVIE})$ admits a unique
solution $X(\cd)\in C_\dbF^p([0,T];\dbR^n)$, and estimate
$(\ref{|X|-estimate})$ holds. Now for $i=1,2$, let
$A_0^i(\cd),A_1^i(\cd),C_0^i(\cd)$, $C_1^i(\cd)$ satisfy
{\rm(L1)$'$}, $\f_i(\cd)\in C^p_\dbF([0,T];\dbR^n)$, and
$X_i(\cd)\in C^p_\dbF([0,T];\dbR^n)$ be the corresponding solutions
to $(\ref{L-MF-FSVIE})$, then for any $2<r<p$,
\bel{|X-X|}\2n\ba{ll}
\ns\ds\sup_{t\in[0,T[}\dbE|X_1(t)-X_2(t)|^r\le
K\sup_{t\in[0,T]}\dbE|\f_1(t)-\f_2(t)|^r\\
\ns\ds\q+K\(1+\sup_{t\in[0,T]}\dbE|\f_1(t)|^p\)^{r\over
p}\Big\{\sup_{t\in[0,T]}\[\dbE\(\1n\int_0^t\3n|A^1_0(t,s)\1n-\1n
A^2_0(t,s)|^{r\over r-1}ds\1n\)\1n^{(r-1)p\over p-r}\]\1n^{p-r\over
p}\\
\ns\ds\q+\sup_{t\in[0,T]}\[\dbE^2\1n\(\1n\int_0^t\2n|C_0^1(t,s)\1n-\1n
C_0^2(t,s)|^{r\over
r-1}ds\1n\)\1n^{(r-1)p\over p-r}\]^{p-r\over p}\\
\ns\ds\q+\sup_{t\in[0,T]}\[\dbE\(\1n\int_0^t\2n|A^1_1(t,s)\1n-\1n
A_1^2(t,s)|^{2r\over
r-2}ds\1n\)\1n^{(r-2)p\over2(p-r)}\]\1n^{p-r\over
p}\\
\ns\ds\q+\sup_{t\in[0,T]}\[\dbE^2\(\1n\int_0^t\2n|C_1^1(t,s)\1n-\1n
C^2_1(t,s)|^{2r\over r-2}ds\)^{(r-2)p\over2(p-r)}\]^{p-r\over
p}\1n\Big\}.\ea\ee

\ms

\it Proof. \rm We need only to prove the stability estimate. Let
$X_i(\cd)\in L^p_\dbF(0,T;\dbR^n)$ be the solutions to the linear
MF-FSVIEs corresponding to the coefficients
$(A_0^i(\cd),C^i_0(\cd),A^i_1(\cd),C^i_1(\cd))$ satisfying (L1) and
free term $\f_i(\cd)\in L^p_\dbF(0,T;\dbR^n)$. Then we have
$$\dbE\int_0^T|X_i(s)|^pds\le K\(1+\dbE\int_0^T|\f_i(s)|^p\).$$
Now, for any $2<r<p$,
$$\3n\3n\3n\3n\ba{ll}
\ns\ds\dbE\int_0^T|X_1(t)-X_2(t)|^r\le K\Big\{\dbE\int_0^T|\f_1(t)-\f_2(t)|^rdt\\
\ns\ds+\dbE\int_0^T\[\int_0^t\(|A_0^1(t,s)-A_0^2(t,s)||X_1(s)|+\dbE'[\,|C_0^1(t,s)
-C_0^2(t,s)||X_1(s)|\,]\)
ds\]^rdt\\
\ns\ds+\dbE\int_0^T\[\int_0^t\(|A^1_1(t,s)-A_1^2(t,s)|^2|X_1(s)|^2
+\dbE'[\,|C_1^1(t,s)-C^2_1(t,s)|^2|X_1(s)|^2\,]\)
ds\]^{r\over2}dt\Big\}\\
\ns\ds\le\1n
K\1n\Big\{\dbE\int_0^T|\f^1(t)-\f^2(t)|^rdt+\dbE\int_0^T\(\int_0^t|A^1_0(t,s)-A^2_0(t,s)|
|X_1(s)|ds\)^rdt\\
\ns\ds\q+\dbE^2\3n\int_0^T\3n\1n\(\1n\int_0^t\2n|C_0^1(t,s)\1n-\1n
C_0^2(t,s)||X_1\1n(s)|ds\)^r\2n
dt\1n+\1n\dbE\2n\int_0^T\3n\1n\(\1n\int_0^t\2n|A^1_1(t,s)\1n-\1n
A_1^2(t,s)|^2|X_1\1n(s)|^2
ds\1n\)^{\1n{r\over2}}dt\\
\ns\ds\q+\dbE^2\int_0^T\(\int_0^t|C_1^1(t,s)-C^2_1(t,s)|^2|X_1(s)|^2ds\)^{r\over2}dt\Big\}\\
\ns\ds\le\1n K\1n\Big\{\1n\dbE\3n\int_0^T\3n|\f^1(t)-\f^2(t)|^r\1n
dt\1n+\1n\dbE\2n\int_0^T\3n\(\int_0^t|A^1_0(t,s)-A^2_0(t,s)|^{r\over
r-1}ds\)^{r-1}\(\int_0^t
|X_1(s)|^rds\)dt\\
\ns\ds\q+\dbE^2\int_0^T\(\int_0^t|C_0^1(t,s)-C_0^2(t,s)|^{r\over
r-1}ds\)^{r-1}
\(\int_0^t|X_1(s)|^rds\)dt\\
\ns\ds\q+\dbE\int_0^T\(\int_0^t|A^1_1(t,s)-A_1^2(t,s)|^{2r\over
r-2}ds\)^{r-2\over2}\(\int_0^t|X_1(s)|^r
ds\)dt\\
\ns\ds\q+\dbE^2\int_0^T\(\int_0^t|C_1^1(t,s)-C^2_1(t,s)|^{2r\over
r-2}ds\)^{r-2\over2}\(\int_0^t|X_1(s)|^rds\)dt\Big\}\\
\ns\ds\le
K\dbE\int_0^T|\f^1(t)-\f^2(t)|^rdt+K\dbE\int_0^T\[\dbE\(\int_0^t
|X_1(s)|^rds\)^{p\over r}\]^{r\over p}\\
\ns\ds\q\cd\Big\{\1n\[\dbE\(\1n\int_0^t\2n|A^1_0(t,s)\1n-\1n
A^2_0(t,s)|^{r\over r-1}ds\)\1n^{(r-1)p\over p-r}\]\1n^{p-r\over
p}\3n+\2n\[\dbE^2\(\1n\int_0^t\2n|C_0^1(t,s)\1n-\1n
C_0^2(t,s)|^{r\over
r-1}ds\1n\)\1n^{(r-1)p\over p-r}\]^{p-r\over p}\\
\ns\ds\q+\1n\[\dbE\(\1n\int_0^t\2n|A^1_1(t,s)\1n-\1n
A_1^2(t,s)|^{2r\over r-2}ds\)^{(r-2)p\over2(p-r)}\]\1n^{p-r\over
p}\3n+\2n\[\dbE^2\1n\(\1n\int_0^t|C_1^1(t,s)\1n-\1n
C^2_1(t,s)|^{2r\over
r-2}ds\1n\)\1n^{(r-2)p\over2(p-r)}\]\1n^{p-r\over
p}\1n\Big\}dt.\ea$$
Then (\ref{stability3}) follows. The case that (L1)$'$ holds case be
proved similarly. \endpf

\ms

We point out that linear MF-FSVIE (\ref{L-MF-FSVIE}) is general
enough in some sense. To see this, let us formally look at the
variational equation of (\ref{MF-FSVIE1}). More precisely, let
$X^\d(\cd)$ be the unique solution of (\ref{MF-FSVIE1}) with
$\f(\cd)$ replaced by $\f(\cd)+\d\bar\f(\cd)$. We formally let
$$\bar X(t)=\lim_{\d\to0}{X^\d(t)-X(t)\over\d}.$$
Then $\bar X(\cd)$ should satisfy the following linear MF-FSVIE:
\bel{bar X}\ba{ll}
\ns\ds\bar X(t)=\bar\f(t)+\int_0^t\(b_x(t,s)\bar
X(s)+b_\g(t,s)\dbE'\[\th^b_x(t,s)\bar X(s,\o)+\th^b_{x'}(t,s)\bar
X(s,\o')\]\)ds\\
\ns\ds\qq\q+\int_0^t\(\si_x(t,s)\bar
X(s)+\si_\g(t,s)\dbE'\[\th^\si_x(t,s)\bar
X(s,\o)+\th^\si_{x'}(t,s)\bar X(s,\o')\]\)dW(s),\ea\ee
where (with a little misuse of $\g$)
\bel{bar X2}\left\{\ba{ll}
\ns\ds
b_x(t,s)=b_x(t,s,\o,\G^b(t,s,X(s,\o))),\q\th^b_x(t,s)=\th^b_x(t,s,\o,\o',X(s,\o),X(s,\o')),\\
\ns\ds
b_\g(t,s)=b_\g(t,s,\o,\G^b(t,s,X(s,\o))),\q\th^b_{x'}(t,s)=\th^b_{x'}(t,s,\o,\o',X(s,\o),X(s,\o')),\\
\ns\ds
\si_x(t,s)=\si_x(t,s,\o,\G^\si(t,s,X(s,\o))),\q\th^\si_x(t,s)=\th^\si_x(t,s,\o,\o',X(s,\o),X(s,\o')),\\
\ns\ds
\si_\g(t,s)=\si_\g(t,s,\o,\G^\si(t,s,X(s,\o))),\q\th^\si_{x'}(t,s)=\th^\si_{x'}(t,s,\o,\o',X(s,\o),X(s,\o')).
\ea\right.\ee
It is interesting to note that (\ref{bar X}) can be written as
follows:
\bel{bar X3}\ba{ll}
\ns\ds\bar
X(t)=\bar\f(t)+\int_0^t\Big\{\(b_x(t,s)+b_\g(t,s)\dbE'[\th^b_x(t,s)]\)\bar
X(s)+\dbE'\[b_\g(t,s)\th^b_{x'}(t,s)\bar
X(s)\]\Big\}ds\\
\ns\ds\qq\q+\int_0^t\Big\{\(\si_x(t,s)+\si_\g(t,s)\dbE'[\th^\si_x(t,s)]\)\bar
X(s,\o)+\dbE'\[\si_\g(t,s)\th^\si_{x'}(t,s)\bar
X(s,\o')\]\Big\}dW(s),\ea\ee
which is a special case of (\ref{L-MF-FSVIE}).

\ms

Mimicking the above, we see that general linear MF-BSVIE should take
the following form:
\bel{L-MF-BSVIE}\ba{ll}
\ns\ds
Y(t)=\psi(t)+\int_t^T\(\bar A_0(t,s)Y(s)+\bar B_0(t,s)Z(t,s)+\bar C_0(t,s)Z(s,t)\\
\ns\ds\qq\qq+\dbE'\[\bar A_1(t,s)Y(s)+\bar B_1(t,s)Z(t,s)+\bar
C_1(t,s)Z(s,t)\]\)ds-\int_t^TZ(t,s)dW(s).\ea\ee
For the coefficients, we should adopt the following hypothesis.

\ms

{\bf(L2)} The maps
$$\bar A_0,\bar B_0,\bar C_0:\D\times\O\to\dbR^{n\times n},\q
\bar A_1,\bar B_1,\bar C_1:\D\times\O^2\to\dbR^{n\times n}$$
are measurable and uniformly bounded. Moreover, for any $t\in[0,T]$,
$s\mapsto(\bar A_0(t,s),\bar B_0(t,s),\bar C_0(t,s))$ is
$\dbF$-progressively measurable on $[t,T]$, and $s\mapsto(\bar
A_1(t,s),\bar B_1(t,s),\bar C_1(t,s))$ is $\dbF^2$-progressively
measurable on $[t,T]$.

\ms

We expect that under (L2), for reasonable $\psi(\cd)$, the above
(\ref{L-MF-BSVIE}) will have a unique adapted M-solution. Such a
result will be a consequence of the main result of the next section.

\ms

\section{Well-posedness of MF-BSVIEs.}

In this section, we are going to establish the well-posedness of our
MF-BSVIEs. To begin with, let us introduce the following
hypothesis.

\ms

{\bf(H3)$_q$} The map $\th:\D\times\O^2\times\dbR^{6n}\to\dbR^m$
satisfies (H0)$_q$. The map
$g:\D\times\O\times\dbR^{3n}\times\dbR^m\to\dbR^n$ is measurable and
for all $(t,y,z,\h z,\g)\in[0,T]\times\dbR^{3n}\times\dbR^m$, the
map $(s,\o)\mapsto g(t,s,\o,y,z,\h z,\g)$ is $\dbF$-progressively
measurable. Moreover, there exist constants $L>0$ and $q\ge2$ such
that
\bel{g-Lip}\ba{ll}
\ns\ds|g(t,s,\o,y_1,z_1,\h z_1,\g_1)-g(t,s,\o,y_2,z_2,\h z_2,\g_2)|\\
\ns\ds\le L\(|y_1-y_2|+|z_1-z_2|+|\h z_1-\h
z_2|+|\g_1-\g_2|\),\\
\ns\ds\qq\qq\qq\qq\forall(t,s,\o)\in\D\times\O,(y_i,z_i,\h
z_i,\g_i)\in\dbR^{3n}\times\dbR^m,i=1,2,\ea\ee
and
\bel{g-growth}\ba{ll}
\ns\ds|g(t,s,\o,y,z,\h z,\g)|\le L\(1+|y|+|z|+|\h z|^{2\over q}+|\g|\),\\
\ns\ds\qq\qq\qq\qq\qq\forall(t,s,\o)\in\D\times\O,(y,z,\h
z,\g)\in\dbR^{3n}\times\dbR^m.\ea\ee

Similar to (H0)$_\infty$ in the previous section, (H3)$_\infty$ is
understood that (\ref{th-growth2}) holds and (\ref{g-growth}) is
replaced by the following:
\bel{g-growth2}\ba{ll}
\ns\ds|g(t,s,\o,y,z,\h z,\g)|\le L\(1+|y|+|z|+|\g|\),\\
\ns\ds\qq\qq\qq\qq\qq\forall(t,s,\o)\in\D\times\O,(y,z,\h
z,\g)\in\dbR^{3n}\times\dbR^m.\ea\ee

\subsection{A special MF-BSVIE.}

In this subsection, we firstly consider the following special
MF-BSVIE:
\bel{MF-BSVIE-s}Y(t)=\psi(t)+\int_t^T\wt g(t,s,Z(t,s),\wt
\G(t,s,Z(t,s)))ds-\int_t^TZ(t,s)dW(s),\q t\in[0,T],\ee
where
$$\left\{\ba{ll}
\ns\ds\wt g(t,s,Z,\g)=g(t,s,y(s),Z,z(s,t),\g),\\
\ns\ds\wt\G(t,s,Z)=\G(t,s,y(s),Z,z(s,t))\equiv\dbE'\[\th(t,s,y(s),Z,z(s,t),
y',z',\h z')\]_{(y',z',\h z')=(y(s),Z,z(s,t))},\ea\right.$$
for some given $(y(\cd),z(\cd\,,\cd))\in\cM^p[0,T]$. Therefore,
$$\wt g(t,s,Z(t,s),\wt\G(t,s,Z(t,s)))=g(t,s,y(s),Z(t,s),z(s,t),\G(t,s,y(s),Z(t,s),z(s,t))).$$
Note that we may take much more general $\wt g(\cd)$ and
$\wt\G(\cd)$. But the above is sufficient for our purpose, and by
restricting such a case, we avoid stating a lengthy assumption
similar to (H3)$_q$. We now state and prove the following result
concerning MF-BSVIE (\ref{MF-BSVIE-s}).

\ms

\bf Proposition 3.1. \sl Let {\rm(H3)$_q$} hold. Then for any $p>1$
and $\psi(\cd)\in L^p_{\cF_T}(0,T;\dbR^n)$, MF-BSVIE
$(\ref{MF-BSVIE-s})$ admits a unique M-solution
$(Y(\cd),Z(\cd\,,\cd))\in\cM^p[0,T]$. Moreover, the following
estimate holds:
\bel{Lp-s}\ba{ll}
\ns\ds\dbE\[|Y(t)|^p+\(\int_t^T|Z(t,s)|^2ds\)^{p\over2}\]\le
K\dbE\[|\psi(t)|^p+\(\int_t^T|\wt
g(t,s,0,\wt\G(t,s,0))|ds\)^p\].\ea\ee
Further, for $i=1,2$, let $\psi_i(\cd)\in L^p_{\cF_T}(0,T;\dbR^n)$,
$(y_i(\cd),z_i(\cd\,,\cd))\in\cM^p[0,T]$, and
$$\left\{\ba{ll}
\ns\ds\wt
g_i(t,s,Z(t,s),\wt\G_i(t,s,Z(t,s))=g_i(t,s,y_i(s),Z(t,s),z_i(s,t),\G_i(t,s,y_i(s),Z(t,s),z_i(s,t)),\\
\ns\ds\G_i(t,s,Y,Z,\h Z)=\dbE'\[\th_i(t,s,y_i(s),Z,z_i(s,t),y',z',\h
z')\]_{(y',z',\hat z')=(y_i(s),Z,z_i(s,t))}\ea\right.$$
with $g_i(\cd)$ and $\th_i(\cd)$ satisfying {\rm(H3)$_q$}. Then the
corresponding M-solutions $(Y_i(\cd),Z_i(\cd))$ satisfy the
following stability estimate:
\bel{stability-s}\ba{ll}
\ns\ds\dbE\[|Y_1(t)-Y_2(t)|^p+\(\int_t^T|Z_1(t,s)-Z_2(t,s)|^2ds\)^{p\over2}\]\le
K\dbE\[|\psi_1(t)-\psi_2(t)|^p\\
\ns\ds\qq+\(\int_t^T|\wt g_1(t,s,Z_1(t,s),\wt\G_1(t,s,Z_1(t,s))-\wt
g_2(t,s,Z_1(t,s),\wt\G_2(t,s,Z_1(t,s))|ds\)^p\].\ea\ee

\ms

\it Proof. \rm Fix $t\in[0,T)$. Consider the following MF-BSDE
(parameterized by $t$):
\bel{BSVIE-s}\eta(r)=\psi(t)+\int_r^T\wt g(t,s,\z(s),\wt
\G(t,s,\z(s)))ds-\int_r^T\z(s)dW(s),\q r\in[t,T].\ee
If $p\in(1,2]$, it follows from (H3)$_{q}$ that
$$ \ba{ll} \ns\ds
\dbE\int_0^T\(\int_t^T|\wt g(t,s,0,\wt\G(t,s,0))|ds\)^pdt\\
\ns\ds \le K \dbE\int_0^T\(\int_t^T(|y(s)|+|z(s,t)|
+\dbE|y(s)|+\dbE|z(s,t)|)ds\)^pdt + K \\
\ns\ds \le
K \dbE\int_0^T\int_t^T|y(s)|^pdsdt + K \dbE\int_0^T\(\int_0^s |z(s,t)|^2dt\)^{p\over
2}ds + K <\infty, \ea $$
As to the case of $p=q>2$, similarly we have
$$ \ba{ll} \ns\ds
\dbE\int_0^T\(\int_t^T|\wt g(t,s,0,\wt\G(t,s,0))|ds\)^qdt\\
\ns\ds \le K \dbE\int_0^T\(\int_t^T(|y(s)|+|z(s,t)|^{2\over
q}+\dbE|y(s)|+\dbE|z(s,t)|^{2\over q})ds\)^qdt + K\\
\ns\ds \le
K \dbE\int_0^T\int_t^T|y(s)|^qdsdt + K \dbE\int_0^T \int_t^T|z(s,t)|^2ds dt + K <\infty.\ea $$
Similar to a standard argument for BSDEs, making use of contraction
mapping theorem, we can show that the above MF-BSDE admits a unique
adapted solution
$$(\eta(\cd),\z(\cd))\equiv\(\eta(\cd\,;t,\psi(t)),\z(\cd\,;t,\psi(t))\).$$
Moreover, the following estimate holds:
\bel{Lp-BSDE}\ba{ll}
\ns\ds\dbE\[\sup_{r\in[t,T]}|\eta(r;t,\psi(t))|^p+\(\int_t^T|\z(s;t,\psi(t))|^2ds\)^{p\over2}\]\\
\ns\ds\le K\dbE\[|\psi(t)|^p+\(\int_t^T|\wt
g(t,s,0,\wt\G(t,s,0))|ds\)^p\].\ea\ee
Further, for $i=1,2$, let $\psi_i(\cd)\in L^p_{\cF_T}(0,T;\dbR^n)$,
$(y_i(\cd),z_i(\cd\,,\cd))\in\cM^p[0,T]$, and
$$\left\{\ba{ll}
\ns\ds\wt
g_i(t,s,Z(t,s),\wt\G_i(t,s,Z(t,s))=g_i(t,s,y_i(s),Z(t,s),z_i(s,t),\G_i(t,s,y_i(s),Z(t,s),z_i(s,t)),\\
\ns\ds\G_i(t,s,Y,Z,\h Z)=\dbE'\[\th_i(t,s,y_i(s),Z,z_i(s,t),y',z',\h
z')\]_{(y',z',\hat z')=(y_i(s),Z,z_i(s,t))}\ea\right.$$
with $g_i(\cd)$ and $\th_i(\cd)$ satisfying (H3)$_q$. Then let
$(\eta_i(\cd),\z_i(\cd))$ be the adapted solutions of the
corresponding BSDE. It follows that
\bel{stability-BSDE}\ba{ll}
\ns\ds\dbE\[\sup_{r\in[t,T]}|\eta_1(r)-\eta_2(r)|^p\]
+\dbE\(\int_t^T|\z_1(s)-\z_2(s)|^2ds\)^{p\over2}\\
\ns\ds\le\2n
K\dbE\[|\psi_1(t)-\psi_2(t)|^p\2n+\2n\(\1n\int_t^T\3n|\wt
g_1(t,s,\z_1(s),\wt\G_1(t,s,\z_1(s))-\wt
g_2(t,s,\z_1(s),\wt\G_2(t,s,\z_1(s))|ds\)^p\].\ea\ee
Now, we define
$$Y(t)=\eta(t;t,\psi(t)),\q Z(t,s)=\z(s;t,\psi(t)),\qq(t,s)\in\D,$$
and $Z(t,s)$ on $\D^c$ through the martingale representation:
$$Y(t)=\dbE Y(t)+\int_0^tZ(t,s)dW(s),\qq t\in[0,T].$$
Then $(Y(\cd),Z(\cd\,,\cd))\in\cM^p[0,T]$ is the unique M-solution
to (\ref{MF-BSVIE-s}). Estimates (\ref{Lp-s}) and
(\ref{stability-s}) follows easily from (\ref{Lp-BSDE}) and
(\ref{stability-BSDE}), respectively. \endpf

\ms

Note that the cases that we are interested in are $p=2,q$. We will
use them below.

\subsection{The general case.}

Now, we consider our MF-BSVIEs. For convenience, let us rewrite
(\ref{MF-BSVIE}) here:
\bel{MF-BSVIE2}\ba{ll}
\ns\ds
Y(t)=\psi(t)+\int_t^Tg(t,s,Y(s),Z(t,s),Z(s,t),\G(t,s,Y(s),Z(t,s),Z(s,t)))ds\\
\ns\ds\qq\qq\qq-\int_t^TZ(t,s)dW(s),\qq\qq t\in[0,T],\ea\ee
with
\bel{G(t,s)}\ba{ll}
\ns\ds\G(t,s,Y(s),Z(t,s),Z(s,t))=\dbE'\[\th(t,s,Y(s),Z(t,s),Z(s,t))\]\\
\ns\ds=\int_\O\th(t,s,\o',\o,Y(s,\o'),Z(t,s,\o'),Z(s,t,\o'),Y(s,\o),
Z(t,s,\o),Z(s,t,\o))\dbP(d\o').\ea\ee
Our main result of this section is the following.

\ms

\bf Theorem 3.2. \sl Let {\rm(H3)$_q$} hold with $2\le q<\infty$.
Then for any $\psi(\cd)\in L^q_{\cF_T}(0,T;\dbR^n)$, MF-BSVIE
$(\ref{MF-BSVIE2})$ admits a unique adapted M-solution
$(Y(\cd),Z(\cd\,,\cd))\in\cM^q[0,T]$, and the following estimate
holds:
\bel{Lp-estimate}\|(Y(\cd),Z(\cd\,,\cd))\|_{\cM^q[0,T]}\le
K\(1+\|\psi(\cd)\|_{L^q_{\cF_T}(0,T;\dbR^n)}\).\ee
Moreover, for $i=1,2$, let $g_i(\cd)$ and $\th_i(\cd)$ satisfy
{\rm(H3)$_q$}, and $\psi_i(\cd)\in L^q_{\cF_T}(0,T;\dbR^n)$. Let
$(Y_i(\cd),Z_i(\cd\,,\cd))\in\cM^q[0,T]$ be the corresponding
adapted M-solutions. Then
\bel{stability}\ba{ll}
\ns\ds\|(Y_1(\cd),Z_1(\cd\,,\cd))-(Y_2(\cd),Z_2(\cd\,,\cd))\|^2_{\cM^2[0,T]}\\
\ns\ds\le
K\dbE\Big\{\int_0^T|\psi_1(t)-\psi_2(t)|^2dt+\int_0^T\(\int_t^T|(g_1-g_2)(t,s)|ds\)^2dt\Big\},\ea\ee
where
$$\ba{ll}
\ns\ds(g_1-g_2)(t,s)=g_1(t,s,Y_1(s),Z_1(t,s),Z_1(s,t),\G_1(t,s,Y_1(s),Z_1(t,s),Z_1(s,t)))\\
\ns\ds\qq\qq\qq\qq-g_2(t,s,Y_1(s),Z_1(t,s),Z_1(s,t),\G_2(t,s,Y_1(s),Z_1(t,s),Z_1(s,t))),\ea$$
with
$$\ba{ll}
\ns\ds\G_i(t,s,Y_i(s),Z_i(t,s),Z_i(s,t))\\
\ns\ds=\dbE'\[\th_i(t,s,Y_i(s),Z_i(t,s),Z_i(s,t),y,z,\h
z)\]_{(y,z,\hat z)
=(Y_i(s),Z_i(t,s),Z_i(s,t))}\\
\ns\ds\equiv\int_\O\th_i(t,s,\o',\o,Y_i(s,\o'),Z_i(t,s,\o'),Z_i(s,t,\o'),Y_i(s,\o),
Z_i(t,s,\o),Z_i(s,t,\o))\dbP(d\o').\ea$$

\ms

\it Proof. \rm We split the proof into several steps.

\ms

\it Step 1. \rm Existence and uniqueness of M-solutions of
(\ref{MF-BSVIE2}) in $\cM^p[0,T]$ with $p\in(1,2]$.

\ms

Let $\psi(\cd)\in L_{\cF_T}^p(0,T;\dbR^n)$ be given. For any
$(y(\cd),z(\cd\,,\cd))\in\cM^p[0,T]$, we consider the following
MF-BSVIE
\bel{BSVIE(yz)}\ba{ll}
\ns\ds Y(t)=\psi(t)+\int_t^Tg(t,s,y(s),Z(t,s),z(s,t),\G(t,s,y(s),Z(t,s),z(s,t)))ds\\
\ns\ds\qq\qq\qq-\int_t^TZ(t,s)dW(s),\qq\qq t\in[0,T].\ea\ee
According to Proposition 3.1, there exists a unique adapted
M-solution $(Y(\cd),Z(\cd\,,\cd))\in\cM^p[0,T]$. Moreover, the
following estimate holds: (making use of (\ref{2.7} ) and (\ref{2.8} ))
\bel{3.12}\ba{ll}
\ns\ds\dbE\Big\{|Y(t)|^p+\(\int_t^T|Z(t,s)|^2ds\)^{p\over2}\Big\}\\
\ns\ds\le
K\dbE\Big\{|\psi(t)|^p+\(\int_t^T|g(t,s,y(s),0,z(s,t),\G(t,s,y(s),0,z(s,t))|ds\)^p\Big\}\\
\ns\ds\le K\dbE\Big\{1+|\psi(t)|^p+\[\int_t^T\(|y(s)|+|z(s,t)|+|\G(t,s,y(s),0,z(s,t))|\)ds\]^p\Big\}\\
\ns\ds\le
K\dbE\Big\{\1n1\1n+\1n|\psi(t)|^p\2n+\1n\[\int_t^T\2n\(|y(s)|+|z(s,t)|+\dbE|y(s)|+\dbE|z(s,t)|+|y(s)|+|z(s,t)|\)ds\]^p\Big\}\\
\ns\ds\le
K\dbE\Big\{1+|\psi(t)|^p+\[\int_t^T\(|y(s)|+|z(s,t)|\)ds\]^p\Big\}\\
\ns\ds\le
K\dbE\Big\{1+|\psi(t)|^p+\int_t^T|y(s)|^pds+\int_t^T|z(s,t)|^pds\Big\}.\ea\ee
Consequently, (making use of (\ref{2.14}) for
$(y(\cd),z(\cd\,,\cd))\in\cM^p[0,T]$)
$$\ba{ll}
\ns\ds\|(Y(\cd),Z(\cd\,,\cd))\|^p_{\cH^p[0,T]}\equiv\dbE\Big\{\int_0^T|Y(t)|^pdt
+\int_0^T\(\int_0^T|Z(t,s)|^2ds\)^{p\over2}dt\Big\}\\
\ns\ds\le
K\dbE\Big\{1+\int_0^T|\psi(t)|^pdt+\int_0^T\int_t^T|y(s)|^pdsdt+\int_0^T\int_t^T|z(s,t)|^pdsdt\Big\}\\
\ns\ds\le
K\dbE\Big\{1+\int_0^T|\psi(t)|^pdt+\int_0^T|y(t)|^pdt+\int_0^T\int_0^t|z(t,s)|^pdsdt\Big\}\\
\ns\ds\le
K\dbE\Big\{1+\int_0^T|\psi(t)|^pdt+\int_0^T|y(t)|^pdt+\int_0^T\(\int_0^t|z(t,s)|^2ds\)^{p\over2}dt\Big\}\\
\ns\ds\le
K\dbE\Big\{1+\int_0^T|\psi(t)|^pdt+\int_0^T|y(t)|^pdt\Big\}\\
\ns\ds\le
K\Big\{1+\|\psi(\cd)\|^p_{L^p_{\cF_T}(0,T;\dbR^n)}+\|(y(\cd),z(\cd\,,\cd))\|_{\cM^p[0,T]}\Big\}.\ea$$
Hence, if we define
$\Th(y(\cd),z(\cd\,,\cd))=(Y(\cd),Z(\cd\,,\cd))$, then $\Th$ maps
from $\cM^p[0,T]$ to itself. We now show that the mapping $\Th$ is
contractive. To this end, take any
$(y_i(\cd),z_i(\cd\,,\cd))\in\cM^p[0,T]$ ($i=1,2$), and let
$$(Y_i(\cd),Z_i(\cd\,,\cd))=\Th(y_i(\cd),z_i(\cd\,,\cd)).$$
Then by Proposition 3.1, we have (note (\ref{2.8}))
$$\ba{ll}
\ns\ds\dbE\[|Y_1(t)-Y_2(t)|^p+\(\int_t^T|Z_1(t,s)-Z_2(t,s)|^2ds\)^{p\over2}\]\\
\ns\ds\le K\dbE\(\int_t^T|g(t,s,y_1(s),Z_1(t,s),z_1(s,t),\G(t,s,y_1(s),Z_1(t,s),z_1(s,t)))\\
\ns\ds\qq-
g(t,s,y_2(s),Z_1(t,s),z_2(s,t),\G(t,s,y_2(s),Z_1(t,s),z_2(s,t)))|ds\)^p\\
\ns\ds\le K\dbE\[\int_t^T\(|y_1(s)-y_2(s)|+|z_1(s,t)-z_2(s,t)|\\
\ns\ds\qq+|\G(t,s,y_1(s),Z_1(t,s),z_1(s,t))-\G(t,s,y_2(s),Z_1(t,s),z_2(s,t))|\)ds\]^p\\
\ns\ds\le K\dbE\[\int_t^T\(|y_1(s)-y_2(s)|+|z_1(s,t)-z_2(s,t)|\\
\ns\ds\qq+|y_1(s)-y_2(s)|+|z_1(s,t))-z_2(s,t)|+\dbE|y_1(s)-y_2(s)|+\dbE|z_1(s,t)-z_2(s,t)|\)ds\]^p\\
\ns\ds\le
K\dbE\[\int_t^T\(|y_1(s)-y_2(s)|+|z_1(s,t)-z_2(s,t)|\)ds\]^p\\
\ns\ds\le
K\dbE\[\int_t^T|y_1(s)-y_2(s)|^pds+\(\int_t^T|z_1(s,t)-z_2(s,t)|ds\)^p\].\ea$$
Hence,
\bel{3.15}\ba{ll}
\ns\ds\|\Th(y_1(\cd),z_1(\cd\,,\cd))-\Th(y_2(\cd),z_2(\cd\,,\cd))\|^p_{\cM^p_\b[0,T]}\\
\ns\ds\equiv
\int_0^Te^{\b t}\dbE\[|Y_1(t)-Y_2(t)|^p+\(\int_t^T|Z_1(t,s)-Z_2(t,s)|^2ds\)^{p\over2}\]dt\\
\ns\ds\le K\int_0^Te^{\b
t}\dbE\[\int_t^T|y_1(s)-y_2(s)|^pds+\(\int_t^T|z_1(s,t)-z_2(s,t)|ds\)^p\]dt\\
\ns\ds=\1n K\dbE\[\1n\int_0^T\3n\(\1n\int_0^s\3n e^{\b
t}dt\)|y_1(s)\1n-\1n y_2(s)|^pds\1n+\2n\int_0^T\3n e^{\b
t}\(\1n\int_t^T\3n e^{-{\b\over p}s}e^{{\b\over p}s}|z_1(s,t)\1n-\1n z_2(s,t)|ds\)^pdt\]\\
\ns\ds\le\1n K\dbE\[{1\over\b}\2n\int_0^T\3n e^{\b t}|y_1(t)\1n-\1n
y_2(t)|^pdt\1n+\3n\int_0^T\3n e^{\b t}\(\1n\int_t^T\3n e^{-q\b
s\over
p}ds\)^{p\over q}\(\int_t^T\3n e^{\b s}|z_1(s,t)-z_2(s,t)|^pds\)dt\]\\
\ns\ds\le{K\over\b}\dbE\int_0^Te^{\b
t}|y_1(t)-y_2(t)|^pdt+K\({p\over\b q}\)^{p\over
q}\dbE\int_0^T\int_0^se^{\b
s}|z_1(s,t)-z_2(s,t)|^pdtds\\
\ns\ds\le{K\over\b}\dbE\int_0^Te^{\b
t}|y_1(t)-y_2(t)|^pdt+K\({1\over\b}\)^{p\over q}\dbE\int_0^Te^{\b
s}\(\int_0^s|z_1(s,t)-z_2(s,t)|^2dt\)^{p\over2}ds\\
\ns\ds\le\({K\over\b}+K\({1\over\b}\)^{p\over q}\)
\dbE\int_0^Te^{\b t}|y_1(t)-y_2(t)|^pdt\\
\ns\ds\le\({K\over\b}+K\({1\over\b}\)^{p\over
q}\)\|(y_1(\cd),z_1(\cd\,,\cd))-(y_2(\cd),z_2(\cd\,,\cd))\|^p_{\cM^p_\b[0,T]}.\ea\ee
Since the constant $K>0$ appears in the right hand side of the above
is independent of $\b$, by choosing $\b>0$ large, we obtain that
$\Th$ is a contraction. Hence, there exists a unique fixed point
$(Y(\cd),Z(\cd\,,\cd))\in\cM^p[0,T]$, which is the unique adapted
M-solution of (\ref{MF-BSVIE}).

\ms

\it Step 2. \rm The adapted M-solution
$(Y(\cd),Z(\cd\,,\cd))\in\cM^q[0,T]$ if $\psi(\cd)\in
L^q_{\cF_T}(0,T;\dbR^n)$.

\ms

Let $\psi(\cd)\in L^q_{\cF_T}(0,T;\dbR^n)\subseteq
L^2_{\cF_T}(0,T;\dbR^n)$. According to Step 1, there exists a unique
adapted M-solution $(Y(\cd),Z(\cd\,,\cd))\in\cM^2[0,T]$. We want to
show that in the current case, $(Y(\cd),Z(\cd\,,\cd))$ is actually
in $\cM^q[0,T]$. To show this, for the obtained adapted M-solution
$(Y(\cd),Z(\cd\,,\cd))$, let us consider the following MF-BSVIE:
\bel{MF-BSVIE3}\ba{ll}
\ns\ds
\wt Y(t)=\psi(t)+\int_t^Tg(t,s,\wt Y(s),\wt Z(t,s),Z(s,t),\G(t,s,\wt Y(s),\wt Z(t,s),Z(s,t))ds\\
\ns\ds\qq\qq\qq-\int_t^T\wt Z(t,s)dW(s),\qq\qq t\in[0,T].\ea\ee
For any $(y(\cd),z(\cd\,,\cd))\in\cM^p[0,T]$, by Proposition 3.1
(with $p=q$), the following MF-BSVIE admits a unique adapted
M-solution $(\wt Y(\cd),\wt Z(\cd\,,\cd))\in\cM^q[0,T]$:
\bel{MF-BSVIE3s}\ba{ll}
\ns\ds
\wt Y(t)=\psi(t)+\int_t^Tg(t,s,y(s),\wt Z(t,s),Z(s,t),\G(t,s,y(s),\wt Z(t,s),Z(s,t))ds\\
\ns\ds\qq\qq\qq-\int_t^T\wt Z(t,s)dW(s),\qq\qq t\in[0,T].\ea\ee
Thus, if we define $\wt\Th(y(\cd),z(\cd\,,\cd))=(\wt Y(\cd),\wt
Z(\cd\,,\cd))$, then $\wt\Th:\cM^q[0,T]\to\cM^q[0,T]$. We now show
that $\wt\Th$ is a contraction on $\cM^q[0,T]$ (Compare that $\Th$
in Step 1 is a contraction on $\cM^2[0,T]$). To this end, let
$(y_i(\cd),z_i(\cd\,,\cd))\in\cM^q[0,T]$ and let
$$(\wt Y_i(\cd),\wt
Z_i(\cd\,,\cd))=\wt\Th(y_i(\cd),z_i(\cd\,,\cd)),\qq i=1,2.$$
Then by Proposition 3.1 (with $p=q$), we have
$$\ba{ll}
\ns\ds\dbE\[|\wt Y_1(t)-\wt Y_2(t)|^q+\(\int_t^T|\wt Z_1(t,s)-\wt Z_2(t,s)|^2ds\)^{q\over2}\]\\
\ns\ds\le K\dbE\(\int_t^T|g(t,s,y_1(s),\wt Z_1(t,s),Z(s,t),\G(t,s,y_1(s),\wt Z_1(t,s),Z(s,t)))\\
\ns\ds\qq-
g(t,s,y_2(s),\wt Z_1(t,s),Z(s,t),\G(t,s,y_2(s),\wt Z_1(t,s),Z(s,t)))|ds\)^q\\
\ns\ds\le K\dbE\[\int_t^T\(|y_1(s)-y_2(s)|
+|\G(t,s,y_1(s),\wt Z_1(t,s),Z(s,t))-\G(t,s,y_2(s),\wt Z_1(t,s),Z(s,t))|\)ds\]^q\\
\ns\ds\le
K\dbE\[\int_t^T\(|y_1(s)-y_2(s)|+\dbE|y_1(s)-y_2(s)|\)ds\]^q\le
K\dbE\int_t^T|y_1(s)-y_2(s)|^qds.\ea$$
Then
$$\ba{ll}
\ns\ds\dbE\[\int_0^Te^{\b t}|\wt Y_1(t)-\wt
Y_2(t)|^qdt+\int_0^Te^{\b t}
\(\int_t^T|\wt Z_1(t,s)-\wt Z_2(t,s)|^2ds\)^{q\over2}dt\]\\
\ns\ds\le K\dbE\int_0^Te^{\b
t}\int_t^T|y_1(s)-y_2(s)|^qdsdt=K\dbE\int_0^T\int_0^se^{\b
t}|y_1(s)-y_2(s)|^qdtds\\
\ns\ds\le{K\over\b}\int_0^Te^{\b t}|y_1(t)-y_2(t)|^qdt.\ea$$
Hence, $\wt\Th$ is a contraction on $\cM^q[0,T]$ (with the
equivalent norm). Hence, (\ref{MF-BSVIE3s}) admits a unique adapted
M-solution $(\wt Y(\cd),\wt
Z(\cd\,,\cd))\in\cM^q[0,T]\subseteq\cM^2[0,T]$. Then by the
uniqueness of adapted solutions in $\cM^2[0,T]$ of
(\ref{MF-BSVIE3s}), it is necessary that
$$(Y(\cd),Z(\cd\,,\cd))=(\wt Y(\cd),\wt
Z(\cd\,,\cd))\in\cM^q[0,T].$$

\ms

\it Step 3. \rm Some estimates.

\nobreak

According to Proposition 3.1, we have
$$\ba{ll}
\ns\ds\dbE\Big\{|Y(t)|^q+\(\int_t^T|Z(t,s)|^2ds\)^{q\over2}\Big\}\\
\ns\ds\le
K\dbE\Big\{|\psi(t)|^q+\(\int_t^T|g(t,s,Y(s),0,Z(s,t),\G(t,s,Y(s),0,Z(s,t))|ds\)^q\Big\}\\
\ns\ds\le
K\dbE\Big\{1+|\psi(t)|^q+\[\int_t^T\(|Y(s)|+|Z(s,t)|^{2\over q}
+|\G(t,s,Y(s),0,Z(s,t))|\)ds\]^q\Big\}\\
\ns\ds\le
K\dbE\Big\{1+|\psi(t)|^q+\[\int_t^T\(|Y(s)|+|Z(s,t)|^{2\over q}\)ds\]^q\Big\}\\
\ns\ds\le
K\dbE\Big\{1+|\psi(t)|^q+\int_t^T|Y(s)|^qds+\(\int_t^T|Z(s,t)|^{2\over
q}ds\)^q\Big\}.\ea$$
Then
$$\ba{ll}
\ns\ds\dbE\Big\{\int_0^Te^{\b t}|Y(t)|^qdt
+\int_0^Te^{\b t}\(\int_t^T|Z(t,s)|^2ds\)^{q\over2}dt\Big\}\\
\ns\ds\le K\dbE\Big\{\int_0^Te^{\b
t}\(1+|\psi(t)|^q\)dt+\int_0^Te^{\b
t}\int_t^T|Y(s)|^qdsdt+\int_0^Te^{\b t}\(\int_t^T|Z(s,t)|^{2\over
q}ds\)^qdt\Big\}.\ea$$
Note that
$$\ba{ll}
\ns\ds\dbE\int_0^Te^{\b t}\(\int_t^T|Z(s,t)|^{2\over
q}ds\)^qdt=\dbE\int_0^Te^{\b t}\(\int_t^Te^{-{\b\over
q}s}e^{{\b\over
q}s}|Z(s,t)|^{2\over q}ds\)^qdt\\
\ns\ds\le\dbE\int_0^Te^{\b t}\(\int_t^Te^{-{\b\over
q-1}s}ds\)^{q-1}\(\int_t^Te^{\b s}|Z(s,t)|^2ds\)dt\\
\ns\ds=\dbE\int_0^Te^{\b t}\[{q-1\over\b}\(e^{-\b t\over
q-1}-e^{-{\b
T\over q-1}}\)\]^{q-1}\(\int_t^Te^{\b s}|Z(s,t)|^2ds\)dt\\
\ns\ds\le{(q-1)^{q-1}\over\b^{q-1}}\dbE\int_0^T\int_0^se^{\b
s}|Z(s,t)|^2dtds\le{(q-1)^{q-1}\over\b^{q-1}}\dbE\int_0^Te^{\b t
}|Y(t)|^2dt.\ea$$
Hence,
$$\ba{ll}
\ns\ds\dbE\Big\{\int_0^Te^{\b t}|Y(t)|^qdt
+\int_0^Te^{\b t}\(\int_t^T|Z(t,s)|^2ds\)^{q\over2}dt\Big\}\\
\ns\ds\le K\dbE\Big\{\int_0^Te^{\b
t}\(1+|\psi(t)|^q\)dt+{1\over\b}\int_0^Te^{\b
t}|Y(t)|^qdt+{1\over\b^{q-1}}\int_0^Te^{\b t}|Y(t)|^2dt\Big\}.\ea$$
By choosing $\b>0$ large, we obtain
$$\ba{ll}
\ns\ds\dbE\Big\{\int_0^Te^{\b t}|Y(t)|^qdt +\int_0^Te^{\b
t}\(\int_t^T|Z(t,s)|^2ds\)^{q\over2}dt\Big\}\le
K\dbE\(1+\int_0^Te^{\b t}|\psi(t)|^qdt\).\ea$$
Then (\ref{Lp-estimate}) follows.

\ms

Finally, let $\psi_i(\cd)\in L^q_{\cF_T}(0,T;\dbR^n)$, $g_i(\cd)$
and $\th_i(\cd)$ satisfy (H3)$_q$. Observe that
$$\ba{ll}
\ns\ds
Y_1(t)=\psi_1(t)+\int_t^Tg_1(t,s,Y_1(s),Z_1(t,s),Z_1(s,t),\G_1(t,s,Y_1(s),Z_1(t,s),Z_1(s,t)))ds\\
\ns\ds\qq\qq\qq-\int_t^TZ_1(t,s)dW(s)\\
\ns\ds=\psi_1(t)+\int_t^Tg_1(t,s,Y_1(s),Z_1(t,s),Z_1(s,t),\G_1(t,s,Y_1(s),Z_1(t,s),Z_1(s,t)))ds\\
\ns\ds\qq\qq-\int_t^Tg_2(t,s,Y_1(s),Z_1(t,s),Z_1(s,t),\G_2(t,s,Y_1(s),Z_1(t,s),Z_1(s,t)))ds\\
\ns\ds\qq\qq+\int_t^Tg_2(t,s,Y_1(s),Z_1(t,s),Z_1(s,t),\G_2(t,s,Y_1(s),Z_1(t,s),Z_1(s,t)))ds\\
\ns\ds\qq\qq-\int_t^Tg_2(t,s,Y_2(s),Z_1(t,s),Z_2(s,t),\G_2(t,s,Y_2(s),Z_1(t,s),Z_2(s,t)))ds\\
\ns\ds\qq\qq+\int_t^Tg_2(t,s,Y_2(s),Z_1(t,s),Z_2(s,t),\G_2(t,s,Y_2(s),Z_1(t,s),Z_2(s,t)))ds\\
\ns\ds\qq\qq\qq-\int_t^TZ_1(t,s)dW(s)\\
\ns\ds\equiv\h\psi_1(t)+\int_t^Tg_2(t,s,Y_2(s),Z_1(t,s),Z_2(s,t),\G_2(t,s,Y_2(s),Z_1(t,s),Z_2(s,t)))ds\\
\ns\ds\qq\qq\qq-\int_t^TZ_1(t,s)dW(s),\ea$$
with $\h\psi_1(\cd)$ defined in an obvious way. Then by Proposition
3.1 (with $p=2$), we obtain
$$\ba{ll}
\ns\ds\dbE\[|Y_1(t)-Y_2(t)|^2+\int_t^T|Z_1(t,s)-Z_2(t,s)|^2ds\]\le K\dbE|\h\psi_1(t)-\psi_2(t)|^2\\
\ns\ds\le K\dbE\Big\{|\psi_1(t)-\psi_2(t)|^2\\
\ns\ds\qq+\(\int_t^T|g_1(t,s,Y_1(s),Z_1(t,s),Z_1(s,t),\G_1(t,s,Y_1(s),Z_1(t,s),Z_1(s,t)))\\
\ns\ds\qq\qq-g_2(t,s,Y_1(s),Z_1(t,s),Z_1(s,t),\G_2(t,s,Y_1(s),Z_1(t,s),Z_1(s,t)))|ds\)^2\\
\ns\ds\qq+\(\int_t^T|g_2(t,s,Y_1(s),Z_1(t,s),Z_1(s,t),\G_2(t,s,Y_1(s),Z_1(t,s),Z_1(s,t)))\\
\ns\ds\qq\qq-g_2(t,s,Y_2(s),Z_1(t,s),Z_2(s,t),\G_2(t,s,Y_2(s),Z_1(t,s),Z_2(s,t)))|ds\)^2\Big\}\\
\ns\ds\le K\dbE\Big\{|\psi_1(t)-\psi_2(t)|^2+\(\int_t^T|(g_1-g_2)(t,s)|ds\)^2\\
\ns\ds\qq+\[\int_t^T\(|Y_1(s)-Y_2(s)|+|Z_1(s,t)-Z_2(s,t)|\)ds\]^2\Big\}.\ea$$
Then similar to the proof of the contraction for $\Th$, we can
obtain our stability estimate (\ref{stability}). \endpf

\ms

Let us make some remarks on the above result, together with its
proof.

\ms

First of all, we have seen that the growth of the maps
\bel{3.17}\hat z\mapsto g(t,s,y,z,\hat z),\q(\hat z,\hat
z')\mapsto\th(t,s,y,z,\hat z,y',z',\hat z')\ee
plays an important role in proving the well-posedness of MF-BSVIEs,
especially for the case of $p>2$. When $p\in(1,2]$, the adapted
M-solutions for BSVIEs was discussed in \cite{Wang 2011}. It is
possible to adopt the idea of \cite{Wang 2011} to treat MF-BSVIEs
for $p\in(1,2)$. If (H3)$_\infty$ holds, then for any $p>1$, as long
as $\psi(\cd)\in L^p_{\cF_T}(0,T;\dbR^n)$, (\ref{MF-BSVIE2}) admits
a unique adapted M-solution $(Y(\cd),Z(\cd\,,\cd))\in\cM^p[0,T]$. On
the other hand, if the maps in (\ref{3.17}) grow linearly, the
adapted M-solution $(Y(\cd),Z(\cd\,,\cd))$ of (\ref{MF-BSVIE2}) may
not be in $\cM^p[0,T]$ for $p>2$, even if $\psi(\cd)\in
L^p_{\cF_T}(0,T;\dbR^n)$. This can be seen from the following
example.

\ms

\bf Example 3.3. \rm Consider BSVIE:
\bel{3.20}Y(t)=\psi(t)+\int_t^TZ(s,t)ds-\int_t^TZ(t,s)dW(s),\q
t\in[0,T].\ee
Let
$$\psi(t) \equiv \int_0^T\psi_1(s)dW(s), \qq \forall t\in[0,T],$$
with $\psi_1(\cd)$ being deterministic and
$$\psi_1(\cd)\in
L^2(0,T)\setminus\bigcup_{p>2}L^p(0,T-\d;\dbR),$$
for some fixed $\d\in(0,T;\dbR)$. Thus, for any $p>1$,
$$\ba{ll}
\ns\ds\dbE\int_0^T|\psi(t)|^pds=T\dbE\Big|\int_0^T\psi_1(s)dW(s)\Big|^p\le
C\(\int_0^T|\psi_1(s)|^2ds\)^{p\over2},\ea$$
which means $\psi(\cd)\in L^p_{\cF_T}(0,T;\dbR)$ for any $p>1$. If
we define
$$\left\{\ba{ll}
\ns\ds Y(t)=\int_0^t\psi_1(s)dW(s)+\psi_1(t)(T-t),\qq
t\in[0,T],\\
\ns\ds Z(t,s)=\psi_1(s),\qq(s,t)\in[0,T]^2,\ea\right.$$
then
$$\ba{ll}
\ns\ds Y(t)=\int_0^t\psi_1(s)dW(s)+\psi_1(t)(T-t)\\
\ns\ds\qq=\psi(t)-\int_t^T\psi_1(s)dW(s)+\int_t^T\psi_1(t)ds\\
\ns\ds\qq=\psi(t)-\int_t^TZ(s,t)ds+\int_t^TZ(t,s)dW(s).\ea$$
This means that $(Y(\cd),Z(\cd\,,\cd))$ is the adapted M-solution of
(\ref{3.20}). We claim that $Y(\cd)\notin L^p_\dbF(0,T;\dbR)$, for
any $p>2$. In fact, if $Y(\cd)\in L^p_\dbF(0,T;\dbR)$ for some
$p>2$, then
$$\ba{ll}
\ns\ds\d^p\dbE\int_0^{T-\d}|\psi_1(t)|^pdt\le\int_0^T(T-t)^p|\psi_1(t)|^pdt\le2^{p-1}\dbE\int_0^T\(|Y(t)|^p
+\Big|\int_0^t\psi_1(s)dW(s)\Big|^p\)dt\\
\ns\ds\qq\qq\qq\qq\q\le
K\Big\{\int_0^T\dbE|Y(t)|^pdt+\(\int_0^T|\psi(s)|^2ds\)^{p\over2}\)<\infty.\ea$$
This is a contradiction.

\ms

Next, if
$$g_i(t,s,y,z,\hat z)=g_i(t,s,y,z),\qq\th_i(t,s,y,z,\hat
z,y',z',\hat z')=\th_i(t,s,y,z,y',z'),$$
then the stability estimate (\ref{stability}) can be improved to
\bel{stability2}\ba{ll}
\ns\ds\|(Y_1(\cd),Z_1(\cd\,,\cd))-(Y_2(\cd),Z_2(\cd\,,\cd))\|^q_{\cM^q[0,T]}\\
\ns\ds\le
K\dbE\Big\{\int_0^T|\psi_1(t)-\psi_2(t)|^qdt+\int_0^T\(\int_t^T|(g_1-g_2)(t,s)|ds\)^qdt\Big\},\ea\ee
for any $q>2$.

\ms

We point out that even for the special case of BSVIEs, the proof we
provided here significantly simplifies that given in \cite{Yong
2008}. The key is that we have a better understanding of the term
$Z(s,t)$ in the drift, and find a new way to treat it (see
(\ref{3.15})).

\ms

Now, let us look at linear MF-BSVIE (\ref{L-MF-BSVIE}). It is not
hard to see that under (L2), we have (H3)$_q$ with $q=2$. Hence, we
have the following corollary.

\ms

\bf Corollary 3.4. \rm Let {\rm(L2)} hold. Then for any
$\psi(\cd)\in L^2_{\cF_T}(0,T;\dbR^n)$, (\ref{L-MF-BSVIE}) admits a
unique adapted M-solution $(Y(\cd),Z(\cd\,,\cd))\in\cM^2[0,T]$.

\ms

\section{Duality Principles.}

In this section, we are going to establish two duality principles
between linear MF-FSVIEs and linear MF-BSVIEs. Let us first consider
the following linear MF-FSVIE (\ref{L-MF-FSVIE}) which is rewritten
below (for convenience):
\bel{4.1}\ba{ll}
\ns\ds X(t)=\f(t)+\int_0^t\(A_0(t,s)X(s)+\dbE'\[C_0(t,s)X(s)\]\)ds\\
\ns\ds\qq\qq\q+\int_0^t\(A_1(t,s)X(s)+\dbE'\[C_1(t,s)X(s)\]\)dW(s),\qq
t\in[0,T].\ea\ee
Let (L1) hold and $\f(\cd)\in L^2_\dbF(0,T;\dbR^n)$. Then by
Corollary 2.7, (\ref{4.1}) admits a unique solution $X(\cd)\in
L^2_\dbF(0,T;\dbR^n)$. Now, let $(Y(\cd),Z(\cd\,,\cd))\in\cM^2[0,T]$
be undetermined, and we observe the following:
$$\ba{ll}
\ns\ds\dbE\int_0^T\lan Y(t),\f(t)\ran dt=\dbE\int_0^T\lan
Y(t),X(t)-\int_0^t\(A_0(t,s)X(s)+\dbE'[C_0(t,s)X(s)]\)ds\ran dt\\
\ns\ds\qq\qq-\dbE\int_0^T\lan
Y(t),\int_0^t\(A_1(t,s)X(s)+\dbE'[C_1(t,s)X(s)]\)dW(s)\ran dt\\
\ns\ds\qq\qq\equiv\dbE\int_0^T\lan Y(t),X(t)\ran
dt-\sum_{i=1}^4I_i.\ea$$
We now look at each term $I_i$. First, for $I_1$, we have
$$I_1=\dbE\int_0^T\int_0^t\lan Y(t),A_0(t,s)X(s)\ran
dsdt=\dbE\int_0^T\lan X(t),\int_t^TA_0(s,t)^TY(s)ds\ran dt.$$
Next, for $I_2$, let us pay some extra attention on $\o$ and $\o'$,
$$\ba{ll}
\ns\ds I_2=\dbE\int_0^T\3n\int_0^t\2n\lan
Y(t),\dbE'[C_0(t,s)X(s)]\ran
dsdt=\dbE'\dbE\int_0^T\3n\int_s^T\2n\lan
C_0(t,s,\o,\o')^TY(t,\o),X(s,\o')\ran dtds\\
\ns\ds\q=\dbE\dbE^*\3n\int_0^T\3n\int_s^T\3n\lan
C_0(t,s,\o^*,\o)^TY(t,\o^*),X(s,\o)\ran
dtds\1n=\1n\dbE\2n\int_0^T\3n\lan
X(t),\int_t^T\2n\dbE^*[C_0(s,t)^TY(s)]ds\ran dt.\ea$$
Here, we have introduced the notation $\dbE^*$, whose definition is
obvious from the above, to distinguish $\dbE$ (and $\dbE'$). For
$I_3$, we have
$$\ba{ll}
\ns\ds I_3=\dbE\int_0^T\lan Y(t),\int_0^tA_1(t,s)X(s)dW(s)\ran
dt\\
\ns\ds\q=\dbE\int_0^T\lan\dbE
Y(t)+\int_0^tZ(t,s)dW(s),\int_0^tA_1(t,s)X(s)dW(s)\ran dt\\
\ns\ds\q=\dbE\int_0^T\int_0^t\lan Z(t,s),A_1(t,s)X(s)\ran
dsdt=\dbE\int_0^T\lan X(t),\int_t^TA_1(s,t)^TZ(s,t)ds\ran dt.\ea$$
Finally, we look at $I_4$.
$$\ba{ll}
\ns\ds I_4=\dbE\int_0^T\lan
Y(t),\int_0^t\dbE'[C_1(t,s)X(s)]dW(s)\ran dt\\
\ns\ds\q=\dbE\int_0^T\int_0^t\lan Z(t,s),\dbE'[C_1(t,s)X(s)]\ran dsdt\\
\ns\ds\q=\dbE'\dbE\int_0^T\int_s^T\lan
Z(t,s,\o),C_1(t,s,\o,\o')X(s,\o')\ran dtds\\
\ns\ds\q=\dbE\int_0^T\lan
X(t),\int_t^T\dbE^*[C_1(s,t)^TZ(s,t)]ds\ran dt.\ea$$
Hence, we obtain
$$\ba{ll}
\ns\ds\dbE\int_0^T\lan Y(t),\f(t)\ran dt=\dbE\int_0^T\lan
X(t),Y(t)-\int_t^T\(A_0(s,t)^TY(s)+A_1(s,t)^TZ(s,t)\\
\ns\ds\qq\qq\qq\qq\qq\qq+\dbE^*\[C_0(s,t)^TY(s)+C_1(s,t)^TZ(s,t)\]\)ds\ran
dt.\ea$$
On the other hand, suppose (L1)$'$ holds and $\f(\cd)\in
C^p_\dbF([0,T];\dbR^n)$. Then $X(\cd)\in C^p_\dbF([0,T];\dbR^n)$.

Consequently, we obtain the following {\it duality principle for
MF-FSVIEs} whose proof is clear from the above.

\ms

\bf Theorem 4.1. \sl Let {\rm(L1)} hold, and $\f(\cd),\psi(\cd)\in
L^2_\dbF(0,T;\dbR^n)$. Let $X(\cd)\in L^2_\dbF(0,T;\dbR^n)$ be the
solution to the linear MF-FSVIE $(\ref{4.1})$, and
$(Y(\cd),Z(\cd\,,\cd))\in\cM^2[0,T]$ be the adapted M-solution to
the following linear MF-BSVIE:
\bel{4.5}\ba{ll}
\ns\ds Y(t)=\psi(t)+\int_t^T\(A_0(s,t)^TY(s)+A_1(s,t)^TZ(s,t)\\
\ns\ds\qq\qq\qq\qq+\dbE^*\[C_0(s,t)^TY(s)+C_1(s,t)^TZ(s,t)\]\)ds-\int_t^TZ(t,s)dW(s).\ea\ee
Then
\bel{duality2}\dbE\int_0^T\lan X(t),\psi(t)\ran dt=\dbE\int_0^T\lan
Y(t),\f(t)\ran dt.\ee

\ms

\rm

We call (\ref{4.5}) the adjoint equation of (\ref{4.1}). The above
duality principle will be used in establishing Pontryagin's type
maximum principle for optimal controls of MF-FSVIEs.

\ms

Next, different from the above, we want to start from the followng
linear MF-BSVIE:
\bel{L-MF-BSVIE1}\ba{ll}
\ns\ds
Y(t)=\psi(t)+\int_t^T\(\bar A_0(t,s)Y(s)+\bar C_0(t,s)Z(s,t)+\dbE'[\bar A_1(t,s)Y(s)+\bar C_1(t,s)Z(s,t)]\)ds\\
\ns\ds\qq\qq\qq\qq-\int_t^TZ(t,s)dW(s),\qq\qq t\in[0,T].\ea\ee
This is a special case of (\ref{L-MF-BSVIE}) in which
$$\bar B_0(t,s)=0,\qq\bar B_1(t,s)=0.$$
Under (L2), by Corollary 3.4, for any $\psi(\cd)\in
L^2_\dbF(0,T;\dbR^n)$, (\ref{L-MF-BSVIE1}) admits a unique adapted
M-solution $(Y(\cd),Z(\cd\,,\cd))\in\cM^2[0,T]$. We point out here
that for each $t\in[0,T)$, the maps
$$s\mapsto\bar C_0(t,s),\q s\mapsto\bar C_1(t,s)$$
are $\dbF$-progressively measurable and $\dbF^2$-progressively
measurable on $[t,T]$, respectively. Now, we let a process
$X(\cd)\in L^2_\dbF(0,T;\dbR^n)$ be undetermined, and make the
following calculation:
$$\ba{ll}
\ns\ds\dbE\int_0^T\lan X(t),\psi(t)\ran dt=\dbE\int_0^T\lan
X(t),Y(t)-\int_t^T\(\bar A_0(t,s)Y(s)+\bar C_0(t,s)Z(s,t)\\
\ns\ds\qq\q+\dbE'[\bar A_1(t,s)Y(s)+\bar
C_1(t,s)Z(s,t)]\)ds-\int_t^TZ(t,s)dW(s)\ran
dt\\
\ns\ds=\dbE\int_0^T\lan X(t),Y(t)\ran dt-\dbE\int_0^T\int_t^T\lan
X(t),\bar A_0(t,s)Y(s)\ran dsdt\\
\ns\ds\qq-\dbE\int_0^T\int_t^T\lan X(t),\bar C_0(t,s)Z(s,t)\ran
dsdt-\dbE\int_0^T\int_t^T\lan X(t),\dbE'[\bar A_1(t,s)Y(s)]\ran
dsdt\\
\ns\ds\qq-\dbE\int_0^T\int_t^T\lan X(t),\dbE'[\bar
C_1(t,s)Z(s,t)]\ran dsdt\equiv\dbE\int_0^T\lan X(t),Y(t)\ran
dt-\sum_{i=1}^4I_i.\ea$$
Similar to the above, we now look at the terms $I_i$ ($i=1,2,3,4$)
one by one. First, we look at $I_1$:
$$\ba{ll}
\ns\ds I_1=\dbE\int_0^T\int_t^T\lan X(t),\bar A_0(t,s)Y(s)\ran
dsdt=\dbE\int_0^T\int_0^s\lan\bar A_0(t,s)^TX(t),Y(s)\ran dtds\\
\ns\ds\qq=\dbE\int_0^T\lan\int_0^t\bar A_0(s,t)^TX(s)ds,Y(t)\ran
dt.\ea$$
Next, for $I_2$, one has
$$\ba{ll}
\ns\ds I_2=\dbE\int_0^T\int_t^T\lan X(t),\bar C_0(t,s)Z(s,t)\ran
dsdt=\dbE\int_0^T\int_0^s\lan\bar C_0(t,s)^TX(t),Z(s,t)\ran dtds\\
\ns\ds\q=\int_0^T\dbE\int_0^t\lan\bar C_0(s,t)^TX(s),Z(t,s)\ran
dsdt\\
\ns\ds\q=\int_0^T\dbE\lan\int_0^t\dbE[\bar
C_0(s,t)^T\bigm|\cF_s]X(s)dW(s),\int_0^tZ(t,s)dW(s)\ran
dt\\
\ns\ds\q=\int_0^T\dbE\lan\int_0^t\dbE[\bar
C_0(s,t)^T\bigm|\cF_s]X(s)dW(s),Y(t)-\dbE
Y(t)\ran dt\\
\ns\ds\q=\dbE\int_0^T\lan\int_0^t\dbE[\bar
C_0(s,t)^T\bigm|\cF_s]X(s)dW(s),Y(t)\ran dt.\ea$$
Now, for $I_3$,
$$\ba{ll}
\ns\ds I_3=\dbE\int_0^T\int_t^T\lan X(t),\dbE'[\bar
A_1(t,s)Y(s)]\ran
dsdt=\dbE\dbE'\int_0^T\int_0^s\lan\bar A_1(t,s,\o,\o')^TX(t,\o),Y(s,\o')\ran dtds\\
\ns\ds\q=\dbE'\int_0^T\lan\int_0^t\dbE[\bar
A_1(s,t,\o,\o')^TX(s,\o)]ds,Y(t,\o')\ran
dt\\
\ns\ds\q=\dbE\int_0^T\lan\int_0^t\dbE^*[\bar
A_1(s,t,\o^*,\o)^TX(s,\o^*)]ds,Y(t,\o)\ran
dt\\
\ns\ds\q\equiv\dbE\int_0^T\lan\int_0^t\dbE^*[\bar
A_1(s,t)^TX(s)]ds,Y(t)\ran dt.\ea$$
Finally, similar to the above, one has
$$\ba{ll}
\ns\ds I_4=\dbE\int_0^T\int_t^T\lan X(t),\dbE'[\bar
C_1(t,s)Z(s,t)]\ran dsdt\\
\ns\ds\q=\dbE\dbE'\int_0^T\int_0^s\lan\bar C_1(t,s,\o,\o')^TX(t,\o),Z(s,t,\o')\ran dtds\\
\ns\ds\q=\dbE'\int_0^T\int_0^t\lan \dbE[\bar
C_1(s,t,\o,\o')^TX(s,\o)],Z(t,s,\o')\ran
dsdt\\
\ns\ds\q=\dbE\int_0^T\int_0^t\lan
\dbE^*[\bar C_1(s,t)^TX(s)],Z(t,s)\ran dsdt\\
\ns\ds\q=\int_0^T\dbE\int_0^t\lan
\dbE\[\dbE^*[\bar C_1(s,t)^TX(s)]\bigm|\cF_s\],Z(t,s)\ran dsdt\\
\ns\ds\q=\int_0^T\dbE\lan\int_0^t
\dbE^*\[\dbE[\bar C_1(s,t)^T\bigm|\cF_s]X(s)\]dW(s),\int_0^tZ(t,s)dW(s)\ran dt\\
\ns\ds\q=\int_0^T\dbE\lan\int_0^t
\dbE^*\[\dbE[\bar C_1(s,t)^T\bigm|\cF_s]X(s)\]dW(s),Y(t)-\dbE Y(t)\ran dt\\
\ns\ds\q=\dbE\int_0^T\lan\int_0^t\dbE^*\[\dbE[\bar
C_1(s,t)^T\bigm|\cF_s]X(s)\]dW(s),Y(t)\ran dt.\ea$$
Combining the above, we obtain
\bel{4.2}\ba{ll}
\ns\ds\dbE\int_0^T\lan X(t),\psi(t)\ran dt=\dbE\int_0^T\lan
X(t),Y(t)\ran dt-\sum_{i=1}^4I_i\\
\ns\ds=\dbE\int_0^T\lan
Y(t),X(t)-\int_0^t\(\bar A_0(s,t)^TX(s)+\dbE^*[\bar A_1(s,t)^TX(s)]\)ds\\
\ns\ds\qq-\int_0^t\(\dbE[\bar
C_0(s,t)^T\bigm|\cF_s]X(s)+\dbE^*\[\dbE[\bar
C_1(s,t)^T\bigm|\cF_s]X(s)\]\)dW(s)\ran dt.\ea\ee
Now, we are at the position to state and prove the following {\it
duality principle for MF-BSVIEs}.

\ms

\bf Theorem 4.2. \sl Let {\rm(L2)} hold and $\psi(\cd)\in
L^2_{\cF_T}(0,T;\dbR^n)$. Let $(Y(\cd),Z(\cd\,,\cd))\in\cM^2[0,T]$
be the unique adapted M-solution of linear MF-BSVIE
$(\ref{L-MF-BSVIE1})$. Further, let $\f(\cd)\in
L_\dbF^2(0,T;\dbR^n)$ and $X(\cd)\in L^2_\dbF(0,T;\dbR^n)$ be the
solution to the following linear MF-FSVIE:
\bel{L-MF-FSVIE2}\ba{ll}
\ns\ds X(t)=\f(t)+\int_0^t\(\bar A_0(s,t)^TX(s)+\dbE^*[\bar A_1(s,t)^TX(s)]\)ds\\
\ns\ds\qq\qq+\int_0^t\(\dbE[\bar
C_0(s,t)^T\bigm|\cF_s]X(s)+\dbE^*\[\dbE[\bar
C_1(s,t)^T\bigm|\cF_s]X(s)\]\) dW(s),\qq t\in[0,T].\ea\ee
Then
\bel{duality}\dbE\int_0^T\lan Y(t),\f(t)\ran dt=\dbE\int_0^T\lan
X(t),\psi(t)\ran dt.\ee

\ms

\it Proof. \rm For linear MF-FSVIE (\ref{L-MF-FSVIE2}), when (L2)
holds, we have (L1). Hence, for any $\f(\cd)\in
L^2_\dbF(0,T;\dbR^n)$, (\ref{L-MF-FSVIE2}) admits a unique solution
$X(\cd)\in L^2_\dbF(0,T;\dbR^n)$. Then (\ref{duality}) follows from
(\ref{4.2}) immediately. \endpf

\ms

We call MF-FSVIE (\ref{L-MF-FSVIE2}) the adjoint equation of
MF-BSVIE (\ref{L-MF-BSVIE1}). Such a duality principle will be used
to establish comparison theorems for MF-BSVIEs. Note that since for
$s<t$, $\bar C_0(s,t)^T$ is $\cF_t$-measurable and not necessarily
$\cF_s$-measurable, we have
\bel{4.8}\dbE[\bar C_0(s,t)^T\bigm|\cF_s]\ne\bar C_0(s,t),\qq
t\in(s,T],\ee
in general. Likewise, in general,
\bel{4.9}\dbE[\bar C_1(s,t)^T\bigm|\cF_s]\ne\bar C_1(s,t),\qq
t\in(s,T].\ee

 \ms

\ms

We now make some comparison between Theorems 4.1 and 4.2.

First, we begin with linear MF-FSVIE (\ref{4.1}) which is rewritten
here for convenience:
\bel{4.7}\ba{ll}
\ns\ds X(t)=\f(t)+\int_0^t\(A_0(t,s)X(s)+\dbE'[C_0(t,s)X(s)]\)ds\\
\ns\ds\qq\qq\q+\int_0^t\(A_1(t,s)X(s)+\dbE'[C_1(t,s)X(s)]\)dW(s),\qq
t\in[0,T].\ea\ee
According to Theorem 4.1, the adjoint equation of (\ref{4.7}) is
MF-BSVIE (\ref{4.5}). Now, we want to use Theorem 4.2 to find the
adjoint equation of (\ref{4.5}) which is regarded as
(\ref{L-MF-BSVIE1}) with
$$\left\{\ba{ll}
\ns\ds\bar A_0(t,s)=A_0(s,t)^T,\qq\bar A_1(t,s,\o,\o')=C_0(s,t,\o',\o)^T,\\
\ns\ds\bar C_0(t,s)=A_1(s,t)^T,\qq\bar
C_1(t,s,\o,\o')=C_1(s,t,\o',\o)^T.\ea\right.$$
Then, by Theorem 4.2, we obtain the adjoint equation
(\ref{L-MF-FSVIE2}) with the coefficients:
$$\left\{\ba{ll}
\ns\ds\bar A_0(s,t)^T=A_0(t,s),\qq\bar
A_1(s,t,\o',\o)^T=C_0(t,s,\o,\o'),\\
\ns\ds\dbE[\bar
C_0(s,t)^T\bigm|\cF_s]=\dbE[A_1(t,s)\bigm|\cF_s]=A_1(t,s),\\
\ns\ds\dbE[\bar
C_1(s,t,\o',\o)^T\bigm|\cF_s]=\dbE[C_1(t,s,\o,\o')\bigm|\cF_s]=C_1(t,s,\o,\o').\ea\right.$$
Hence, (\ref{4.7}) is the adjoint equation of (\ref{4.5}). Thus, we
have the following conclusion:
$$\hb{\sl Twice adjoint equation of a linear MF-FSVIE is itself.}$$

\ms

Next, we begin with linear MF-BSVIE (\ref{L-MF-BSVIE1}). From
Theorem 4.2, we know that the adjoint equation is linear MF-FSVIE
(\ref{L-MF-FSVIE2}). Now, we want to use Theorem 4.1 to find the
adjoint equation of (\ref{L-MF-FSVIE2}) which is regarded as
(\ref{4.7}) with
$$\left\{\ba{ll}
\ns\ds A_0(t,s)=\bar A_0(s,t)^T,\qq C_0(t,s,\o,\o')=\bar A_1(s,t,\o',\o)^T,\\
\ns\ds A_1(t,s)=\dbE[\bar C_0(s,t)^T\bigm|\cF_s],\qq
C_1(t,s,\o,\o')=\dbE[\bar C_1(s,t,\o',\o)^T\bigm|\cF_s].\ea\right.$$
Then by Theorem 4.2, the adjoint equation is given by (\ref{4.5})
with coefficients:
$$\left\{\ba{ll}
\ns\ds A_0(s,t)^T=\bar A_0(t,s),\qq C_0(s,t,\o',\o)^T=\bar A_1(t,s,\o,\o'),\\
\ns\ds A_1(s,t)^T=\dbE[\bar C_0(t,s)\bigm|\cF_t],\qq
C_1(s,t,\o',\o)=\dbE[\bar C_1(t,s,\o,\o')\bigm|\cF_t].\ea\right.$$
In another word, the twice adjoint equation of linear MF-BSVIE
(\ref{L-MF-BSVIE1}) is the following:
\bel{}\ba{ll}
\ns\ds Y(t)=\psi(t)+\int_t^T\(\bar A_0(t,s)Y(s)+\dbE[\bar
C_0(t,s)\bigm|\cF_t]Z(s,t)\\
\ns\ds\qq\qq+\dbE'\[\bar A_1(t,s)Y(s)+\dbE[\bar
C_1(t,s)\bigm|\cF_t]Z(s,t)\]\)ds-\int_t^TZ(t,s)dW(s),\q
t\in[0,T],\ea\ee
which is different from (\ref{L-MF-BSVIE1}), unless $\bar C_0(t,s)$
and $\bar C_1(t,s)$ are $\cF_t$-measurable for all $(t,s)\in\D$.
Thus, we have the following conclusion:
$$\hb{\sl Twice adjoint of a linear MF-BSVIE is not necessarily itself.}$$

\section{Comparison Theorems.}

\rm

In this section, we are going to establish some comparison theorems
for MF-FSVIEs and MF-BSVIEs, allowing the dimension to be larger
than 1. Let
$$\dbR^n_+=\Big\{(x_1,\cds,x_n)\in\dbR^n\bigm|x_i\ge0,~1\le i\le n\Big\}.$$
When $x\in\dbR^n_+$, we also denote it by $x\ge0$, and say that $x$
is {\it nonnegative}. By $x\le0$ and $x\ge y$ (if $x,y\in\dbR^n$),
we mean $-x\ge0$ and $x-y\ge0$, respectively. Moreover, if $X(\cd)$
is a process, then by $X(\cd)\ge0$, we mean
$$X(t)\ge0,\qq t\in[0,T],\q\as$$
Also, $X(\cd)$ is said to be {\it nondecreasing} if it is
componentwise nondecreasing. Likewise, we may define $X(\cd)\le0$
and $X(\cd)\ge Y(\cd)$ (if both $X(\cd)$ and $Y(\cd)$ are
$\dbR^n$-valued processes), and so on.

\ms

In what follows, we let $e_i\in\dbR^n$ be the vector that the $i$-th
entry is 1 and all other entries are zero. Also, we let
$$\left\{\ba{ll}
\ns\ds\dbM^n_+=\Big\{A=(a_{ij})\in\dbR^{n\times
n}\bigm|a_{ij}\ge0,~i\ne j\Big\}\equiv \Big\{A\in\dbR^{n\times
n}\bigm|\lan Ae_i,e_j\ran\ge0,~i\ne j\Big\},\\
\ns\ds\h\dbM^{n\times m}_+=\Big\{A=(a_{ij})\in\dbR^{n\times
m}\bigm|a_{ij}\ge0,~1\le i\le n,~1\le j\le m\Big\},\\
\ns\ds\dbM^n_0=\Big\{A=(a_{ij})\in\dbR^{n\times
n}\bigm|a_{ij}=0,~i\ne j\Big\}\equiv \Big\{A\in\dbR^{n\times
n}\bigm|\lan Ae_i,e_j\ran=0,~i\ne j\Big\}.\ea\right.$$
Note that $\h\dbM^{n\times m}_+$ is the set of all $(n\times m)$
matrices with all the entries being nonnegative, $\dbM^n_+$ is the
set of all $(n\times n)$ matrices with all the off-diagonal entries
being nonnegative, and $\dbM^n_0$ is actually the set of all
$(n\times n)$ diagonal matrices. Clearly, $\dbM^n_+$ and
$\h\dbM^{n\times m}_+$ are closed convex cones of $\dbR^{n\times n}$
and $\dbR^{n\times m}$, respectively, and $\dbM^n_0$ is a proper
subspace of $\dbR^{n\times n}$. Whereas, for $n=m=1$, one has
\bel{5.1}\dbM^1_+=\dbM^1_0=\dbR,\qq\h\dbM^{1\times1}_+=\dbR_+\equiv[0,\infty).\ee
We have the following simple result which will be useful below and
whose proof is obvious.

\ms

\bf Proposition 5.1. \sl Let $A\in\dbR^{n\times m}$. Then
$A\in\h\dbM^{n\times m}_+$ if and only if
\bel{}Ax\ge0,\qq\forall x\in\dbR^m_+.\ee

\rm

In what follows, we will denote $\h\dbM^n_+=\h\dbM^{n\times n}_+$.

\subsection{Comparison of solutions to MF-FSVIEs.}

In this subsection, we would like to discuss comparison of solutions
to linear MF-FSVIEs. There are some positive and also negative
results. To begin with, let us first present an example of MF-FSDEs.

\ms

\bf Example 5.2. \rm Consider the following one-dimensional linear
MF-FSDE, written in the integral form:
$$X(t)=1+\int_0^t\dbE X(s)dW(s),\qq t\in[0,T].$$
Taking expectation, we have
$$\dbE X(t)=1,\qq\forall t\in[0,T].$$
Consequently, the solution $X(\cd)$ is given by
$$X(t)=1+\int_0^tdW(s)=1+W(t),\qq t\in[0,T].$$
Thus, although $X(0)=1>0$, the following fails:
$$X(t)\ge0,\qq t\in[0,T],\q\as$$

The above example shows that if the diffusion contains a nonlocal
term in an MF-FSDE, we could not get an expected comparison of
solutions, in general. Therefore, for linear MF-FSDEs, one had
better only look at the following:
\bel{LSDE}\ba{ll}
\ns\ds
X(t)=x+\int_0^t\(A_0(s)X(s)+\dbE'[C_0(s)X(s)]\)ds+\int_0^tA_1(s)X(s)dW(s),\q
t\in[0,T],\ea\ee
with the diffusion does not contain a nonlocal term. For the above,
we make the following assumption.

\ms

{\bf(C1)} The maps
$$A_0,A_1:[0,T]\times\O\to\dbR^{n\times n},\q
C_0:[0,T]\times\O^2\to\dbR^{n\times n},$$
are uniformly bounded, and they are $\dbF$-progressively measurable,
and $\dbF^2$-progressively measurable, respectively.

\ms

Note that, due to (\ref{5.1}), the above (C1) is always true if
$n=1$. We now present the following comparison theorem for linear
MF-FSDEs.

\ms

\bf Proposition 5.3. \sl Let {\rm(C1)} hold and moreover,
\bel{5.4}A_0(s,\o)\in\dbM^n_+,\q C_0(s,\o,\o')\in\h\dbM^n_+,\q
A_1(s,\o)\in\dbM^n_0,\qq s\in[0,T],~~\as\o,\o'\in\O.\ee
Let $X(\cd)\in L^2_\dbF(0,T;\dbR^n)$ be the solution of linear
MF-FSDE $(\ref{LSDE})$ with $x\ge0$. Then
\bel{ge0}X(t)\ge0,\qq\forall t\in[0,T],\q\as\ee

\ms

\it Proof. \rm It is known from Theorem 2.6 that as a special case
of MF-FSVIE, the linear MF-FSDE (\ref{LSDE}) admits a unique
solution $X(\cd)\in L^p_\dbF(0,T;\dbR^n)$ for any $x\in\dbR^n$, and
any $p\ge2$. Further, it is not hard to see that $X(\cd)$ has
continuous paths. Since the equation is linear, it suffices to show
that $x\le0$ implies
\bel{le0}X(t)\le0,\qq t\in[0,T],\q\as\ee
To prove (\ref{le0}), we define a convex function
$$f(x)=\sum_{i=1}^n(x_i^+)^2,\qq\forall
x=(x_1,x_2,\cds,x_n)\in\dbR^n,$$
where $a^+=\max\{a,0\}$ for any $a\in\dbR$. Applying It\^o's formula
to $f(X(t))$, we get
$$\ba{ll}
\ns\ds f(X(t))-f(x)=\int_0^t\[\lan
f_x(X(s)),A_0(s)X(s)+\dbE'[C_0(s)X(s)]\ran\\
\ns\ds\qq+{1\over2}\lan
f_{xx}(X(s))A_1(s)X(s),A_1(s)X(s)\ran\]ds+\int_0^t\lan
f_x(X(s)),A_1(s)X(s)\ran dW(s).\ea$$
We observe the following: (noting $A_0(s)\in\dbM^n_+$)
$$\ba{ll}
\ns\ds\lan f_x(X(s)),A_0(s)X(s)\ran=\sum_{i,j=1}^n2X_i(s)^+\lan
e_i,A_0(s)e_j\ran X_j(s)\\
\ns\ds\qq=\sum_{i=1}^n2X_i(s)^+\lan e_i,A_0(s)e_i\ran
X_i(s)+\sum_{i\ne j}2X_i(s)^+\lan
e_i,A_0(s)e_j\ran X_j(s)\\
\ns\ds\qq\le\sum_{i=1}^n2[X_i(s)^+]^2\lan
e_i,A_0(s)e_i\ran+\sum_{i\ne j}2\lan e_i,A_0(s)e_j\ran X_i(s)^+
X_j(s)^+\le Kf(X(s)).\ea$$
Also, one has (making use of $C_0(s)\in\h\dbM^n_+$)
$$\ba{ll}
\ns\ds\dbE\lan
f_x(X(s)),\dbE'[C_0(s)X(s)]\ran\\
\ns\ds=2\int_{\O^2}\sum_{i,j=1}^nX_i(s,\o)^+\lan
e_i,C_0(s,\o,\o')e_j\ran X_j(s,\o')\dbP(d\o)\dbP(d\o')\\
\ns\ds\le2\int_{\O^2}\sum_{i,j=1}^nX_i(s,\o)^+\lan
e_i,C_0(s,\o,\o')e_j\ran X_j(s,\o')^+\dbP(d\o)\dbP(d\o')\\
\ns\ds\le K\(\dbE\[\sum_{i=1}^nX_i(s)^+\]\)^2\le K\dbE f(X(s)).\ea$$
Next, we have (noting $A_1(\cd)$ and $f_{xx}(\cd)$ are diagonal)
$$\ba{ll}
\ns\ds{1\over2}\dbE\lan
f_{xx}(X(s))A_1(s)X(s),A_1(s)X(s)\ran={1\over2}\dbE\sum_{i=1}^nI_{(X_i(s)\ge0)}\(\lan
A_1(s)e_i,e_i\ran X_i(s)\)^2\\
\ns\ds={1\over2}\dbE\sum_{i=1}^n\lan
A_1(s)e_i,e_i\ran{}^2[X_i(s)^+]^2\le Kf(X(s)).\ea$$
Consequently,
$$\dbE f(X(t))\le f(x)+K\int_0^t\dbE f(X(s))ds,\qq t\in[0,T].$$
Hence, by Gronwall's inequality, we obtain
$$\sum_{i=1}^n\dbE|X_i(t)^+|^2\le K\sum_{i=1}^n|x_i^+|^2,\qq t\in[0,T].$$
Therefore, if $x\le0$ (component-wise), then
$$\sum_{i=1}^n\dbE|X_i(t)^+|^2=0,\qq\forall t\in[0,T].$$
This leads to (\ref{le0}). \endpf

\ms

We now make some observations on condition (\ref{5.4}).

\ms

\it 1. Let $C_0(\cd)=0$, $A_1(\cd)=0$, and $A_0(\cd)$ be continuous
and for some $i\ne j$,
$$\lan A_0(0)e_i,e_j\ran<0,$$
\rm i.e., at least one off-diagonal entry of $A_0(0)$ is negative.
Then by letting $x=e_i$, we have
$$X_j(t)=\lan X(t),e_j\ran=\int_0^t\lan A_0(s)X(s),e_j\ran ds=\lan A_0(0)e_i,e_j\ran t+o(t)<0,$$
for $t>0$ small. Thus, $X(0)\ge0$ does not imply $X(t)\ge0$.

\ms

\it 2. Let $A_0(\cd)=0$, $A_1(\cd)=0$, and $C_0(\cd)$ be continuous
and for some $i\ne j$,
$$\lan C_0(0)e_i,e_j\ran<0,$$
\rm i.e., at least one off-diagonal entry of $C_0(0)$ is negative.
Then by a similar argument as above, we have that $X(0)\ge0$ does
not imply $X(t)\ge0$.

\ms

\it 3. Let $A_0(\cd)=0$, $C_0(\cd)=0$ and for some $i\ne j$,
$$\int_0^T\dbP\(\lan A_1(s)e_i,e_j\ran\ne0\)ds>0,$$
\rm i.e., at least one off-diagonal entry of $A_1(\cd)$ is not
identically zero. Then by letting $x=e_i$, we have
$$X_j(t)=\int_0^t\lan A_1(s)X(s),e_j\ran dW(s)\not\equiv0,\qq t\in[0,T].$$
On the other hand,
$$\dbE X_j(t)=0,\qq t\in[0,T].$$
Hence,
$$X_j(t)\ge0,\qq\forall t\in[0,T],~\as$$
must fail.

\ms

\it 4. Let $n=1$, $A_0(\cd)=A_1(\cd)=0$ and $C_0(\cd)$ bounded,
$\dbF$-adapted with
$$C_0(s)\ne0,\q\dbE C_0(s)=0,\qq s\in[0,T].$$
\rm This means that ``$C_0(s)\ge0,~\forall s\in[0,T],~\as$'' fails
(or diagonal elements of $C_0(\cd)$ are not all nonnegative).
Consider the following MF-FSDE:
$$X(t)=1+\int_0^tC_0(s)\dbE X(s)ds,\qq t\in[0,T].$$
Then
$$\dbE X(t)=1,\qq t\in[0,T].$$
Hence,
$$X(t)=1+\int_0^tC_0(s)ds,\qq t\in[0,T].$$
It is easy to choose a $C_0(\cd)$ such that
$$X(t)\ge0,\qq\forall t\in[0,T],~\as$$
is violated.

\ms

The above observations show that, in some sense, conditions assumed
in (\ref{5.4}) are sharp for Proposition 5.3.

\ms

Based on the above, let us now consider the following linear
MF-FSVIE:
\bel{5.7}\ba{ll}
\ns\ds X(t)=\f(t)+\int_0^t\(A_0(t,s)X(s)+\dbE'\[C_0(t,s)X(s)\]\)ds\\
\ns\ds\qq\qq\qq\qq+\int_0^tA_1(s)X(s)dW(s),\qq\qq t\in[0,T].\ea\ee
Note that $A_1(\cd)$ is independent of $t$ here. According to
\cite{Tudor 1989}, we know that for (linear) FSVIEs (without the
nonlocal term, i.e., $C_0(\cd\,,\cd)=0$ in (\ref{5.7})), if the
diffusion depends on both $(t,s)$ and $X(\cd)$, i.e., $A_1(t,s)$
really depends on $(t,s)$, a comparison theorem will fail in
general. Next, let us look at an example which is concerned with the
free term $\f(\cd)$.

\ms

\bf Example 5.4. \rm Consider the following one-dimensional FSVIE:
$$X(t)=T-t+\int_0^tbX(s)ds+\int_0^t\si X(s)dW(s),\qq t\in[0,T],$$
for some $b,\si\in\dbR$. The above is equivalent to the following:
$$\left\{\ba{ll}
\ns\ds dX(t)=[bX(t)-1]dt+\si X(t)dW(t),\qq t\in[0,T],\\
\ns\ds X(0)=T.\ea\right.$$
The solution to the above is explicitly given by the following:
$$X(t)=e^{(b-{\si^2\over2})t+\si
W(t)}\[T-\int_0^te^{-(b-{\si^2\over2})s-\si W(s)}ds\],\qq
t\in[0,T].$$
We know that as long as $\si\ne0$, for any $t>0$ small and any
$K>0$,
$$\dbP\(\int_0^te^{-(b-{\si^2\over2})s-\si W(s)}ds\ge K\)>0.$$
Therefore, we must have
$$\dbP(X(t)<0)>0,\qq\forall t>0~\hb{(small)}.$$
On the other hand, if $\si=0$, then
$$X(t)=e^{bt}\[T-\int_0^te^{-bs}ds\],\qq t\in[0,T].$$
Thus, when $b=0$, one has
$$X(t)=T-t,\qq t\in[0,T],$$
and when $b\ne0$,
$$X(t)=e^{bt}T+{1\over b}(1-e^{bt})={e^{bt}\over b}\(e^{-bt}-1+bT\),\qq t\in[0,T].$$
Since
$$e^\l-1-\l>0,\qq\forall\l\ne0,$$
we have that $b<0$ implies
$$X(T)<0.$$

The above example tells us that when $\si\ne0$, or $\si=0$ and
$b<0$, although the free term $\f(t)=T-t$ is nonnegative on $[0,T]$,
the solution $X(\cd)$ of the FSVIE (\ref{5.7}) does not necessarily
remain nonnegative on $[0,T]$. Consequently, nonnegativity of the
free term is not enough for the solution of the MF-FSVIE to be
nonnegative. Thus, besides the nonnegativity of the free term, some
additional conditions are needed.

\ms

To present positive results, we introduce the following assumption.

\ms

{\bf(C2)} The maps
$$A_0:\D^*\times\O\to\dbR^{n\times n},\q A_1:[0,T]\times\O\to\dbR^{n\times
n},\q C_0:\D^*\times\O^2\to\dbR^{n\times n},$$
are measurable and uniformly bounded. For any $t\in[0,T]$,
$s\mapsto(A_0(t,s),A_1(s))$ is $\dbF$-progressively measurable on
$[0,t]$, and $s\mapsto C_0(t,s)$ is $\dbF^2$-progressively
measurable on $[0,t]$.

\ms

We now present the following result which is simple but will be
useful later.

\ms

\bf Proposition 5.5. \sl Let {\rm(C2)} hold. Further,
\bel{}A_0(t,s,\o),C_0(t,s,\o,\o')\in\h\dbM^n_+,\q
A_1(s,\o)=0,\q\ae(t,s)\in\D^*,~\as\o,\o'\in\O.\ee
Let $X(\cd)$ be the solution to $(\ref{5.7})$, with $\f(\cd)\in
L^2_\dbF(0,T;\dbR^n)$ and $\f(\cd)\ge0$. Then
\bel{5.10}X(t)\ge\f(t)\ge0,\qq t\in[0,T].\ee

\rm

\it Proof. \rm Define
$$(\cA X)(t)=\int_0^t\(A_0(t,s)X(s)+\dbE'[C_0(t,s)X(s)]\)ds,\qq t\in[0,T].$$
By our condition, we see that
$$(\cA X)(\cd)\ge0,\qq\forall X(\cd)\in L^2_\dbF(0,T;\dbR^n),~X(\cd)\ge0.$$
Now, we define the following sequence
$$\left\{\ba{ll}
\ns\ds X_0(\cd)=\f(\cd),\\
\ns\ds X_k(\cd)=\f(\cd)+(\cA X_{k-1})(\cd),\qq k\ge1.\ea\right.$$
It is easy to see that
$$X_k(\cd)\ge\f(\cd),\qq\forall k\ge0,$$
and
$$\lim_{k\to\infty}\|X_k(\cd)-X(\cd)\|_{L^2_\dbF(0,T;\dbR^n)}=0,$$
with $X(\cd)$ being the solution to (\ref{5.7}). Then it is easy to
see that (\ref{5.10}) holds. \endpf

\ms

For the case that the diffusion is nonzero in the equation, we have
the following result.

\ms

\bf Proposition 5.6. \sl Let {\rm(C2)} hold. Suppose
\bel{}\ba{ll}
\ns\ds A_0(t,s,\o)\in\dbM_+^n,\q C_0(t,s,\o,\o')\in\h\dbM_+^n,\q
A_1(s,\o)\in\dbM_0^n,\\
\ns\ds\qq\qq\qq\qq\qq\qq\qq\ae(t,s)\in\D^*,~\as\o,\o'\in\O.\ea\ee
Moreover, let $t\mapsto(\f(t),A_0(t,s),C_0(t,s))$ be continuous, and
$\f(\cd)\in L^p_\dbF(0,T;\dbR^n)$ for some $p>2$. Further,
\bel{5.11}\ba{ll}
\ns\ds\f(t_1)\ge\f(t_0)\ge0,\q A_0(t_1,s)\h x\ge A_0(t_0,s)\h x,\q
C_0(t_1,s)\h x\ge C_0(t_0,s)\h x,\\
\ns\ds\qq\qq\qq\qq\forall\,s\le t_0<t_1\le T,~s\in[0,T],~\h
x\in\dbR^n_+,~\as\ea\ee
Let $X(\cd)$ be the solution of linear MF-FSVIE $(\ref{5.7})$. Then
\bel{5.12}X(t)\ge0,\qq t\in[0,T],~\as\ee

\ms

\it Proof. \rm Let $\Pi=\{\t_k,0\le k\le N\}$ be an arbitrary set of
finitely many $\dbF$-stopping times with $0=\t_0<\t_1<\cds<\t_N=T$,
and we define its mesh size by
$$\|\Pi\|=\esssup_{\o\in\O}\max_{1\le k\le N}|\t_k-\t_{k-1}|.$$
Let
$$\left\{\ba{ll}
\ns\ds
A^\Pi_0(t,s)=\sum_{k=0}^{N-1}A_0(\t_k,s)I_{[\t_k,\t_{k+1})}(t),\qq
C^\Pi_0(t,s)=\sum_{k=0}^{N-1}C_0(\t_k,s)I_{[\t_k,\t_{k+1})}(t),\\
\ns\ds\f^\Pi(t)=\sum_{k=0}^{N-1}
\f(\t_k)I_{[\t_k,\t_{k+1})}(t).\ea\right.$$
Clearly, each $A_0(\t_k,\cd)$ is an $\dbF$-adapted bounded process,
each $C_0(\t_k,\cd)$ is an $\dbF^2$-adapted bounded process, and
each $\f(\t_k)$ is an $\cF_{\t_k}$-measurable random variable.
Moreover, for each $k\ge0$,
\bel{}A_0(\t_k,s)\in\dbM^n_+,\q C_0(\t_k,s)\in\h\dbM^n_+,\qq
s\in[\t_k,\t_{k+1}],~~\as,\ee
and
\bel{5.13}0\le\f(\t_k)\le\f(\t_{k+1}),\qq\as\ee
Now, we let $X^\Pi(\cd)$ be the solution to the following MF-FSVIE:
\bel{}\ba{ll}
\ns\ds X^\Pi(t)=\f^\Pi(t)+\int_0^t\(A_0^\Pi(t,s)X^\Pi(s)+\dbE'\[C_0^\Pi(t,s)X^\Pi(s)\]\)ds\\
\ns\ds\qq\qq\qq\qq\qq\qq+\int_0^tA_1(s)X^\Pi(s)dW(s),\qq
t\in[0,T].\ea\ee
Then on interval $[0,\t_1)$, we have
$$
X^\Pi(t)=\f(0)+\int_0^t\(A_0(0,s)X^\Pi(s)+\dbE'\[C_0(0,s)X^\Pi(s)\]\)ds+\int_0^tA_1(s)X^\Pi(s)dW(s),$$
which is an MF-FSDE, and $X^\Pi(\cd)$ has continuous paths. From
Proposition 5.3, we have
$$X^\Pi(t)\ge0,\qq t\in[0,\t_1),~\as$$
In particular,
\bel{5.16}\ba{ll}
\ns\ds
X^\Pi(\t_1-0)=\f(0)+\int_0^{\t_1}\(A_0(0,s)X^\Pi(s)+\dbE'\[C_0(0,s)X^\Pi(s)\]\)ds\\
\ns\ds\qq\qq\qq\qq\qq+\int_0^{\t_1}A_1(s)X^\Pi(s)dW(s)\ge0.\ea\ee
Next, on $[\t_1,\t_2)$, we have (making use the above)
$$\ba{ll}
\ns\ds
X^\Pi(t)=\f(\t_1)+\int_0^{\t_1}\(A_0(\t_1,s)X^\Pi(s)+\dbE'\[C_0(\t_1,s)X^\Pi(s)\]\)ds+\int_0^{\t_1}
A_1(s)X^\Pi(s)dW(s)\\
\ns\ds\qq\qq+\int_{\t_1}^t\(A_0(\t_1,s)X^\Pi(s)+\dbE'\[C_0(\t_1,s)X^\Pi(s)\]\)ds
+\int_{\t_1}^tA_1(s)X^\Pi(s)dW(s)\\
\ns\ds\qq=\f(\t_1)-\f(0)+X^\Pi(\t_1-0)\\
\ns\ds\qq\qq+\int_0^{\t_1}\Big\{\(A_0(\t_1,s)-A_0(0,s)\)X^\Pi(s)+\dbE'\[\(C_0(\t_1,s)-C_0(0,s)\)
X^\Pi(s)\]\Big\}ds\\
\ns\ds\qq\qq+\int_{\t_1}^t\(A_0(\t_1,s)X^\Pi(s)+\dbE'\[C_0(\t_1,s)X^\Pi(s)\]\)ds
+\int_{\t_1}^tA_1(s)X^\Pi(s)dW(s)\\
\ns\ds\qq\equiv\wt
X(\t_1)+\int_{\t_1}^t\(A_0(\t_1,s)X^\Pi(s)+\dbE'\[C_0(\t_1,s)X^\Pi(s)\]\)ds
+\int_{\t_1}^tA_1(s)X^\Pi(s)dW(s),\ea$$
where, by our conditions assumed in (\ref{5.11}), and noting
(\ref{5.16}),
$$\ba{ll}
\ns\ds\wt X(\t_1)\equiv\f(\t_1)-\f(0)+X^\Pi(\t_1-0)\\
\ns\ds\qq\qq+\int_0^{\t_1}\Big\{\(A_0(\t_1,s)-A_0(0,s)\)X^\Pi(s)+\dbE'\[\(C_0(\t_1,s)-C_0(0,s)\)
X^\Pi(s)\]\Big\}ds\ge0.\ea$$
Hence, by Proposition 5.3 again, one obtains
$$X^\Pi(t)\ge0,\qq t\in[\t_1,\t_2).$$
By induction, we see that
$$X^\Pi(t)\ge0,\qq t\in[0,T],~\as$$
On the other hand, it is ready to show that
$$\lim_{\|\Pi\|\to0}\|X^\Pi(\cd)-X(\cd)\|_{L^2_\dbF(0,T;\dbR^n)}=0,$$
Then (\ref{5.12}) follows from the stability estimate in Corollary
2.7. \endpf

\ms

We now look at the following (nonlinear) MF-FSVIEs with $i=0,1$:
\bel{5.5}\ba{ll}
\ns\ds
X_i(t)=\f_i(t)+\int_0^tb_i(t,s,X_i(t),\G^b_i(t,s,X_i(s)))ds+\int_0^t\si(s,X_i(s))dW(s),\q
t\in[0,T],\ea\ee
where
\bel{}\G^b_i(t,s,X_i(s))=\int_\O\th^b_i(t,s,\o,\o',X_i(s,\o),X_i(s,\o'))\dbP(d\o').\ee
Note that $\si(\cd)$ does not contain a nonlocal term, and it is
independent of $t\in[0,T]$, as well as $i=0,1$. The following result
can be regarded as an extension of \cite{Tudor 1989} from FSVIEs to
MF-FSVIEs.

\ms

\bf Theorem 5.7. \sl For $i=0,1$, let
$b_i(\cd),\si(\cd),\th^b_i(\cd)$ appeared in $(\ref{5.5})$ satisfy
{\rm(H1)--(H2)} and $\f_i(\cd)\in L^2_\dbF(0,T;\dbR^n)$. Further,
for all $x,\bar x,x'\in\dbR^n$, $\h x\in\dbR^n_+$,
$\g\in\dbR^{m_1}$, almost all $(t,s)\in\D^*$ and almost sure
$\o,\o'\in\O$,
\bel{5.19}(b_0)_\g(t,s,\o,x,\g)\in\h\dbM^{n\times
m_1}_+,\qq\si_x(s,\o,x)\in\dbM^n_0,\ee
and maps
\bel{5.20}\ba{ll}
\ns\ds t\mapsto(b_0)_x(t,s,\o,x,\g)\h x,\\
\ns\ds t\mapsto(b_0)_\g(t,s,\o,x,\g)(\th_0^b)_x(t,s,\o,\bar
x,x')\h x,\\
\ns\ds t\mapsto(b_0)_\g(t,s,\o,x,\g)(\th_0^b)_{x'}(t,s,\o,\bar
x,x')\h x,\\
\ns\ds t\mapsto b_1(t,s,\o,x,\g)-b_0(t,s,\o,x,\g),\\
\ns\ds
t\mapsto\th^b_1(t,s,\o,\o',x,x')-\th^b_0(t,s,\o,\o',x,x'),\\
\ns\ds
t\mapsto(b_0)_\g(t,s,\o,x,\g)\[\th^b_1(t,s,\o,\o',x,x')-\th^b_0(t,s,\o,\o',x,x')\],\\
\ns\ds t\mapsto\f_1(t)-\f_0(t)\ea\ee
are continuous, nonnegative and nondecreasing on $[s,T]$. Let
$X_i(\cd)\in L^p_\dbF(0,T;\dbR^n)$ be the solutions to the
corresponding equation $(\ref{5.5})$. Then
\bel{5.21}X_0(t)\le X_1(t),\qq\forall t\in[0,T],~\as\ee

\ms

\it Proof. \rm From the equations satisfied by $X_0(\cd)$ and
$X_1(\cd)$, we have the following:
$$\ba{ll}
\ns\ds X_1(t)-X_0(t)=\f_1(t)-\f_0(t)\\
\ns\ds\qq+\int_0^t\[b_1(t,s,X_1(s),\G^b_1(t,s,X_1(s)))-b_0(t,s,X_0(s),\G^b_0(t,s,X_0(s)))\]ds\\
\ns\ds\qq+\int_0^t\[\si(s,X_1(s))-\si(s,X_0(s))\]dW(s)\\
\ns\ds\qq=\h\f_1(t)-\h\f_0(t)\\
\ns\ds\qq+\int_0^t\[b_0(t,s,X_1(s),\G^b_0(t,s,X_1(s)))-b_0(t,s,X_0(s),\G^b_0(t,s,X_0(s)))\]ds\\
\ns\ds\qq+\int_0^t\[\si(s,X_1(s))-\si(s,X_0(s))\]dW(s),\ea$$
where (making use of Proposition 5.1 and (\ref{5.19})--(\ref{5.20}))
$$\ba{ll}
\ns\ds\h\f_1(t)-\h\f_0(t)=\f_1(t)-\f_0(t)\\
\ns\ds\qq+\int_0^t\[b_1(t,s,X_1(s),\G^b_1(t,s,X_1(s)))-b_0(t,s,X_1(s),\G^b_0(t,s,X_1(s)))\]ds\\
\ns\ds=\f_1(t)-\f_0(t)+\int_0^t\[b_1(t,s,X_1(s),\G^b_1(t,s,X_1(s)))-b_0(t,s,X_1(s),\G^b_1(t,s,X_1(s)))\]ds\\
\ns\ds\qq+\int_0^t\[\int_0^1(b_0)_\g(t,s,X_1(s),\wt\G^b_\l(t,s))d\l\]\(\G^b_1(t,s,X_1(s)))
-\G^b_0(t,s,X_1(s)))\)ds\ge0,\ea$$
and nondecreasing in $t$, where
$$\wt\G^b_\l(t,s)=(1-\l)\G^b_0(t,s,X_1(s))+\l\G^b_1(t,s,X_1(s)).$$
Now, we look at the following:
$$\ba{ll}
\ns\ds b_0(t,s,X_1(s),\G^b_0(t,s,X_1(s)))-b_0(t,s,X_0(s),\G^b_0(t,s,X_0(s)))\\
\ns\ds=\[\int_0^1(b_0)_x(t,s,X_\l(s),\G^b_\l(t,s))d\l\]\(X_1(s)-X_0(s)\)\\
\ns\ds\qq+\[\int_0^1(b_0)_\g(t,s,X_\l(s),\G^b_\l(t,s))d\l\]\(\G^b_0(t,s,X_1(s))-\G^b_0(t,s,X_0(s))\)\\
\ns\ds\equiv(b_0)_x(t,s)\(X_1(s)-X_0(s)\)+(b_0)_\g(t,s)\(\G^b_0(t,s,X_1(s))-\G^b_0(t,s,X_0(s))\),\ea$$
where
\bel{Xl}\left\{\ba{ll}
\ns\ds X_\l(s)=(1-\l)X_0(s)+\l X_1(s),\\
\ns\ds\G^b_\l(t,s)=(1-\l)\G^b_0(t,s,X_0(s))+\l\G^b_0(t,s,X_1(s)).\ea\right.\ee
and
$$\left\{\ba{ll}
\ns\ds(b_0)_x(t,s)=\int_0^1(b_0)_x(t,s,X_\l(s),\G^b_\l(t,s))d\l,\\
\ns\ds(b_0)_\g(t,s)=\int_0^1(b_0)_\g(t,s,X_\l(s),\G^b_\l(t,s))d\l.\ea\right.$$
Moreover,
$$\ba{ll}
\ns\ds\G^b_0(t,s,X_1(s))-\G^b_0(t,s,X_0(s))\\
\ns\ds=\int_\O\[\th_0^b(t,s,\o,\o',X_1(s,\o),X_1(s,\o'))-\th_0^b(t,s,\o,\o',X_0(s,\o),X_0(s,\o'))\]
\dbP(d\o')\\
\ns\ds=\Big\{\int_\O\[\int_0^1(\th_0^b)_x(t,s,\o,\o',X_\l(s,\o),X_\l(s,\o'))d\l\]
\dbP(d\o')\Big\}\(X_1(s,\o)-X_0(s,\o)\)\\
\ns\ds\qq+\int_\O\[\int_0^1(\th_0^b)_{x'}(t,s,\o,\o',X_\l(s,\o),X_\l(s,\o'))d\l\]
\(X_1(s,\o')-X_0(s,\o')\)\dbP(d\o')\\
\ns\ds=\dbE'\[(\th_0^b)_x(t,s)\]\(X_1(s)-X_0(s)\)+\dbE'\[(\th^b_0)_{x'}(t,s)\(X_1(s,\o')-X_0(s,\o')\)
\],\ea$$
where
\bel{}\left\{\ba{ll}
\ns\ds(\th^b_0)_x(t,s)=\int_0^1(\th^b_0)_x(t,s,\o,\o',X_\l(s,\o),
X_\l(s,\o'))d\l,\\
\ns\ds(\th^b_0)_{x'}(t,s)=\int_0^1(\th^b_0)_{x'}(t,s,\o,\o',X_\l(s,\o),X_\l(s,\o'))d\l,\ea\right.\ee
and $X_\l(\cd)$ is defined as (\ref{Xl}). Thus,
$$\ba{ll}
\ns\ds b_0(t,s,X_1(s),\G^b_0(t,s,X_1(s)))-b_0(t,s,X_0(s),\G^b_0(t,s,X_0(s)))\\
\ns\ds=\Big\{(b_0)_x(t,s)+\dbE'\[(b_0)_\g(t,s)(\th^b_0)_x(t,s)\]\Big\}\(X_1(s)-X_0(s)\)\\
\ns\ds\qq\qq+\dbE'\[(b_0)_\g(t,s)(\th^b_0)_{x'}(t,s)\(X_1(s,\o')-X_0(s,\o')\)\]\\
\ns\ds\equiv
A_0(t,s)\(X_1(s)-X_0(s)\)+\dbE'\[C_0(t,s)\(X_1(s)-X_0(s)\)\],\ea$$
where
\bel{5.24}\left\{\ba{ll}
\ns\ds A_0(t,s)=(b_0)_x(t,s)+\dbE'\[(b_0)_\g(t,s)(\th^b_0)_x(t,s)\]\in\dbM^n_+,\\
\ns\ds
C_0(t,s)=(b_0)_\g(t,s)(\th^b_0)_{x'}(t,s)\in\h\dbM^n_+,\ea\right.\qq(t,s)\in\D^*,~\as\ee
Similarly,
$$\ba{ll}
\ns\ds\si(s,X_1(s))-\si(s,X_0(s))\equiv
A_1(s)\(X_1(s)-X_0(s)\),\ea$$
where
\bel{5.25}A_1(s)\equiv\int_0^1\si_x(s,X_\l(s))d\l\in\dbM^n_0,\qq(t,s)\in\D^*,~\as\ee
Then we have
$$\ba{ll}
\ns\ds X_1(t)-X_0(t)=\h\f_1(t)-\h\f_0(t)+\int_0^t\Big\{A_0(t,s)\(X_1(s)-X_0(s)\)\\
\ns\ds\qq\qq\qq+\dbE'\[C_0(t,s)\(X_1(s)-X_0(s)\)\]\Big\}ds+\int_0^tA_1(s)\(X_1(s)-X_0(s)\)dW(s).\ea$$
From (\ref{5.19})--(\ref{5.20}), we see that the coefficients of the
above linear MF-FSVIE satisfy (C2), and $\h\f_1(\cd)-\h\f_0(\cd)$ is
nonnegative and nondecreasing. Then (\ref{5.21}) follows from
Proposition 5.6. \endpf

\ms

From the above proof, we see that one may replace $b_0(\cd)$ in
conditions (\ref{5.19}) by $b_1(\cd)$. Also, by an approximation
argument, we may replace the derivatives in (\ref{5.19}) of
$b_0(\cd)$ and $\si(\cd)$ by the corresponding difference quotients.

\subsection{Comparison theorems for MF-BSVIEs.}

In this subsection, we discuss comparison property for MF-BSVIEs.
First, we consider the following linear MF-BSVIE:
\bel{5.26}\ba{ll}
\ns\ds Y(t)=\psi(t)+\int_t^T\(\bar A_0(t,s)Y(s)+\bar
C_0(t)Z(s,t)+\dbE'\[\bar A_1(t,s)Y(s)\]\)ds\\
\ns\ds\qq\qq\qq\qq-\int_t^TZ(t,s)dW(s),\qq t\in[0,T].\ea\ee
Note that $Z(t,s)$ does not appear in the whole drift term, and
$Z(s,t)$ does not appear in the nonlocal term. Further, the
coefficient of $Z(s,t)$ is independent of $s$. Let us introduce the
following assumption.

\ms

{\bf(C3)} The maps
$$\bar A_0:\D\times\O\to\dbR^{n\times n},\q\bar C_0:[0,T]\times\O\to\dbR^{n\times n},\q
\bar A_1:\D\times\O^2\to\dbR^{n\times n}$$
are uniformly bounded, $\bar C_0(\cd)$ is $\dbF$-progressively
measurable, and for each $t\in[0,T]$, $s\mapsto\bar A_0(t,s)$ and
$s\mapsto\bar A_1(t,s)$ are $\dbF$-progressively measurable and
$\dbF^2$-progressively measurable on $[t,T]$, respectively.

\ms

We have the following result.

\ms

\bf Theorem 5.8. \sl Let {\rm (C3)} hold. In addition, suppose
\bel{}\ba{ll}
\ns\ds\bar A_0(t,s,\o)\in\dbM_+^n,\q\bar
A_1(t,s,\o,\o')\in\h\dbM_+^n,\q
\bar C_0(s,\o)\in\dbM_0^n,\\
\ns\ds\qq\qq\qq\qq\qq\qq\qq\ae(t,s)\in\D^*,~\as\o,\o'\in\O.\ea\ee
Moreover, let $t\mapsto(\bar A_0(s,t),\bar C_0(s,t))$ be continuous,
and
\bel{5.28}\ba{ll}
\ns\ds\bar A_0(s,t_1)^Tx\ge\bar A_0(s,t_0)^Tx,\q
\bar A_1(s,t_1)^Tx\ge\bar A_1(s,t_0)^Tx,\\
\ns\ds\qq\qq\qq\qq\forall\,s\le t_0<t_1\le
T,~s\in[0,T],~x\in\dbR^n_+,~\as\ea\ee
Let $(Y(\cd),Z(\cd\,,\cd))$ be the adapted M-solution to
$(\ref{5.26})$ with $\psi(\cd)\in L^2_{\cF_T}(0,T;\dbR^n)$,
$\psi(\cd)\ge0$. Then
\bel{5.29}\dbE\[\int_t^TY(s)ds\bigm|\cF_t\]\ge0,\qq\forall
t\in[0,T],~\as\ee

\it Proof. \rm We consider the following linear MF-FSVIE:
\bel{}\ba{ll}
\ns\ds X(t)=\f(t)+\int_0^t\(\bar A_0(s,t)^TX(s)+\dbE^*[\bar
A_1(s,t)^TX(s)]\)ds\\
\ns\ds\qq\qq\qq\qq\qq+\int_0^t\(\bar C_0(s)^TX(s)\)dW(s),\qq
t\in[0,T],\ea\ee
where
$$\f(t)=\int_0^t\eta(s)ds,\qq t\in[0,T],$$
for some $\eta(\cd)\in L^2_\dbF(0,T;\dbR^n)$ with $\eta(\cd)\ge0$.
By our conditions on $\bar A_0(\cd\,,\cd)$ and $\bar
A_1(\cd\,,\cd)$, using Proposition 5.6, we have
$$X(\cd)\ge0.$$
Then by Theorem 4.2, one obtains
$$\ba{ll}
\ns\ds0\le\dbE\int_0^T\lan\psi(t),X(t)\ran
dt=\dbE\int_0^T\lan\f(t),Y(t)\ran dt\\
\ns\ds\q=\dbE\int_0^T\int_0^t\lan\eta(s),Y(t)\ran
dsdt=\dbE\int_0^T\lan\eta(s),\int_s^TY(t)dt\ran ds.\ea$$
This proves (\ref{5.29}). \endpf

\ms

Since the conditions assumed in Proposition 5.6 are very close to
necessary conditions, we feel that it is very difficult (if not
impossible) to get better comparison results for general MF-BSVIEs.
However, if the drift term does not contain $Z(\cd\,,\cd)$, we are
able to get a much better looking result. Let us now make it
precise. For $i=0,1$, we consider the following (nonlinear)
MF-BSVIEs:
\bel{5.31}\ba{ll}
\ns\ds
Y_i(t)=\psi_i(t)+\int_t^Tg_i(t,s,Y_i(s),\G_i(t,s,Y_i(s)))ds-\int_t^TZ_i(t,s)dW(s),\q
t\in[0,T],\ea\ee
where
\bel{5.26a}\ba{ll}
\ns\ds\G_i(t,s,Y_i(s))=\dbE'\[\th_i(t,s,Y_i(s),Y_i(s,\o'))\]\\
\ns\ds\qq\qq\q~~\equiv\int_\O\th_i(t,s,\o,\o',Y_i(s,\o),Y_i(s,\o'))\dbP(d\o').\ea\ee
Note that in the above, $Z_i(\cd\,,\cd)$ does not appear in the
drift term.

\ms

\bf Theorem 5.9. \sl Let
$g_i:\D\times\O\times\dbR^n\times\dbR^m\to\dbR^n$ and
$\th_i:\D\times\O^2\times\dbR^n\times\dbR^n\to\dbR^m$ satisfy
{\rm(H3)$_q$} for some $q\ge2$. Moreover, for all $y,y'\in\dbR^n$,
$\g\in\dbR^m$, almost all $(t,s)\in\D$, and almost surely
$\o,\o'\in\O$, the following hold:
\bel{}\left\{\ba{ll}
\ns\ds(g_0)_\g(t,s,\o,y,\g)\in\h\dbM^{n\times
m}_+,\qq(\th_0)_{y'}(t,s,\o,\o',y,y')\in\h\dbM_+^{m\times n},\\
\ns\ds(g_0)_y(t,s,\o,y,\g)\in\h
M_+^n,\q(\th_0)_y(t,s,\o,\o',y,y')\in\h\dbM_+^{m\times
n},\ea\right.\ee
and
\bel{}\left\{\ba{ll}
\ns\ds g_1(t,s,\o,y,\g)\ge g_0(t,s,\o,y,\g),\\
\ns\ds\th_1(t,s,\o,\o',y,y')\ge\th_0(t,s,\o,\o',y,y').\ea\right.\ee
Let $\psi_i(\cd)\in L^2_{\cF_T}(0,T;\dbR^n)$ with
\bel{}\psi_0(t)\le\psi_1(t),\qq\forall t\in[0,T],~\as,\ee
and $(Y_i(\cd),Z_i(\cd\,,\cd))$ be the adapted M-solutions to the
corresponding MF-BSVIEs $(\ref{5.31})$. Then
\bel{5.36}Y_0(t)\le Y_1(t),\qq t\in[0,T],~\as\ee

\ms

\rm

\it Proof. \rm From the MF-BSVIEs satisfied by
$(Y_i(\cd),Z_i(\cd\,,\cd))$, we have
$$\ba{ll}
\ns\ds
Y_1(t)-Y_0(t)=\psi_1(t)-\psi_0(t)+\int_t^T\[g_1(t,s,Y_1(s),\G_1(t,s,Y_1(s)))\\
\ns\ds\qq\qq\qq-g_0(t,s,Y_0(s),\G_0(t,s,Y_0(s)))\]ds
-\int_t^T\[Z_1(t,s)-Z_0(t,s)\]dW(s)\\
\ns\ds=\h\psi_1(t)-\h\psi_0(t)+\int_t^T\[g_0(t,s,Y_1(s),\G_0(t,s,Y_1(s)))
-g_0(t,s,Y_0(s),\G_0(t,s,Y_0(s)))\]ds\\
\ns\ds\qq\qq\qq-\int_t^T\[Z_1(t,s)-Z_0(t,s)\]dW(s),\ea$$
where (making use of our condition)
$$\ba{ll}
\ns\ds\h\psi_1(t)-\h\psi_0(t)=\psi_1(t)-\psi_0(t)+\int_t^T\(g_1(t,s,Y_1(s),
\G_1(t,s,Y_1(s)))-g_0(t,s,Y_1(s),\G_0(t,s,Y_1(s)))\)ds\\
\ns\ds=\psi_1(t)-\psi_0(t)+\int_t^T\(g_1(t,s,Y_1(s),\G_1(t,s,Y_1(s)))
-g_0(t,s,Y_1(s),\G_1(t,s,Y_1(s)))\)ds\\
\ns\ds\qq\qq+\int_t^T\(g_0(t,s,Y_1(s),\G_1(t,s,Y_1(s)))-g_0(t,s,Y_1(s),\G_0(t,s,Y_1(s)))\)ds\\
\ns\ds=\psi_1(t)-\psi_0(t)+\int_t^T\(g_1(t,s,Y_1(s),\G_1(t,s,Y_1(s)))
-g_0(t,s,Y_1(s),\G_1(t,s,Y_1(s)))\)ds\\
\ns\ds\qq\qq+\int_t^T(\wt
g_0)_\g(t,s)\(\G_1(t,s,Y_1(s))-\G_0(t,s,Y_1(s))\)ds\ge0,\ea$$
with
$$\ba{ll}
\ns\ds(\wt
g_0)_\g(t,s)=\int_0^1(g_0)_\g(t,s,Y_1(s),\G_\l(t,s,Y_1(s)))d\l\in\h\dbM^{n\times m}_+,\\
\ns\ds\G_\l(t,s,Y_1(s))=(1-\l)\G_0(t,s,Y_1(s))+\l\G_1(t,s,Y_1(s)).\ea$$
Next, we note that
$$\ba{ll}
\ns\ds g_0(t,s,Y_1(s),\G_0(t,s,Y_1(s)))
-g_0(t,s,Y_0(s),\G_0(t,s,Y_0(s)))\\
\ns\ds=\int_0^1\Big\{(g_0)_y(t,s,Y_\l(s),\G_\l(t,s))\[Y_1(s)-Y_0(s)\]\\
\ns\ds\qq+(g_0)_\g(t,s,Y_\l(s),\G_\l(t,s))\[\G_0(t,s,Y_1(s))
-\G_0(t,s,Y_0(s))\]\Big\}d\l\\
\ns\ds\equiv(g_0)_y(t,s)\[Y_1(s)-Y_0(s)\]+(g_0)_\g(t,s)\[\G_0(t,s,Y_1(s))
-\G_0(t,s,Y_0(s))\],\ea$$
where
$$\left\{\ba{ll}
\ns\ds Y_\l(s)=(1-\l)Y_0(s)+\l Y_1(s),\\
\ns\ds\G_\l(t,s)=(1-\l)\G_0(t,s,Y_0(s))+\l\G_0(t,s,Y_1(s)),\ea\right.$$
and
$$\left\{\ba{ll}
\ns\ds(g_0)_y(t,s)=\int_0^1(g_0)_y(t,s,Y_\l(s),\G_\l(t,s))d\l\in\dbM^n_+,\\
\ns\ds(g_0)_\g(t,s)=\int_0^1(g_0)_\g(t,s,Y_\l(s),\G_\l(t,s))d\l\in\h\dbM^{n\times
m}_+.\ea\right.$$
Also,
$$\ba{ll}
\ns\ds\G_0(t,s,Y_1(s))-\G_0(t,s,Y_0(s))\\
\ns\ds=\dbE'\[\th_0(t,s,Y_1(s),
Y_1(s,\o'))-\th_0(t,s,Y_0(s),Y_0(s,\o'))\]\\
\ns\ds=\dbE'\int_0^1\Big\{(\th_0)_y(t,s,Y_\l(s),Y_\l(s,\o'))\(Y_1(s)-Y_0(s)\)\\
\ns\ds\qq\qq+(\th_0)_{y'}(t,s,Y_\l(s),Y_\l(s,\o'))\(Y_1(s,\o')-Y_0(s,\o')\)\Big\}d\l\\
\ns\ds\equiv\dbE'\[(\th_0)_y(t,s)\]\(Y_1(s)-Y_0(s)\)+\dbE'\[(\th_0)_{y'}(t,s)\(Y_1(s,\o')-Y_0(s,\o')\)\],\ea$$
with
$$\left\{\ba{ll}
\ns\ds(\th_0)_y(t,s)=\int_0^1(\th_0)_y(t,s,Y_\l(s),Y_\l(s,\o'))d\l,\\
\ns\ds(\th_0)_{y'}(t,s)=\int_0^1(\th_0)_{y'}(t,s,Y_\l(s),Y_\l(s,\o'))d\l.\ea\right.$$
Thus,
\bel{5.30}\ba{ll}
\ns\ds Y_1(t)-Y_0(t)=\h\psi_1(t)-\h\psi_0(t)+\int_t^T\Big\{\bar
A_0(t,s)\(Y_1(s)-Y_0(s)\)\\
\ns\ds\qq+\dbE'\[\bar A_1(t,s)
\(Y_1(s)-Y_0(s)\)\]\Big\}ds-\int_t^T\(Z_1(t,s)-Z_0(s,t)\)dW(s),\q
t\in[0,T],\ea\ee
with
\bel{5.38}\left\{\ba{ll}
\ns\ds\bar
A_0(t,s)=(g_0)_y(t,s)+\dbE'\[(g_0)_\g(t,s)(\th_0)_y(t,s)\]\in\h\dbM^n_+,\\
\ns\ds\bar
A_1(t,s)=(g_0)_\g(t,s)(\th_0)_{y'}(t,s)\in\h\dbM^n_+,\ea\right.\qq(t,s)\in\D,~\as\ee
Now, for any $\f(\cd)\in L^2_\dbF(0,T;\dbR^n)$, let $X(\cd)$ be the
solution to the following linear MF-FSVIE:
\bel{}X(t)=\f(t)+\int_0^t\(\bar A_0(s,t)^TX(s)+\dbE^*\[\bar
A_1(s,t)^TX(s)\]\)ds,\qq t\in[0,T].\ee
By Proposition 5.5, we know that $X(\cd)\ge0$. Then by Theorem 4.2,
we have
$$\ba{ll}
\ns\ds0\le\dbE\int_0^T\lan\h\psi_1(t)-\h\psi_0(t),X(t)\ran
dt=\dbE\int_0^T\lan\f(t),Y_1(t)-Y_0(t)\ran dt.\ea$$
Hence, (\ref{5.36}) follows. \endpf

\ms

Combining the above two results, we are able to get a comparison
theorem for the following MF-BSVIE:
\bel{5.34}\ba{ll}
\ns\ds
Y_i(t)=\psi_i(t)+\int_t^T\(g_i(t,s,Y_i(s),\G_i(t,s,Y_i(s)))+\bar
C_0(t)Z_i(s,t)\)ds\\
\ns\ds\qq\qq-\int_t^TZ_i(t,s)dW(s),\qq t\in[0,T],\ea\ee
where $\G_i(\cd)$ is as that in (\ref{5.31}). Under proper
conditions, we will have the following comparison:
\bel{}\dbE\[\int_t^TY_0(s)ds\bigm|\cF_t\]\le\dbE\[\int_t^TY_1(s)ds\bigm|\cF_t\],\qq\forall
t\in[0,T],~\as\ee
We omit the details here.

\ms

We note that in Proposition 5.6, monotonicity conditions for
$\f(\cd)$, $A_0(\cd\,,\cd)$ and $C_0(\cd\,,\cd)$ play a crucial
role. These kind of conditions were overlooked in \cite{Yong 2006,
Yong 2007, Yong 2008}. The following example shows that in general
(\ref{5.36}) might be false.

\ms

\bf Example 5.10. \rm Consider
$$Y_0(t)=-\int_t^TY_0(s)ds,\qq t\in[0,T],$$
and
$$Y_1(t)=t-\int_t^TY_1(s)ds,\qq t\in[0,T].$$
Then
$$Y_0(t)=0,\qq t\in[0,T],$$
and the equation for $Y_1(\cd)$ is equivalent to the following:
$$\dot Y_1(t)=Y_1(t)+1,\qq Y_1(T)=T,$$
whose solution is given by
$$Y_1(t)=e^{t-T}(T+1)-1,\qq t\in[0,T].$$
It is easy to see that
$$Y_1(t)<0=Y_0(t),\qq\forall t\in[0,T-\ln(T+1)).$$
Hence, (\ref{5.36}) fails.

\ms

To conclude this section, we would like to pose the following open
question: For general MF-BSVIEs, under what conditions on the
coefficients, one has a nice-looking comparison theorem?

\ms

We hope to be able to report some results concerning the above
question before long.

\section{An Optimal Control Problem for MF-SVIEs.}

In this section, we will briefly discuss a simple optimal control
problem for MF-FSVIEs. This can be regarded as an application of
Theorem 4.1, a duality principle for MF-FSVIEs. The main clue is
similar to the relevant results presented in \cite{Yong 2006, Yong
2008}. We will omit some detailed derivations. General optimal
control problems for MF-FSVIEs will be much more involved and we
will present systematic results for that in our forthcoming
publications.

\ms

Let $U$ be a non-empty bounded convex set in $\dbR^m$, and let $\cU$
be the set of all $\dbF$-adapted processes $u:[0,T]\times\O\to U$.
Since $U$ is bounded, we see that $\cU\subseteq
L^\infty_\dbF(0,T;\dbR^m)$. For any $u(\cd)\in\cU$, consider the
following controlled MF-FSVIE:
\bel{6.1}\ba{ll}
\ns\ds
X(t)=\f(t)+\int_0^tb(t,s,X(s),u(s),\G^b(t,s,X(s),u(s)))ds\\
\ns\ds\qq\qq\qq+\int_0^t\si(t,s,X(s),u(s),\G^\si(t,s,X(s),u(s)))dW(s),\qq
t\in[0,T],\ea\ee
where
$$\left\{\ba{ll}
\ns\ds b:\D^*\times\O\times\dbR^n\times U\times\dbR^{m_1}\to\dbR^n,\\
\ns\ds \si:\D^*\times\O\times\dbR^n\times U\times\dbR^{m_2}\to
\dbR^n,\ea\right.$$
and
$$\left\{\ba{ll}
\ns\ds\G^b(t,s,X(s),u(s))=\int_\O\th^b(t,s,\o,\o',X(s,\o),u(s,\o),X(s,\o'),u(s,\o'))\dbP(d\o')\\
\ns\ds\qq\qq\qq\qq\equiv\dbE'\[\th^b(t,s,X(s),u(s),x',u')\]_{(x',u')=(X(s),u(s))},\\
\ns\ds\G^\si(t,s,X(s),u(s))=\int_\O\th^\si(t,s,\o,\o',X(s,\o),u(s,\o),X(s,\o'),u(s,\o'))\dbP(d\o')\\
\ns\ds\qq\qq\qq\qq\equiv\dbE'\[\th^\si(t,s,X(s),u(s),x',u')\]_{(x',u')=(X(s),u(s))},\ea\right.$$
with
$$\left\{\ba{ll}
\ns\ds\th^b:\D^*\times\O^2\times\dbR^n\times U\times\dbR^n\times U\to\dbR^{m_1},\\
\ns\ds\th^\si:\D^*\times\O^2\times\dbR^n\times U\times\dbR^n\times
U\to\dbR^{m_2}.\ea\right.$$
In the above, $X(\cd)$ is referred to as the {\it state process} and
$u(\cd)$ as the {\it control process}. We introduce the following
assumptions for the state equation (Comparing with (H1)--(H2)):

\ms

{\bf(H1)$''$} The maps
$$\left\{\ba{ll}
\ns\ds b:\D^*\times\O\times\dbR^n\times
U\times\dbR^{m_1}\to\dbR^n,\\
\ns\ds\si:\D^*\times\O\times\dbR^n\times
U\times\dbR^{m_2}\to\dbR^n,\ea\right.$$
are measurable, and for all
$(t,x,u,\g,\g')\in[0,T]\times\dbR^n\times
U\times\dbR^{m_1}\times\dbR^{m_2}$, the map
$$(s,\o)\mapsto(b(t,s,\o,x,u,\g),\si(t,s,\o,x,u,\g'))$$
is $\dbF$-progressively measurable on $[0,t]$. Moreover, for all
$(t,s,\o,\o')\in\D^c\times\O$, the map
$$(x,u,\g,\g')\mapsto(b(t,s,\o,x,u,\g),\si(t,s,x,u,\g'))$$
is continuously differentiable and there exists some constant $L>0$
such that
\bel{b-si-Lip2}\ba{ll}
\ns\ds
|b_x(t,s,\o,x,u,\g)|+|b_u(t,s,\o,x,u,\g)|+|b_\g(t,s,\o,x,u,\g)|\\
\ns\ds+|\si_x(t,s,\o,x,u,\g')|+|\si_u(t,s,\o,x,u,\g')|+|\si_{\g'}(t,s,\o,x,u,\g')|\le L,\\
\ns\ds\qq\qq(t,s,\o,x,u,\g,\g')\in\D^*\times\O\times \dbR^n\times
U\times\dbR^{m_1}\times\dbR^{m_2}.\ea\ee
Further,
\bel{b-si-growth2}\ba{ll}
\ns\ds|b(t,s,\o,x,u,\g)|+|\si(t,s,\o,x,u,\g')|\le
L(1+|x|+|\g|+|\g'|),\\
\ns\ds\qq\qq\qq(t,s,\o,x,u,\g,\g')\in\D^*\times\O\times\dbR^n\times
U\times\dbR^{m_1}\times\dbR^{m_2}.\ea\ee

{\bf(H2)$''$} The maps
$$\left\{\ba{ll}
\ns\ds\th^b:\D^*\times\O^2\times\dbR^n\times U\times\dbR^n\times
U\to\dbR^{m_1},\\
\ns\ds\th^\si:\D^*\times\O^2\times\dbR^n\times U\times\dbR^n\times
U\to\dbR^{m_2},\ea\right.$$
are measurable, and for all $(t,x,u,x',u')\in[0,T]\times\dbR^n\times
U\times\dbR^n\times U$, the map
$$(s,\o,\o')\mapsto(\th^b(t,s,\o,\o',x,u,x',u'),\th^\si(t,s,\o,\o',x,u,x',u'))$$
is $\dbF^2$-progressively measurable on $[0,t]$. Moreover, for any
$(t,s,\o,\o')\in\D^*\times\O^2$,
$$(x,u,\g,\g')\mapsto(\th^b(t,s,\o,\o',x,u,x',u'),\th^\si(t,s,\o,\o',x,u,x',u'))$$
is continuously differentiable and there exists some constant $L>0$
such that
\bel{th b-si-Lip2}\ba{ll}
\ns\ds
|\th^b_x(t,s,\o,\o',x,u,x',u')|+|\th^b_u(t,s,\o,\o',x,u,x',u')|\\
\ns\ds\q+|\th^b_{x'}(t,s,\o,\o',x,u,x',u')|+|\th^b_{u'}(t,s,\o,\o',x,u,x',u')|\\
\ns\ds
\q+|\th^\si_x(t,s,\o,\o',x,u,x',u')|+|\th^\si_u(t,s,\o,\o',x,u,x',u')|\\
\ns\ds\q+|\th^\si_{x'}(t,s,\o,\o',x,u,x',u')|+|\th^\si_{u'}(t,s,\o,\o',x,u,x',u')|\le L,\\
\ns\ds\qq\qq(t,s,\o,\o',x,u,x',u')\in\D^*\times\O^2\times\dbR^n\times
U\times\dbR^n\times U.\ea\ee
Further,
\bel{th b-si-growth2}\ba{ll}
\ns\ds|\th^b(t,s,\o,\o',x,u,x',u')|+|\th^\si(t,s,\o,\o',x,u,x',u')|\le
L(1+|x|+|x'|),\\
\ns\ds\qq\qq\qq(t,s,\o,\o',x,u,x',u')\in\D^*\times\O^2\times\dbR^n\times
U\times\dbR^n\times U.\ea\ee
It is easy to see that under (H1)$''$--(H2)$''$, for any given
$u(\cd)\in\cU$, the state equation (\ref{6.1}) satisfies (H1)--(H2).
Hence, for any $\f(\cd)\in L^p_\dbF(0,T;\dbR^n)$, (\ref{6.1}) admits
a unique solution $X(\cd)\in L^p(0,T;\dbR^n)$.

\ms

To measure the performance of the control process $u(\cd)$, the
following (Lagrange type) {\it cost functional} is defined:
\bel{6.2}J(u(\cd))=\dbE\int_0^Tg(s,X(s),u(s),\G^g(s,X(s),u(s)))ds,\ee
where
$$g:[0,T]\times\O\times\dbR^n\times
U\times\dbR^\ell\to\dbR,$$
and
$$\ba{ll}
\ns\ds\G^g(s,X(s),u(s))=\int_\O\th^g(s,\o,\o',X(s,\o),u(s,\o),X(s,\o'),u(s,\o'))\dbP(d\o')\\
\ns\ds\qq\qq\qq\q~\equiv\dbE'\[\th^g(s,X(s),u(s),x',u')\]_{(x',u')=(X(s),u(s))},\ea$$
with
$$\th^g:[0,T]\times\O^2\times\dbR^n\times U\times\dbR^n\times
U\to\dbR^\ell.$$
For convenience, we make the following assumptions for the functions
involved in the cost functional.

\ms

{\bf(H1)$'''$} The map $g:[0,T]\times\O\times\dbR^n\times
U\times\dbR^\ell\to\dbR$ is measurable, and for all
$(x,u,\g)\in\dbR^n\times U\times\dbR^\ell$, the map $(t,\o)\mapsto
g(t,\o,x,u,\g)$ is $\dbF$-progressively measurable. Moreover, for
almost all $(t,\o)\in\D^*\times\O$, the map $(x,u,\g)\mapsto
g(t,\o,x,u,\g)$ is continuously differentiable and there exists some
constant $L>0$ such that
\bel{g-Lip2}\ba{ll}
\ns\ds
|g_x(t,\o,x,u,\g)|+|g_u(t,\o,x,u,\g)|+|g_\g(t,\o,x,u,\g)|\le L,\\
\ns\ds\qq\qq(t,\o,x,u,\g)\in[0,T]\times\O\times\dbR^n\times
U\times\dbR^\ell.\ea\ee
Further,
\bel{g-growth2}\ba{ll}
\ns\ds|g(t,\o,x,u,\g)|\le L(1+|x|+|\g|),\\
\ns\ds\qq\qq\qq(t,\o,x,u,\g)\in[0,T]\times\O\times\dbR^n\times
U\times\dbR^\ell.\ea\ee

{\bf(H2)$'''$} The map $\th^g:[0,T]\times\O^2\times\dbR^n\times
U\times\dbR^n\times U\to\dbR^\ell$ is measurable, and for all
$(x,u,x',u')\in\dbR^n\times U\times\dbR^n\times U$, the map
$(s,\o,\o')\mapsto(\th^g(t,\o,\o',x,u,x',u')$ is
$\dbF^2$-progressively measurable. Moreover, for almost all
$(t,\o,\o')\in[0,T]\times\O^2$, the map
$(x,u,x',u')\mapsto\th^g(t,s,\o,\o',x,u,x',u')$ is continuously
differentiable and there exists some constant $L>0$ such that
\bel{th g-si-Lip2}\ba{ll}
\ns\ds
|\th^g_x(t,\o,\o',x,u,x',u')|+|\th^g_u(t,\o,\o',x,u,x',u')|\\
\ns\ds\q+|\th^g_{x'}(t,\o,\o',x,u,x',u')|+|\th^g_{u'}(t,\o,\o',x,u,x',u')|\le L,\\
\ns\ds\qq\qq(t,\o,\o',x,u,x',u')\in[0,T]\times\O^2\times\dbR^n\times
U\times\dbR^n\times U.\ea\ee
Further,
\bel{th g-growth2}\ba{ll}
\ns\ds|\th^g(t,\o,\o',x,u,x',u')|\le
L(1+|x|+|x'|),\\
\ns\ds\qq\qq\qq(t,\o,\o',x,u,x',u')\in[0,T]\times\O^2\times\dbR^n\times
U\times\dbR^n\times U.\ea\ee

Under (H1)$''$--(H2)$''$ and (H1)$'''$--(H2)$'''$, the cost
functional $J(u(\cd))$ is well-defined. Then we can state our
optimal control problem as follows.

\ms

\bf Problem (C). \rm For given $\f(\cd)\in L^p_\dbF(0,T;\dbR^n)$,
find $\bar u(\cd)\in\cU$ such that
\bel{6.11}J(\bar u(\cd))=\inf_{u(\cd)\in\cU}J(u(\cd)).\ee

\ms

Any $\bar u(\cd)\in\cU$ satisfying (\ref{6.11}) is called an {\it
optimal control} of Problem (C), and the corresponding state process
$\bar X(\cd)$ is called an {\it optimal state process}. In this
case, we refer to $(\bar X(\cd),\bar u(\cd))$ as an {\it optimal
pair}.

\ms

We now briefly derive the Pontryagin type maximum principle for any
optimal pair $(\bar X(\cd),\bar u(\cd))$. To this end, we take any
$u(\cd)\in\cU$, let
$$u^\e(\cd)=\bar u(\cd)+\e[u(\cd)-\bar u(\cd)]=(1-\e)\bar u(\cd)+\e
u(\cd)\in\cU.$$
Let $X^\e(\cd)$ be the corresponding state process. Then
$$X_1(\cd)=\lim_{\e\to0}{X^\e(\cd)-\bar X(\cd)\over\e}$$
satisfies the following:
$$\ba{ll}
\ns\ds X_1(t)=\int_0^t\Big\{b_x(t,s)X_1(s)+b_u(t,s)[u(s)-\bar
u(s)]\\
\ns\ds\qq\qq\qq+b_\g(t,s)\dbE'\[\th^b_x(t,s)X_1(s,\o)+\th^b_u(t,s)[u(s,\o)-\bar
u(s,\o)]\\
\ns\ds\qq\qq\qq+\th^b_{x'}(t,s)X_1(s,\o')+\th^b_{u'}(t,s)[u(s,\o')-\bar
u(s,\o')]\]\Big\}ds\\
\ns\ds\qq\qq+\int_0^t\Big\{\si_x(t,s)X_1(s)+\si_u(t,s)[u(s)-\bar
u(s)]\\
\ns\ds\qq\qq\qq+\si_\g(t,s)\dbE'\[\th^\si_x(t,s)X_1(s,\o)+\th^\si_u(t,s)[u(s,\o)-\bar
u(s,\o)]\\
\ns\ds\qq\qq\qq+\th^\si_{x'}(t,s)X_1(s,\o')+\th^\si_{u'}(t,s)[u(s,\o')-\bar
u(s,\o')]\]\Big\}dW(s)\\
\ns\ds\qq=\int_0^t\Big\{\[b_x(t,s)+b_\g(t,s)\dbE'\th^b_x(t,s)\]X_1(s)\\
\ns\ds\qq\qq\qq+\[b_u(t,s)+b_\g(t,s)\dbE'\th^b_u(t,s)\][u(s)-\bar
u(s)]\\
\ns\ds\qq\qq\qq+\dbE'\[b_\g(t,s)\th^b_{x'}(t,s)X_1(s)+b_\g(t,s)\th^b_{u'}(t,s)[u(s)-\bar
u(s)]\]\Big\}ds\\
\ns\ds\qq\qq+\int_0^t\Big\{\[\si_x(t,s)+\si_\g(t,s)\dbE'\th^\si_x(t,s)\]X_1(s)\\
\ns\ds\qq\qq\qq+\[\si_u(t,s)+\si_\g(t,s)\dbE'\th^\si_u(t,s)\][u(s)-\bar
u(s)]\\
\ns\ds\qq\qq\qq+\dbE'\[\si_\g(t,s)\th^\si_{x'}(t,s)X_1(s)+\si_\g(t,s)\th^\si_{u'}(t,s)[u(s)-\bar
u(s)]\]\Big\}dW(s)\\
\ns\ds\q\equiv\int_0^t\Big\{A_0(t,s)X_1(s)+B_0(t,s)[u(s)-\bar
u(s)]+\dbE'\[C_0(t,s)X_1(s)+D_0(t,s)[u(s)-\bar
u(s)]\]\Big\}ds\\
\ns\ds\q\q+\2n\int_0^t\2n\Big\{A_1(t,s)X_1(s)\1n+\1n
B_1(t,s)[u(s)\1n-\1n\bar u(s)]\1n+\1n\dbE'\[C_1(t,s)X_1(s)\1n+\1n
D_1(t,s)[u(s)\1n-\1n\bar
u(s)]\]\Big\}dW(s)\\
\ns\ds\q\equiv\h\f(t)+\int_0^t\Big\{A_0(t,s)X_1(s)+\dbE'\[C_0(t,s)X_1(s)\]\Big\}ds\\
\ns\ds\qq\qq\qq+\int_0^t\Big\{A_1(t,s)X_1(s)+\dbE'\[C_1(t,s)X_1(s)\]\Big\}dW(s),\ea$$
where
$$\left\{\ba{ll}
\ns\ds b_\xi(t,s)=b_\xi(t,s,\bar X(s),\bar u(s),\G^b(t,s,\bar
X(s),\bar
u(s))),\qq\xi=x,u,\g,\\
\ns\ds\th^b_\xi(t,s)=\th^b_\xi(t,s,\o,\o',\bar X(s,\o),\bar
u(s,\o),\bar
X(s,\o'),\bar u(s,\o')),\qq\xi=x,u,x',u',\\
\ns\ds\si_\xi(t,s)=\si_\xi(t,s,\bar X(s),\bar u(s),\G^\si(t,s,\bar
X(s),\bar
u(s))),\qq\xi=x,u,\g,\\
\ns\ds\th^\si_\xi(t,s)=\th^b_\xi(t,s,\o,\o',\bar X(s,\o),\bar
u(s,\o),\bar X(s,\o'),\bar u(s,\o')),\qq\xi=x,u,x',u',\ea\right.$$
and
$$\left\{\ba{ll}
\ns\ds A_0(t,s)=b_x(t,s)+b_\g(t,s)\dbE'\th^b_x(t,s),\qq
B_0(t,s)=b_u(t,s)+b_\g(t,s)\dbE'\th^b_u(t,s),\\
\ns\ds C_0(t,s)=b_\g(t,s)\th^b_{x'}(t,s),\qq
D_0(t,s)=b_\g(t,s)\th^b_{u'}(t,s),\\
\ns\ds A_1(t,s)=\si_x(t,s)+\si_\g(t,s)\dbE'\th^\si_x(t,s),\qq
B_1(t,s)=\si_u(t,s)+\si_\g(t,s)\dbE'\th^\si_u(t,s),\\
\ns\ds C_1(t,s)=\si_\g(t,s)\th^\si_{x'}(t,s),\qq
D_1(t,s)=\si_\g(t,s)\th^\si_{u'}(t,s).\ea\right.$$
Also,
$$\ba{ll}
\ns\ds\h\f(t)=\int_0^t\Big\{B_0(t,s)[u(s)-\bar
u(s)]+\dbE'\[D_0(t,s)[u(s)-\bar
u(s)]\]\Big\}ds\\
\ns\ds\qq\qq+\int_0^t\Big\{B_1(t,s)[u(s)-\bar
u(s)]+\dbE'\[D_1(t,s)[u(s)-\bar u(s)]\]\Big\}dW(s).\ea$$
On the other hand, by the optimality of $(\bar X(\cd),\bar u(\cd))$,
we have
$$\ba{ll}
\ns\ds0\le\lim_{\e\to0}{J(u^\e(\cd))-J(\bar u(\cd))\over\e}\\
\ns\ds\q=\dbE\int_0^T\Big\{g_x(s)X_1(s)+g_u(s)[u(s)-\bar u(s)]\\
\ns\ds\qq\qq+g_\g(s)\dbE'\[\th^g_x(s)X_1(s,\o)+\th^g_u(s)[u(s,\o)-\bar
u(s,\o)]\\
\ns\ds\qq\qq+\th^g_{x'}(s)X_1(s,\o')+\th^g_{u'}(s)[u(s,\o')-\bar
u(s,\o')]\]\Big\}ds\\
\ns\ds\q=\dbE\int_0^T\Big\{\[g_x(s)+g_\g(s)\dbE'\th^g_x(s)\]X_1(s)+\[g_u(s)+g_\g(s)\dbE'\th^g_u(s)\]
[u(s)-\bar u(s)]\\
\ns\ds\qq\qq+\dbE'\[g_\g(s)\th^g_{x'}(s)X_1(s)+g_\g(s)\th^g_{u'}(s)[u(s)-\bar
u(s)]\]\Big\}ds\\
\ns\ds\q=\dbE\int_0^T\Big\{a_0(s)^TX_1(s)+b_0(s)^T[u(s)-\bar u(s)]\\
\ns\ds\qq\qq+\dbE'\[c_0(s)^TX_1(s)+d_0(s)^T[u(s)-\bar
u(s)]\]\Big\}ds\\
\ns\ds\q=\dbE\Big\{\h\f_0+\int_0^T\(a_0(s)^TX_1(s)+\dbE'\[c_0(s)^TX_1(s)\]\)ds\Big\},\ea$$
where
$$\left\{\ba{ll}
\ns\ds g_\xi(s)=g_\xi(s,\bar X(s),\bar u(s),\G^g(s,\bar X(s),\bar
u(s))),\qq\xi=x,u,\g,\\
\ns\ds\th^g_\xi(s)=\th^g_\xi(s,\o,\o',\bar X(s,\o),\bar u(s,\o),\bar
X(s,\o'),\bar u(s,\o')),\qq\xi=x,u,x',u',\ea\right.$$
and
$$\left\{\ba{ll}
\ns\ds a_0(s)^T=g_x(s)+g_\g(s)\dbE'\th^g_x(s),\qq
b_0(s)^T=g_u(s)+g_\g(s)\dbE'\th^g_u(s),\\
\ns\ds c_0(s)^T=g_\g(s)\th^g_{x'}(s),\qq
d_0(s)^T=g_\g(s)\th^g_{u'}(s),\\
\ns\ds\h\f_0=\int_0^T\Big\{b_0(s)^T[u(s)-\bar
u(s)]+\dbE'\[d_0(s)^T[u(s)-\bar u(s)]\]\Big\}ds.\ea\right.$$
Then for any undetermined $(Y(\cd),Z(\cd\,,\cd))\in\cM^2[0,T]$,
similar to the proof of Theorem 4.1, we have
$$\ba{ll}
\ns\ds\dbE\int_0^T\lan Y(t),\h\f(t)\ran dt=\dbE\int_0^T\lan
X_1(t),Y(t)-\int_t^T\(A_0(s,t)^TY(s)+A_1(s,t)^TZ(s,t)\\
\ns\ds\qq\qq\qq\qq\qq\qq+\dbE^*\[C_0(s,t)^TY(s)+C_1(s,t)^TZ(s,t)\]\)ds\ran
dt.\ea$$
Hence,
$$\ba{ll}
\ns\ds0\le\dbE\Big\{\h\f_0+\int_0^T\(a_0(s)^TX_1(s)+\dbE'\[c_0(s)^TX_1(s)\]\)ds\Big\}\\
\ns\ds\q=\dbE\Big\{\h\f_0-\int_0^T\lan Y(t),\h\f(t)\ran
dt+\int_0^T\lan
X_1(t),Y(t)-\int_t^T\(A_0(s,t)^TY(s)\\
\ns\ds\qq\qq+A_1(s,t)^TZ(s,t)+\dbE^*\[C_0(s,t)^TY(s)+C_1(s,t)^TZ(s,t)\]\)ds\ran
dt\\
\ns\ds\qq\qq+\int_0^T\(\lan X_1(t),a_0(t)\ran+\dbE'\[\lan
X_1(t),c_0(t)\ran\]\)dt\Big\}\\
\ns\ds\q=\dbE\Big\{\h\f_0-\int_0^T\lan Y(t),\h\f(t)\ran
dt+\int_0^T\lan
X_1(t),Y(t)+a_0(t)+\dbE^*c_0(t)\\
\ns\ds\qq-\int_t^T\(A_0(s,t)^TY(s)+A_1(s,t)^TZ(s,t)\\
\ns\ds\qq\qq+\dbE^*\[C_0(s,t)^TY(s)+C_1(s,t)^TZ(s,t)\]\)ds\ran
dt\Big\}.\ea$$
We now let $(Y(\cd),Z(\cd\,,\cd))\in\cM^2[0,T]$ be the adapted
M-solution to the following MF-BSVIE:
\bel{adjoint}\ba{ll}
\ns\ds Y(t)=-a_0(t)-\dbE^*c_0(t)+\int_t^T\(A_0(s,t)^TY(s)+A_1(s,t)^TZ(s,t)\\
\ns\ds\qq\qq+\dbE^*\[C_0(s,t)^TY(s)+C_1(s,t)^TZ(s,t)\]\)ds
dt-\int_t^TZ(t,s)dW(s).\ea\ee
Then
$$\ba{ll}
\ns\ds0\le\dbE\Big\{\h\f_0-\int_0^T\lan Y(t),\h\f(t)\ran dt\Big\}\\
\ns\ds\q=\dbE\Big\{\int_0^T\Big\{\lan b_0(t),u(t)-\bar u(t)\ran+\dbE'\[\lan d_0(t),u(t)-\bar u(t)\ran\]\Big\}dt\\
\ns\ds\qq-\int_0^T\lan Y(t),\int_0^t\(B_0(t,s)[u(s)-\bar
u(s)]+\dbE'\[D_0(t,s)[u(s)-\bar
u(s)]\]\)ds\\
\ns\ds\qq\qq+\int_0^t\(B_1(t,s)[u(s)-\bar
u(s)]+\dbE'\[D_1(t,s)[u(s)-\bar u(s)]\]\)dW(s)\ran dt\Big\}\\
\ns\ds\q=\dbE\Big\{\int_0^T\(\lan b_0(t)+[\dbE^*d_0(t)],u(t)-\bar u(t)\ran\)dt\\
\ns\ds\qq-\int_0^T\lan\int_t^T\(B_0(s,t)^TY(s)+\dbE^*[D_0(s,t)^TY(s)]\)ds,u(t)-\bar
u(t)\ran dt\\
\ns\ds\qq-\int_0^T\lan\int_t^T\(B_1(s,t)^TZ(s,t)+\dbE^*[D_1(s,t)^TZ(s,t)]\)ds,u(t)-\bar
u(t)\ran dt\Big\}.\ea$$
Hence, we must have the following variational inequality:
\bel{variational}\ba{ll}
\ns\ds\lan
b_0(t)+[\dbE^*d_0(t)]-\int_t^T\(B_0(s,t)^TY(s)+\dbE^*[D_0(s,t)^TY(s)]\\
\ns\ds\qq\qq\qq\qq+B_1(s,t)^TZ(s,t)+\dbE^*[D_1(s,t)^TZ(s,t)]\)ds,u-\bar
u(t)\ran\ge0,\\
\ns\ds\qq\qq\qq\qq\qq\qq\qq\qq\qq\qq\qq\forall u\in U,~\ae
t\in[0,T],~\as\ea\ee
We now summarize the above derivation.

\ms

\bf Theorem 5.1. \sl Let {\rm(H1)$''$--(H2)$''$} and
{\rm(H1)$'''$--(H2)$'''$} hold and let $(\bar X(\cd),\bar u(\cd))$
be an optimal pair of Problem {\rm(C)}. Then the adjoint equation
$(\ref{adjoint})$ admits a unique adapted M-solution
$(Y(\cd),Z(\cd\,,\cd))\in\cM^2[0,T]$ such that the variational
inequality $(\ref{variational})$ holds.

\ms

\rm

The purpose of presenting a simple optimal control problem of
MF-FSVIEs here is to realize a major motivation of studying
MF-BSVIEs. It is possible to discuss Bolza type cost functional.
Also, some of the assumptions assumed in this section might be
relaxed. However, we have no intention to have a full exploration of
general optimal control problems for MF-FSVIEs in the current paper
since such kind of general problems (even for FSVIEs) are much more
involved and they deserve to be addressed in another paper. We will
report further results along that line in a forthcoming paper.

\end{document}